%% file: ex_article_May_27_21.tex
\documentclass[onefignum,onetabnum]{siamart190516}

\DeclareMathOperator*{\argmin}{arg\,min}
\newcommand{\bx}{\mathbf{x}}

\input{ex_shared}

\ifpdf
\hypersetup{
  pdftitle={ Identification of Differential Equations by Numerical Techniques 
  from A Single Set of Noisy Observation},
  pdfauthor={Y. He, S.H. Kang, W. Liao, H. Liu, Y. Liu}
}
\fi

\begin{document}

\maketitle

\begin{abstract}
  We propose robust methods to identify the underlying Partial Differential Equation (PDE) from a given single  set of noisy time-dependent data.   We assume that the governing equation of the PDE is a linear combination of a few linear and nonlinear differential terms in a prescribed dictionary.  Noisy data make such identification particularly challenging.  Our objective is to develop robust methods against a high level of noise and approximate the underlying noise-free dynamics well.
  We first introduce a Successively Denoised Differentiation (SDD)  scheme to stabilize the amplified noise in numerical differentiation.  SDD effectively denoises the given data and the corresponding derivatives.  Secondly, we present two algorithms for PDE identification: Subspace pursuit Time evolution error (ST) and Subspace pursuit Cross-validation (SC).     Our general strategy is first to find a candidate set using the Subspace Pursuit (SP) greedy algorithm, then choose the best one via time evolution or cross-validation.  ST uses a multi-shooting numerical time evolution and selects the PDE, which yields the least evolution error.  SC evaluates the cross-validation error in the least-squares fitting and picks the PDE that gives the smallest validation error.   
  We present various numerical experiments to validate our methods. Both methods are efficient and robust to noise.
\end{abstract}

\begin{keywords}
  inverse problem, PDE identification, noisy data
\end{keywords}

\begin{AMS}
  35R30, 65Z05, 65M32
\end{AMS}

\section{Introduction}\label{sec:intro}

 Partial Differential Equations (PDEs) are used to model various real-world phenomena in science and engineering.  Numerical solvers for PDEs and analysis of various properties of the solutions have been widely studied in the literature.
In this paper, we focus on the inverse problem: Given a single set of time-dependent noisy data, how to identify the underlying PDE?

Let the given noisy time-dependent discrete data set be
\begin{equation}
\mathbf{D}:=\{U_\mathbf{i}^n\in\mathbb{R}\mid n = 0,\cdots, N;
\mathbf{i} = (i_1,\cdots,i_d) \text{ with } i_j=0,\cdots, M-1,j=1,\cdots,d \}
\label{eqdata}
\end{equation}
for sufficiently large integers $N,M\in\mathbb{N}$, where $\mathbf{i}$ is a $d$-dimensional spatial index of a discretized domain in $\mathbb{R}^d$, and $n$ represents the time index at time $t^n$.   The objective is to find an evolutionary PDE of the form
\begin{align}\label{originalPDE}
\partial_t u = f(u,\partial_{\mathbf{x}}u,\partial_{\mathbf{x}}^2u,\cdots,\partial_{\mathbf{x}}^ku,\cdots)\;,
\end{align}
which represents the dynamics of the given data $\mathbf{D}$.   Here $t$ is the time variable, $\mathbf{x}=[x_1,...,x_d]\in \mathbb{R}^d$ denotes the space variable, and
$\partial_{\mathbf{x}}^ku$ denotes the set of partial derivatives of $u$ with respect to the space variable of order $k$ for $k=0,1,\cdots$, i.e., $\partial_{\mathbf{x}}^ku:=\left\{\frac{\partial^{k}u}{\partial x_1^{k_1}\partial x_2^{k_2}\cdots\partial x_d^{k_d}}\mid k_1,\cdots,k_d\in\mathbb{N},\ \sum_{j=1}^dk_j=k\right\}.$   We assume that $f$ is a polynomial of its arguments so that the right-hand side of \eqref{originalPDE} is a linear combination of linear and nonlinear differential terms. The model in \eqref{originalPDE} includes a class of parametric PDEs where the parameters are the polynomial coefficients in $f$.

Parameter identification in differential equations and dynamical systems has been considered by physicists or applied scientists.
Earlier works include \cite{bock1983recent,muller2004parameter,baake1992fitting,muller2002fitting,bock1981numerical,parlitz2000prediction,bar1999fitting}, and among which,  \cite{muller2004parameter,bar1999fitting} considered the PDE model as in \eqref{originalPDE}.
Two important papers \cite{bongard2007automated,schmidt2009distilling} used symbolic regression to recover the underlying physical systems from experimental data.
Recently, sparse regression and $L_1$-minimization were introduced to promote sparsity in the identification of PDEs or dynamical systems \cite{brunton2016discovering,schaeffer2017learning,rudy2017data,kang2019ident}. In \cite{brunton2016discovering},  Brunton et al. considered the discovery of nonlinear dynamical systems with sparsity-promoting techniques.  The underlying dynamical systems are assumed to be governed by a small number of active terms in a prescribed dictionary, and sparse regression is used to identify these active terms.  The extensions of this sparse regression approach can be found in \cite{kaiser2018sparse,loiseau2018constrained,mangan2017model}.  In \cite{schaeffer2017learning}, Schaeffer considered the problem of PDE identification using the spectral method and focused on the benefit of using $L_1$-minimization for sparse coefficient recovery. The identification of dynamical systems with highly corrupted and undersampled data are considered in \cite{tran2017exact,schaeffer2018extracting}.
In \cite{rudy2017data}, Rudy et al.~proposed to identify PDEs by solving the $L_0$-regularized regression problem followed by a post-processing step of thresholding. 
Sparse Bayesian regression was considered in \cite{zhang2018robust} for the recovery of dynamical systems.
This series of work focused on the benefit of using $L_1$-minimization to resolve dynamical systems or PDEs with specific sparse pattern \cite{schaeffer2013sparse}. In Appendix \ref{Details_app},  we compare some existing methods in terms of the objectives in minimization. Recent works such as~\cite{messenger2020weak} and \cite{gurevich2019robust} introduced PDE learning in a weak formulation to ameliorate the errors due to the instability of numerical differentiation, when the given data are contaminated by noise. This weak formulation gives rise to a robust recovery, while it requires the underlying PDE to possess a weak formulation such that all partial derivatives in the PDE can be transferred to a test function through integration by parts.    Another related problem is to infer the interaction law in a system of agents from the trajectory data.  In \cite{bongini2017inferring,lu2018nonparametric}, nonparametric regression was used to predict the interaction function, and a theoretical guarantee was established.
Another category of methods uses deep learning \cite{long2017pde,long2019pde, raissi2017physics,qin2018data, raissi2018hidden,khoo2018switchnet,lusch2018deep}.

The most closely related work to this paper is \cite{kang2019ident}, where Identifying Differential Equation with Numerical Time evolution (IDENT) was proposed, also for a single set of given data.  It is based on the convergence principle of numerical PDE schemes.  LASSO is used to find a candidate set efficiently, and the correct PDE is identified by computing the numerical Time Evolution Error (TEE).  Among all the PDEs from the candidate set, the one whose numerical solution best matches the given data dynamics is chosen as the identified PDE.
When the given data are contaminated by noise,  the authors used a Least-Square Moving Average method to denoise the data as a pre-processing step. When the coefficients vary in the spatial domain, a Base Element Expansion (BEE) technique was proposed to recover the varying coefficients.

\begin{figure}
	\centering
	\begin{tabular}{ccc}
		(a)&(b)&(c)\\
		\includegraphics[width=0.3\textwidth]{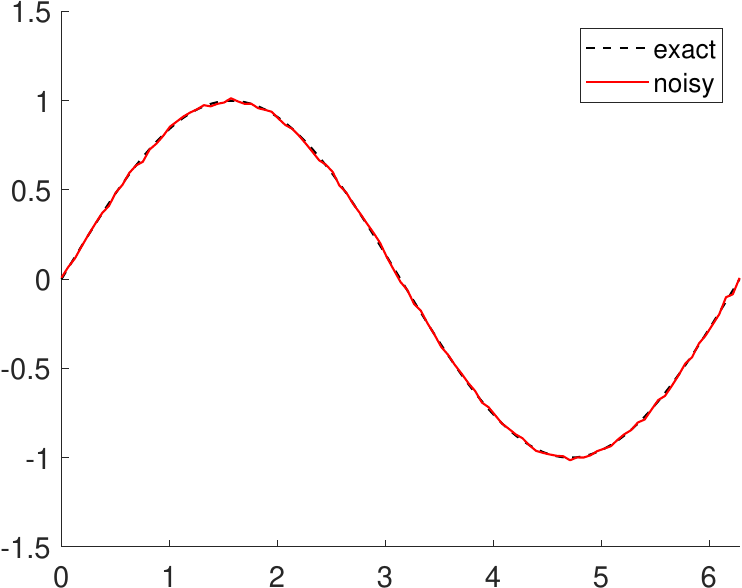}&
		\includegraphics[width=0.3\textwidth]{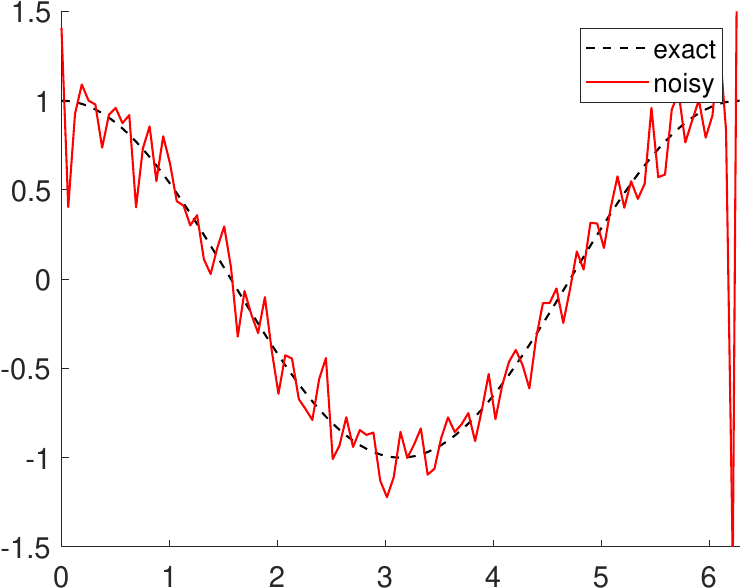}&
		\includegraphics[width=0.3\textwidth]{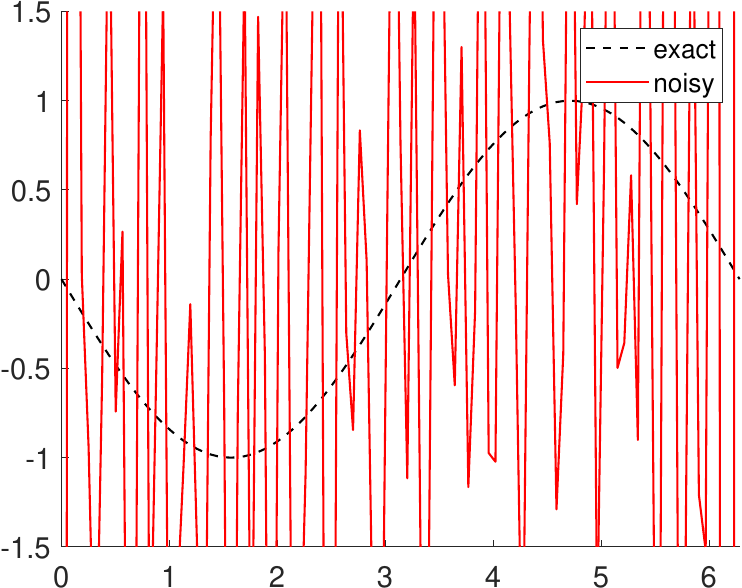}
	\end{tabular}
	\caption{The sensitivity of numerical differentiation to noise. (a) Graph of $\sin( x)$, $0\leq x\leq 2\pi$ (black), and its noisy version (red) with Gaussian noise of mean $0$ and standard deviation $0.01$. (b) The first-order derivatives of the function (black) and the data (red). (c) The second-order derivatives of the function (black) and the data (red).  The derivatives of data in (b) and (c) are computed using the five-point ENO scheme.  As the order of derivative increases, the noise gets amplified.}
	\label{fig.amplif}
\end{figure}

Despite the developments of many useful methods, when the given data are noisy, PDE identification is still challenging.    A small amount of noise can make a recovery unstable, especially for high order PDEs.  It was shown in  \cite{kang2019ident} that the noise to signal ratio for LASSO depends on the order of the underlying PDE, and IDENT can handle a small amount of noise when the PDE contains high order derivatives.
A significant issue is that the numerical differentiation often magnifies noise, which is illustrated by an example in Figure~\ref{fig.amplif}.

In this paper, we propose two robust PDE identification methods that can handle a large amount of noise given a single set of time-dependent data.   Our contributions include:
\begin{enumerate}
	\item{First, we propose a new denoising procedure, called Successively Denoised Differentiation (SDD), to stabilize the numerical differentiation applied to noisy data.  
	}
	\vspace{-0.1cm}
	\item{Second, we present two recovery algorithms which are robust against noise: Subspace pursuit Time evolution (ST) and Subspace pursuit Cross-validation (SC).  Both methods utilize the Subspace Pursuit (SP) greedy algorithm~\cite{dai2008subspace} for selecting a candidate set.  ST considers a multi-shooting numerical time evolution error, and  SC evaluates the cross-validation error in the least-squares fitting.   Both methods are efficient and robust against noise.  
	}
\end{enumerate}

This paper is organized as follows. In Section \ref{sec:data}, we introduce the PDE identification problem and describe the SDD scheme. Our proposed ST and SC algorithms are presented in Section \ref{sec:methods},  and systematic numerical experiments are provided in Section \ref{sec::numexp}.  We conclude the paper in Section \ref{sec::conclusion}, and some details are discussed in the Appendix.

\section{Data Organization and Denoising} \label{sec:data}

\subsection{Data Organization and Notations}\label{sec_data}

Let the time-space domain be $\Omega=[0,T]\times[0,X]^d$ for some $T>0$ and $X>0$.  Suppose the noisy data $\mathbf{D}$ are given as \eqref{eqdata}
on a regular grid in $\Omega$, with time index $n=0,\cdots, N$, $N\in\mathbb{N}$ and spatial index $\mathbf{i}\in\mathbb{I}$, where $\mathbb{I}=\{(i_1,\cdots,i_d)\mid i_j=0,\cdots, M-1,\ j=1,\cdots,d, M\in\mathbb{N}\}$.  Denote $\Delta t :=T/N$ and $\Delta x:=X/(M-1)$ as the time and space spacing in the given data, respectively.

At the time $t^n$ and the location $x_{\mathbf{i}}$, the datum is given as
\begin{equation}
U_\mathbf{i}^n=u(\bx_\mathbf{i},t^n)+\varepsilon_\mathbf{i}^n\;,
\label{Udef}
\end{equation}
where $t^n:=n\Delta t \in[0,T]$,   $\bx_\mathbf{i}:=(i_{1}\Delta x,\cdots,i_{d}\Delta x)\in[0,X]^d$, and $\varepsilon_\mathbf{i}^n$ is i.i.d. random noise with mean $0$.
For $n=0,1,\cdots, N-1$,  we vectorize the data in all spatial domains at time $t_n$, and denote it as $U^n\in\mathbb{R}^{M^d}$. Concatenating the vectors $\{U^n\}_{n=0}^{N-1}$ vertically gives rise to a long vector $U\in\mathbb{R}^{NM^d}$.

The underlying function $f$ in (\ref{originalPDE}) is assumed to be a finite order polynomial of its arguments:
\begin{align}
f(u,\partial_{\mathbf{x}}u,\partial_{\mathbf{x}}^2u,\cdots\partial_{\mathbf{x}}^ku,\cdots)&=c_1+c_{2}\partial_{x_1}u+\cdots+c_m u\partial_{x_1}u+\cdots\;.  \label{eq.map}
\end{align}
where $\partial_{\mathbf{x}}^k$ denotes all $k$-th order partial derivatives and $\partial_{x_j}$ denotes the partial derivative with respect to the $j$-th variable.
We refer to each term, such as  $1,\partial_{x_1}u$, and $u\partial_{x_1}u,\ldots$ in \eqref{eq.map}, as a \textit{feature}.
Since $f$ is a finite order polynomial, only a finite number of features are included. Denote the number of features by $K$.
Under this model, the function $f$ is expressed in a parametric form as a linear combination of $K$ features. Our objective is to recover the parameters, or coefficients,
$$\bc = [c_1 \  \ c_2\  \ldots \  c_m \ \ldots \ \ c_K]^T \in \mathbb{R}^K.$$
where many of the entries may be zero.

From $\mathbf{D}$,  we numerically approximate the time and spatial derivatives of $u$ to obtain the following \textit{approximated time derivative vector} $D_tU\in \mathbb{R}^{NM^d}$ and \textit{approximated feature matrix} $F\in \mathbb{R}^{NM^d \times K}$:
\begin{align}
\begin{tabular}{l}
$D_tU = \begin{bmatrix}
\frac{U^1-U^0}{\Delta t}
\vspace{0.2cm}
\\
\vspace{0.2cm}
\frac{U^2-U^1}{\Delta t}\\
\vdots
\vspace{0.1cm}
\\
\frac{U^{N}-U^{N-1}}{\Delta t}
\end{bmatrix} $\;,\;\; $F =\begin{bmatrix}
\mathbf{1}_{M^d\times 1}&U^0  &\cdots& U^0\circ D_{x_1}U^0 &\cdots\\
\mathbf{1}_{M^d\times 1}&U^1&\cdots& U^1\circ D_{x_1}U^1  &\cdots\\
\vdots& \vdots&\ddots&\vdots&\cdots\\
\mathbf{1}_{M^d\times 1}&U^{N-1}&\cdots& U^{N-1}\circ D_{x_1}U^{N-1}  &\cdots
\end{bmatrix}\;.$
\end{tabular} \label{approxfeaturedef}
\end{align}
In this paper, the time derivatives in $D_tU$ are approximated by the forward difference scheme, and the  spatial derivatives, such as $D_{x_1}U^n$ for $n=0,1,\dots, N-1$ in $F$ are computed using the 5-point ENO scheme \cite{harten1987uniformly}. Other numerical  differentiation schemes can be used here (See \cite{kang2019ident} for an error estimation.) The vector  $\mathbf{1}_{M^d\times 1} \in \mathbb{R}^{M^d}$ denotes the 1-vector of size $M^d$, and the Hadamard product $\circ$ is the element-wise multiplication between two vectors. Each column of $F$ is referred to as a \textit{feature column}. The PDE model in \eqref{originalPDE} suggests that, an optimal coefficient vector $\mathbf{c}$ should satisfy the following approximation:
\begin{equation}
D_t U \approx F \bc\;.
\label{lineareq}
\end{equation}
The objective of this paper is to find the correct set of coefficients in (\ref{eq.map}).  Due to the large size of $K$, the idea of sparsity becomes useful.

The framework of our methods relies on a prescribed dictionary, and  the dictionary should  contain all possible terms in the underlying PDE. If we do not have any a prior knowledge, one strategy 
is to use the pairwise product of the partial derivatives of $u$ up to certain order. 
One can also view the right-hand side of the target PDE in \eqref{originalPDE} as a functional of $u$ and its partial derivatives up to certain order, and then approximate the right-hand side by a Taylor polynomial up to certain degree. 
Our method is capable of identifying this Taylor polynomial, as an approximation to the right-hand side of the underlying PDE.
Another strategy is to estimate the possible features from the given data, which is an open problem to be studied in the future. 

Throughout this paper,  we denote $F_0$ as the true feature matrix whose elements are the exact derivatives evaluated at the corresponding time and space location as those in $F$.
For a vector $\bc$, $\|\bc\|_p := (\sum_{j}|c_j|^p)^{\frac 1 p}$ is the $L_p$ norm of $\bc$. In particular, $\|\bc\|_\infty:=\max_{j}|c_j|$. When $p=0$, $\|\bc\|_0 := \#\{c_j : c_j\neq 0\}$ represents the $L_0$ semi-norm of $\bc$. The support of $\bc$ is denoted by ${\rm supp}(\bc) :=\{j: c_j\neq 0\}$. The vector $\bc$ is said to be $k$-sparse if $\|\bc\|_0=k$ for a non-negative  integer $k$.
For any matrix $A_{m\times n} $ and index sets $\mathcal{L}_1\subseteq\{1,2,\dots,n\}$, $\mathcal{L}_2\subseteq\{1,2,\dots,m\}$, we denote $[A]_{\mathcal{L}_1}$ as the submatrix of $A$ consisting of the columns indexed by $\mathcal{L}_1$, and  $[A]^{\mathcal{L}_2}$ as the submatrix of $A$  consisting of the rows indexed by $\mathcal{L}_2$. $A^T$, $A^*$ and $A^\dagger$ denote the transpose, conjugate transpose and Moore-Penrose pseudoinverse of $A$, respectively. For $x\in \mathbb{R}$, $\lfloor x\rfloor$ denotes the largest integer no larger than $x$.

\subsection{Successively Denoised Differentiation (SDD)}\label{sec:sdd}

As shown in Figure \ref{fig.amplif}, when the given data are contaminated by noise, numerical differentiation amplifies noise. It introduces a large error in the time derivative vector $D_t U$ and the approximated feature matrix $F$.  With random noise, the regularity of the given data is different from the PDE solution's regularity.  Thus, the denoising plays a vital role in PDE identification.

We introduce a smoothing operator $S$ to process the data. Kernel methods are good options for $S$, such as Moving Average \cite{smith1997scientist} and Moving Least Square (MLS)~\cite{lancaster1981surfaces}.
In this paper, the smoothing operator $S$ is chosen as the MLS, where data are locally fit by quadratic polynomials. In MLS, a weighted least squares problem, in the time domain or the spatial domain, is solved at each time $t^n$ and spatial location $x_\mathbf{i}$ as follows:
\begin{align}
&S_{(\mathbf{x})} \left[U\right]_{\mathbf{i}}^n = p^n_\mathbf{i}(\bx_\mathbf{i}),
\text{ with}\;\; p^n_\mathbf{i}=\argmin_{p\in P_2}\sum_{\mathbf{j}\in\mathbb{I}}(p(\bx_\mathbf{j})-U_{\mathbf{j}}^n)^2	\exp\left(-\frac{\|\bx_{\mathbf{i}}-\bx_{\mathbf{j}}\|^2}{h^2}\right),
\label{eqSx}
\\
&S_{(t)} \left[U\right]_{\mathbf{i}}^n = p^n_\mathbf{i}(t^n),
\text{ with}\;\; p^n_\mathbf{i}=\argmin_{p\in P_2}\sum_{0\leq k\leq N}(p(t^k)-U_{\mathbf{i}}^k)^2	\exp\left(-\frac{\|t^n-t^k\|^2}{h^2}\right)\;.
\label{eqSt}
\end{align}
Here $h>0$ is a width parameter of the kernel, and $P_2$ denotes the set of polynomials of degree no more than $2$. 
It is shown in \cite[Theorem 4.1]{wendland2001local} that, for a fixed time index $n$, if the given data $\{U_{\mathbf{i}}^n\}_{\mathbf{i}}$ are sampled from a $C^k$ function $u(\bx,t^n)$, and the $(k-1)$th order polynomials are used in MLS, the output of MLS with a proper choice of the kernel width $h$ gives a $k$th order approximation of $u(\mathbf{x},t^n)$.  This theory demonstrates that, MLS keeps the accuracy of the given data when the data contain no noise and the solution is sufficiently smooth. In practice, the width parameter $h$ is found empirically: as the noise level increases, 
a larger $h$ is used to address the data variability. In our experiments,
we observe that the performance of our method is not sensitive to the choice of $h$ and we use the same $h$ for  different noise levels.

\begin{figure}
	\centering
	\begin{tabular}{ccc}
		(a)&(b)&(c)\\
		\includegraphics[width=0.3\textwidth]{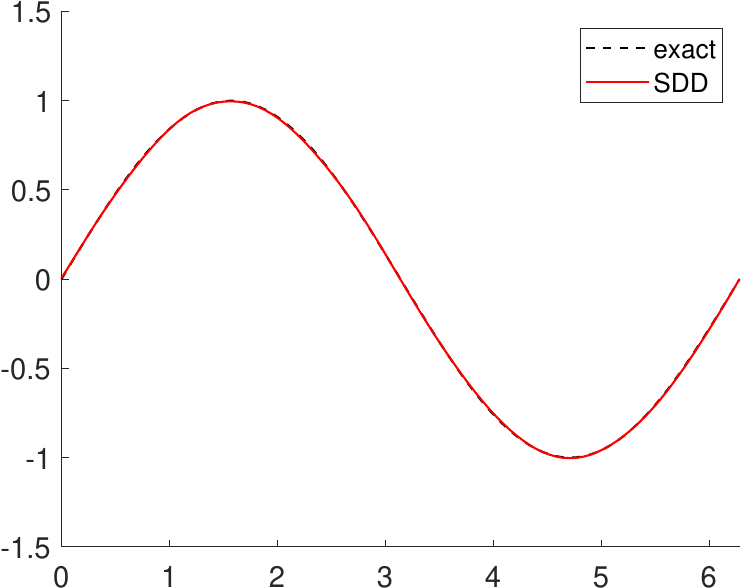}&
		\includegraphics[width=0.3\textwidth]{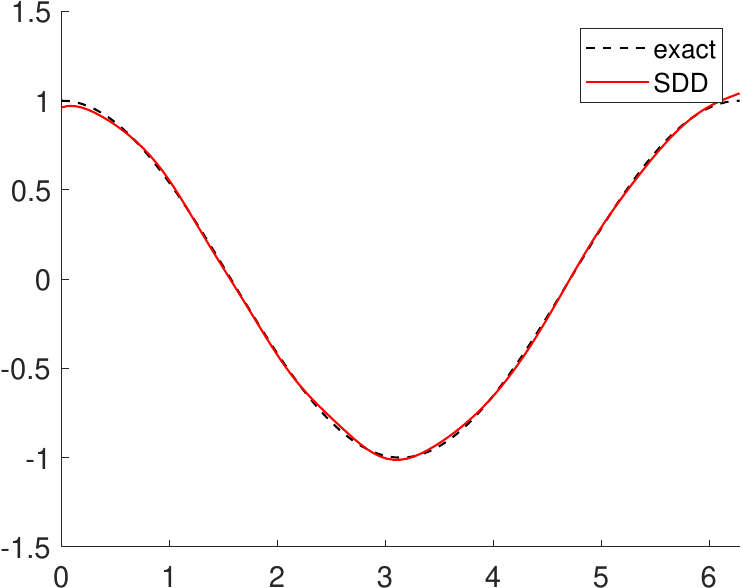}&
		\includegraphics[width=0.3\textwidth]{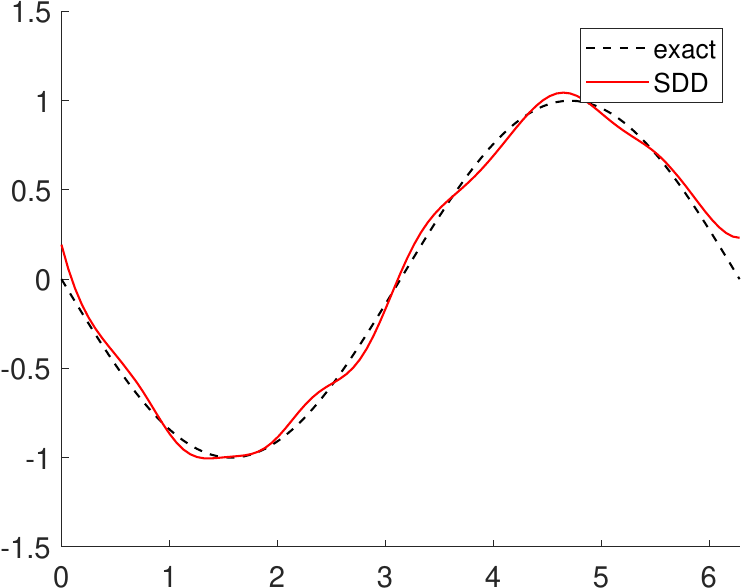}
	\end{tabular}
	\caption{Performance of SDD on the data in Figure \ref{fig.amplif}. (a) Graph of $\sin( x)$, $0\leq x\leq 2\pi$ (black) and the denoised data (red) using MLS.  (b)  First-order derivatives of the function (black) and the denoised data using SDD (red). (c)  Second-order derivatives of the function (black) and the denoised data using SDD (red).  Derivatives are computed by the five-point ENO scheme, and the smoothing operator $S$ is MLS.}
	\label{fig.sddexp1}
\end{figure}

We propose a Successively Denoised Differentiation (SDD) procedure to stabilize the numerical differentiation.    For every derivative approximation, smoothing is applied as described in Table~\ref{tab_smoothing_sum}. 
\begin{table}

\begin{center}
	\begin{tabular}{|c|l| }
		\hline
		Term & Approximation \\
		\hline
		$u$ & $\approx S_{(\mathbf{x})}[U]$ \\
		& $u$ is approximated by spatially MLS denoised data $U$. \\ \hline
		$\partial_{t} u$ & $ \approx S_{(t)}D_{t}S_{(\mathbf{x})}[U]$ \\
		& Time-domain denoising applied after numerical time differentiation. \\  \hline
		$\partial_{\mathbf{x}}^{\mathbf{k}} u$ 
		& $\approx (S_{(\mathbf{x})}D_{x_1})^{k_1}\cdots (S_{(\mathbf{x})} D_{x_d})^{k_d}S_{(\mathbf{x})}[U]$ 
		 where $\mathbf{k} =(k_1,\ldots,k_d)$ \\
		& Spatial denoising applied after every numerical spatial differentiation. \\
		\hline
	\end{tabular}
		\caption{ Examples of SDD: Each (differential) terms are approximated by the spatial and time smoothing operators $S_{(\mathbf{x})}$ and $S_{(t)}$ defined in \eqref{eqSx} and \eqref{eqSt} respectively.  The operator $D_t$ given in \eqref{approxfeaturedef} represents numerical time differentiation by the forward difference scheme, and $D_{x_j}$ for $j=1,\dots, d$ represents numerical spatial differentiation with respect to $x_j$ given by the 5-point ENO scheme \cite{harten1987uniformly}.} \label{tab_smoothing_sum}
\end{center}
\end{table}
The main idea of SDD is to smooth the data at each step (before and after) the numerical differentiation. This simple idea effectively stabilizes numerical differentiation.  Figure \ref{fig.sddexp1} shows the results of SDD for the same data in Figure \ref{fig.amplif}.  The approximations of the first and second-order derivatives of $u$ are greatly improved by SDD. 

When the given data are noiseless and MLS is used in SDD, the following theorem shows that under appropriate assumptions, the estimated partial derivative $S_{(x)}D_x[U]$ has the same accuracy as the estimated derivative without SDD.

\begin{theorem}\label{thm.SDD}
	Let $s$ be a positive integer. Suppose the given data are noiseless and are sampled from a sufficiently smooth function $u(x)$ with spacing $\Delta x$, i.e. $U_i=u(x_i)$ with $x_i = i \Delta x$. Assume polynomials of degree $s$ are used in MLS and the width parameter $h$ is properly chosen.  Let $D_x$ be a linear difference scheme satisfying (a) at every $x_i$, $D_x v=dv/dx + O(\Delta x^s)$ for any sufficiently smooth function $v(x)$; and (b) $D_x v(x_i)=\sum_{-J\le j\le J} d_j v(x_{i+j})$ for some
	positive integer $J$, where $d_j$
	depends only on $\Delta x$ and $d_j=O(1/\Delta x)$, $\forall j$. 
Then for $k=0,1,\cdots, s$, at every $x_i$
$$(S_{(x)}D_x)^k S_{(x)}[u]=d^k u/dx^k + O(\Delta x^{s+1-k})~.$$
\end{theorem}

\begin{proof}
	
	According to \cite[Theorem 4.1]{wendland2001local}, at every $x_i$ 
\begin{equation}
\label{MLS-accuracy}
S_{(x)}[v] = v+O(\Delta x^{s+1})
\end{equation}
for any sufficiently smooth function $v(x)$. Therefore the case $k=0$ has been proved.
Note that $S_{(x)}[v](x_i)=\sum _{-I\le j\le I} a_j v(x_{i+j}) $, for some positive
integer $I$ and coefficients $a_j$.
Also from \cite[Proof of Theorem 4.1]{wendland2001local}, we can deduce that for any
function $w(x)$, at every $x_i$
\begin{equation}
\label{MLS-bound}
|S_{(x)}[w](x_i)|\leq C\max_{-I\le j\le I} |w(x_{i+j})|
\end{equation}
for some constant $C$.
	Now assume at every $x_i$, 
	$$(S_{(x)}D_x)^k S_{(x)}[u]=d^k u/dx^k + O(\Delta x^{s+1-k})$$
	for some $k\in \{0,1,\cdots, s-1\}$.
	We want to show that at every $x_i$, 
	$$(S_{(x)}D_x)^{k+1} S_{(x)}[u]=d^{k+1} u/dx^{k+1} + O(\Delta x^{s-k})~.$$
	Decompose the error at every $x_i$ as
	\begin{equation}
	\begin{array}{l}
		|(S_{(x)}D_x)^{k+1}S_{(x)}[u]-\frac{d^{k+1} u}{dx^{k+1}}|  \nonumber\\
		\le |S_{(x)}D_x[(S_{(x)}D_x)^{k}S_{(x)}[u]-\frac{d^{k} u}{dx^{k}}]|+
		                   |S_{(x)}D_x \frac{d^{k} u}{dx^{k}}-\frac{d^{k+1} u}{dx^{k+1}} | \nonumber\\
		=(A) + (B).
		\end{array}
		\label{eq.proof.decomposition}
	\end{equation}
Using the induction assumption and assumption (b) of $D_x$, we have 
$$|D_x[(S_{(x)}D_x)^{k}S_{(x)}[u]-\frac{d^{k} u}{dx^{k}}]|=O(\Delta x^{s-k})~.
$$
Using (\ref{MLS-bound}), we have $(A)= O(\Delta x^{s-k}) $. 
Using property (a) of $D_x$, we have 
$$D_x \frac{d^{k} u}{dx^{k}} = d^{k+1} u/dx^{k+1}+ O(\Delta x^s)~.$$

Combining property (\ref{MLS-bound}) and (\ref{MLS-accuracy}), we have
$$S_{(x)}D_x \frac{d^{k} u}{dx^{k} }= d^{k+1} u/dx^{k+1}+ O(\Delta x^s)~,$$
thus $(B)= O(\Delta x^s)$.
Finally, we conclude that 
$$(S_{(x)}D_x)^{k+1} S_{(x)}[u]=d^{k+1} u/dx^{k+1} + O(\Delta x^{s-k})$$
and the theorem is proved.
\end{proof}

{\bf Remark.} Suppose $v(x)$ is interpolated by a polynomial $p(x)$ of degree $s$ 
using Lagrangian interpolation on $s+1$ grid points near $x_i$, then $d p/dx$ is a linear difference scheme
and 
satisfies assumption (a) and (b) of Theorem~\ref{thm.SDD}. In particular, the difference scheme used in this
paper is constructed this way.

%
%
%
%
%
%

Theorem \ref{thm.SDD} implies that under proper settings, the estimated derivative by SDD has the same accuracy as the estimated derivative without SDD. Following the proof of Theorem \ref{thm.SDD}, one can easily derive similar results for higher-order derivatives in multidimensions.

In Section \ref{subsec::NoiseAnalysis}, we explore details of SDD when different smoothing operators are used.  We find that MLS has the best performance in terms of preserving the derivative profiles.  Therefore, we set $S$ to be MLS in our numerical experiments.

To simplify the notations, in the rest of this paper, we use $U$ to denote the denoised data $S_{(\mathbf{x})}[U]$, and $D_tU$ as well as $D_{\mathbf{x}}^kU$ to denote the numerical derivatives with SDD applied as above.

\section{Proposed Methods: ST and SC }\label{sec:methods}

Under the parametric model in \eqref{eq.map}, the PDE identification problem can be reduced to solving the linear system \eqref{lineareq} for a sparse vector $\bc$ with few nonzero entries.  Sparse regression can be formulated as the following $L_0$-minimization
\begin{equation}
\min \|\bc\|_0\;, \quad \text{ subject to } \|F \bc - D_t U\| \le \epsilon\;,
\label{l0min}
\end{equation}
for some $\epsilon >0$.
However, the $L_0$-minimization in \eqref{l0min} is NP-hard.  Its approximate solutions have been intensively studied in the literature.
The most popular surrogate for the $L_0$ semi-norm is the $L_1$ norm applied in image and signal processing~\cite{candes2006robust,donoho2006compressed}. The $L_1$-regularized minimization is called Least Absolute Shrinkage and Selection Operator (LASSO) \cite{tibshirani1996regression}, which was used in \cite{kang2019ident,schaeffer2017learning,rudy2017data} for PDE identification. The common strategy in these works is to utilize LASSO to select a candidate set, then refine the results with other techniques.

In this paper, we utilize a greedy algorithm called Subspace Pursuit (SP) \cite{dai2008subspace} to select a candidate set.   Unlike LASSO, SP  takes the sparsity as an input, allowing direct control of the sparsity of the reconstructed coefficient. Let $k$ be a positive integer and denote $\bb = D_t U$. For a fixed sparsity level $k$, SP$(k;F,\bb)$ in Algorithm \ref{SPalgo} gives rise to a $k$-sparse vector whose support is selected in a greedy fashion.  It was proved that SP gives rise to a solution of the $L_0$-minimization \eqref{l0min} under certain conditions of the matrix $F$, such as the restricted isometry property \cite{dai2008subspace}.

\begin{algorithm2e}
	\KwIn{$F\in\mathbb{R}^{NM^d\times K}$, $\bb \in\mathbb{R}^{NM^d}$ and sparsity $k\in\mathbb{N}$.}
	
	\textbf{Initialization:} $j=0$;\\
	$G \leftarrow$ column-normalized version of $F$;
	\\
	$\mathcal{I}^0=\{k$ indices corresponding to the largest magnitude entries in the vector $G^*\mathbf{b}\}$;
	\\
	$\mathbf{b}_{\text{res}}^{0}= \mathbf{b}-G_{\mathcal{I}^0}G_{\mathcal{I}^0}^\dagger\mathbf{b}$.
	
	\While{True}{
		\textbf{Step 1.} $\widetilde{\mathcal{I}}^{j+1}=\mathcal{I}^{j}\cup\{k$ indices corresponding to the largest magnitude entries in the vector $G^*\mathbf{b}_{\text{res}}^{j}\}$\;
		
		\textbf{Step 2.} Set $\mathbf{c}_p=G_{\widetilde{\mathcal{I}}^{j+1}}^{\dagger}\mathbf{b}$\;

		\textbf{Step 3.} $\mathcal{I}^{j+1}=\{k$ indices corresponding to the largest elements of $\mathbf{c}_p\}$\;
		
		\textbf{Step 4.} Compute $\mathbf{b}_{\text{res}}^{j+1}=\mathbf{b}-G_{\mathcal{I}^{j+1}}G_{\mathcal{I}^{j+1}}^\dagger \mathbf{b}$\;
		
		\textbf{Step 5.} If $|\mathbf{b}_{\text{res}}^{j+1}\|_2>\|\mathbf{b}_{\text{res}}^{j}\|_2$, let $\mathcal{I}^{j+1} = \mathcal{I}^j$ and terminate the algorithm; otherwise set $j\gets j+1$ and iterate.
	}
	
	\KwOut{$\widehat{\mathbf{c}}\in\mathbb{R}^{K}$ satisfying $\widehat{\mathbf{c}}_{\mathcal{I}_j}=F_{\mathcal{I}_j}^\dagger \mathbf{b}$ and $\widehat{\mathbf{c}}_{(\mathcal{I}_j)^\complement} = \mathbf{0}.$
	}
	\caption{{\bf Subspace Pursuit $\text{SP}(k;F,\bb)$}  \label{SPalgo}}	
\end{algorithm2e}

We propose two new methods based on SP for PDE identification: Subspace pursuit Time evolution (ST) and Subspace pursuit Cross-validation (SC). 

\subsection{Subspace Pursuit Time Evolution (ST)}\label{subsec::ST}

\begin{figure}
	\centering
	\includegraphics[scale=0.32]{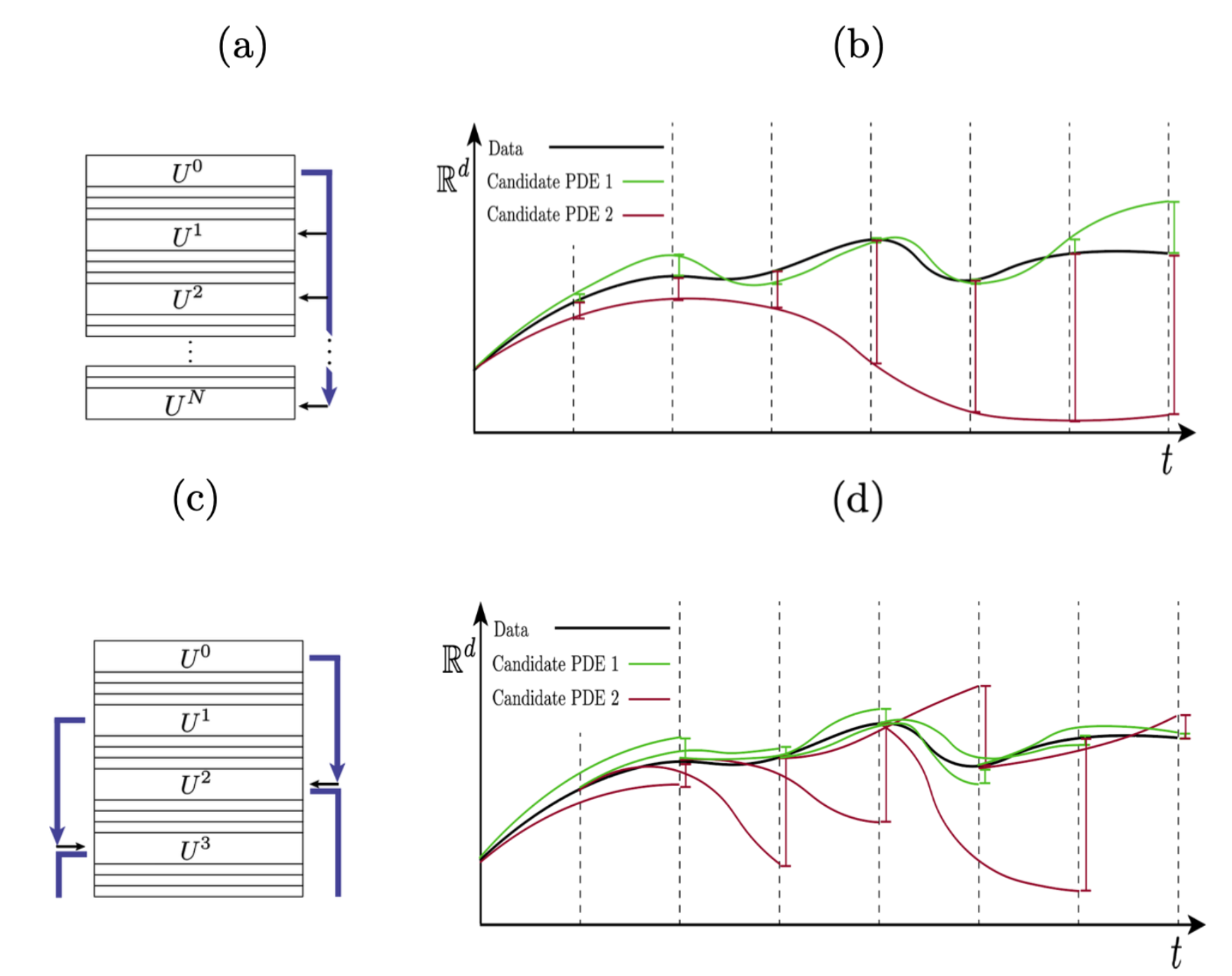}
	\caption{(a) and (b) illustrate the idea of TEE. (c) and (d) illustrate MTEE when $w=2$. The blue arrows in (a) and (c) represent time evolution using the forward Euler scheme on a fine time grid with spacing $\widetilde{\Delta t} \ll {\Delta t}$.  In (b), two different PDEs (green and red) are evolved, and the green one has a smaller TEE.  In (d), the candidate PDEs are evolved from multiple time locations, and their numerical solutions are compared with the denoised data after a time length of $w\Delta t$. }
	\label{fig:TEE/modTEEdemo}
\end{figure}

We first propose a method combining SP and the idea of time evolution.  In  \cite{kang2019ident},  Time Evolution Error (TEE) quantifies the mismatch between the solution simulated from a candidate PDE and the denoised data. Any candidate coefficient vector $\widehat{\bc}=(\widehat{c}_1,\widehat{c}_2\dots)$ defines a candidate PDE: 
$$u_t = \widehat c_1+\widehat c_{2}\partial_{x_1}u+\cdots+ \widehat c_m u\partial_{x_1}u+\cdots.$$  
This PDE is numerically evolved from the initial condition $U^0$ with a smaller time step $\widetilde{\Delta t}\ll\Delta t$. 
Denote $\widehat U^1,\widehat U^2,\ldots, \widehat U^N$ as this numerical solution at the same time-space location as  $U^1, U^2,\ldots, U^N$. The TEE of the candidate PDE given by $\widehat\bc$ is
$${\rm TEE}(\widehat\bc) =\frac{1}{N} \sum_{n=1}^N \|\widehat U^n - U^n \|_2\;, $$
where $U^n$ is the denoised data at time $t^n$. Figure~\ref{fig:TEE/modTEEdemo}~(a) and (b) illustrate the idea of TEE.
When there are several candidate PDEs, the one with the least TEE is picked  \cite{kang2019ident}. This TEE idea is based on the convergence principle that a correct numerical approximation converges to the true solution as the time step $\widetilde{\Delta t}$ goes to zero. The error from the wrongly identified terms grows during this time evolution process, see more details in \cite[Section 2.3]{kang2019ident}. 


In this paper, we propose a \textit{Multi-shooting Time Evolution Error}~(MTEE). The idea is to evolve a candidate PDE from multiple time locations with a time step $\widetilde{\Delta t}\ll\Delta t$  using the forward Euler scheme for a time length of  $w\Delta t$, where $w$ is a positive integer.    
This scheme is stable as long as the PDE is well posed and the solution is smooth, and when the time step is sufficiently small. 
Specifically, if $r$ is the highest order of the spatial derivatives, following the CFL condition, we set the time step as $c(\Delta x)^r$ with some constant $c<1$. 
Let $\widehat{U}^{(n+w)|n}$ be the numerical solution of the candidate PDE at the time $(n+w)\Delta t$, which is evolved from the initial condition $U^n$ at time $t^n = n\Delta t$.
The MTEE is defined as
\begin{align}
\text{MTEE}(\widehat{\mathbf{c}};w) = \frac{1}{N-w}\sum_{n=0}^{ N-1-w}\|\widehat{U}^{(n+w)|n}-U^{n+w}\|_2\label{MTEEeq}\;.
\end{align}
Figure~\ref{fig:TEE/modTEEdemo}~(c) and (d) demonstrate the process of multi-shooting time evolution.   While the TEE evolution starts from the initial condition $U^0$ and ends at $T$, the MTEE evolution starts from various time locations, such as $t^n, n=0,\ldots,N-1-w$, and lasts for a shorter time, e.g., $w\Delta t$ in our case.

MTEE has two advantages over TEE:  
(1) MTEE is more robust against noise in comparison with TEE. If $w \ll N$,  the noise in the initial condition 
accumulates for a smaller amount of time in MTEE, which helps to stabilize numerical solvers.

\begin{figure}
	\centering
	\begin{tabular}{ccc}
		(a)&(b)&(c)\\
		\includegraphics[width=0.3\textwidth]{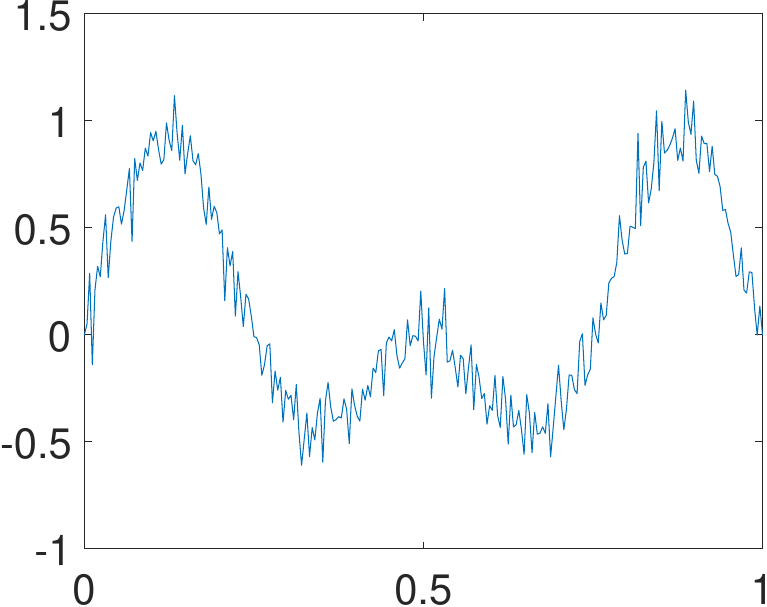}&
		\includegraphics[width=0.3\textwidth]{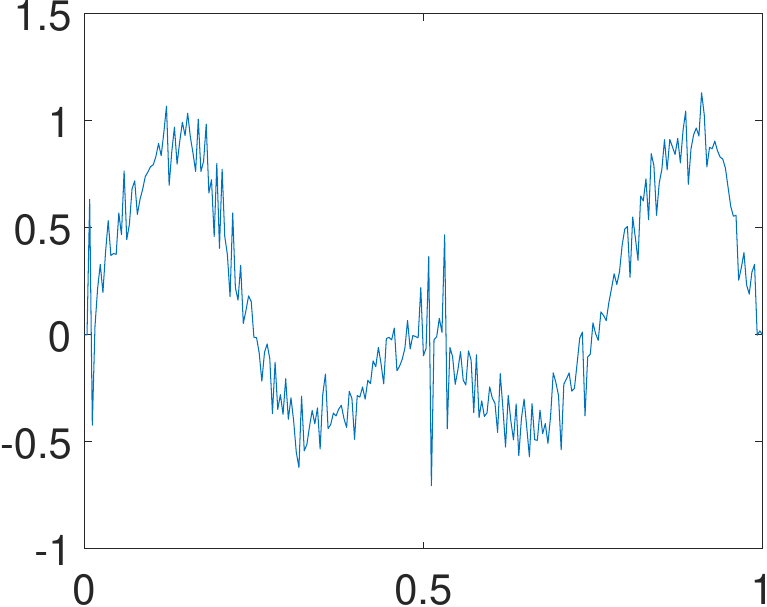}&
		\includegraphics[width=0.3\textwidth]{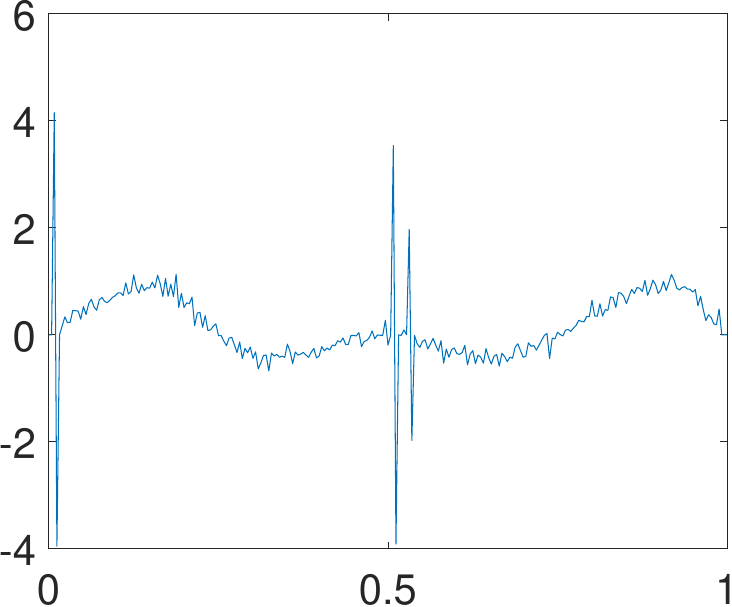}
	\end{tabular}
	\caption{Robustness of MTEE over TEE. (a) A noisy initial condition for the evolution of the Burgers' equation $u_t=-uu_x$. 
		By evolving this noisy initial condition according to $u_t=-uu_x$, (b) shows the numerical solution at $t=0.02$ and (c) shows the numerial solution at $t=0.03$. The numerical solution blows up at $t=0.032$.}
	\label{fig.MTEE}
\end{figure}

For example, consider identifying the Burgers' equation $u_t=-uu_{x}$ from a set of noisy data generated with $T=0.05, \Delta t=0.001, \Delta x=1/256$ (see Figure \ref{fig.MTEE}). 
If one evolves the noisy initial condition in Figure \ref{fig.MTEE}(a), using the correct PDE, i.e., $u_t=-uu_{x}$, the numerical solution blows up at $t=0.032$. The numerical solutions at $t=0.02$ and $t=0.03$ are shown in Figure \ref{fig.MTEE} (b) and (c), respectively.  The TEE at $T =0.05$ is $\infty$ even for the correct PDE since the numerical solution blows up at $t=0.032$.  
On the other hand, MTEE works since we evolve the initial condition for a shorter amount of time, before the numerical solution blows up, such as $t=0.02$ (corresponding to $w=20$ in MTEE) in this example.
(2) MTEE is more flexible, and its computation is parallelizable. The flexibility of MTEE comes from two aspects: (i) The error accumulation time can be controlled by the parameter $w$ such that the PDE is evolved for a time length of $w\Delta t$. (ii) One may assign different weights in the calculation of the evolution errors in different periods. Since each time evolution in the multi-shooting is independent, the computation of MTEE can be parallelized. 

\begin{figure}
	\centering
\includegraphics[height=1.2in]{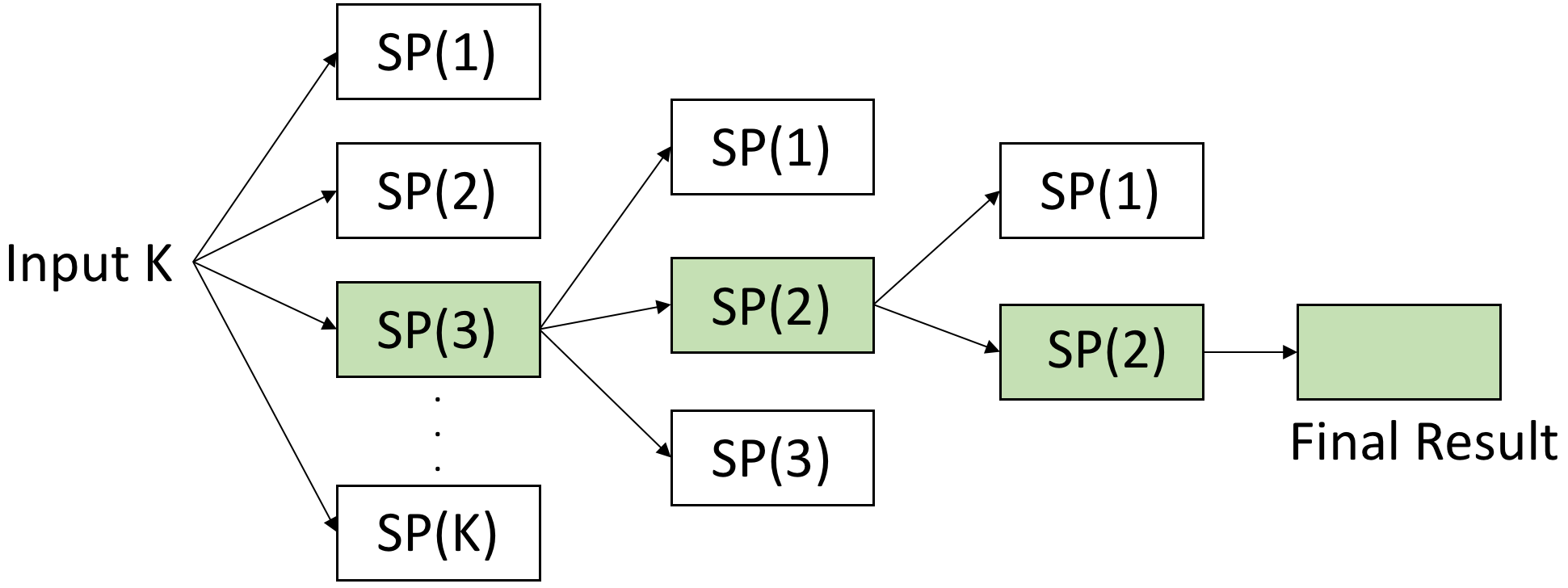}
	\caption{An example of the ST iteration. Starting with the large number $K$, the first iteration gives rise to $K$ candidate coefficients for $k=1,\dots,K$.  The PDE with the smallest MTEE is picked, e.g., SP(3) with cardinality $K_1=3$ and support $\mathcal{A}_1$.   The second iteration gives rise to the candidate coefficients only supported on $\mathcal{A}_1$ using SP(k) with $k=1,2,3$.  The PDE with the smallest MTEE is found, e.g., SP(2) with cardinality $K_2=2$ and support $\mathcal{A}_2$.  The third iteration does not change the support, i.e., $\mathcal{A}_3 = \mathcal{A}_2$, so the final output is the coefficient vector of ${\rm SP}(2)$. }
	\label{fig:TEE/modTEEdemo_2}
\end{figure}

The SP algorithm finds a coefficient vector with a specified sparsity, while the correct sparsity is not known from the given data. Based on SP and MTEE, we propose \textit{Subspace pursuit Time evolution} (ST), which iteratively refines the selection of features.    Figure \ref{fig:TEE/modTEEdemo_2} illustrates the ST iteration: Starting from a large number  $K$ (no more than the number of features), each SP($k$) coefficient vector is computed for all $k=0,\dots, K$.  Among these, the $k$ which gives the minimum MTEE is chosen to be $K_1$.  This procedure continues until two consecutive iterations give the same output or only one feature is left.   This process will terminate after at most $K-1$ iterations.  

More specifically,  as an initial condition, we set $K_0=K$ and $\mathcal{A}_0 = \{1,\ldots,K\}$.  Clearly, this $K$ is bounded by the number of dictionary. 
At the first iteration, all possible sparsity levels are considered upto $K$ in the SP algorithm. For each $k=1, \dots, K$, we run $\text{SP}(k;F,D_tU)$ to obtain a coefficient vector $\widehat{\mathbf{c}}^{(k)}\in\mathbb{R}^{K}$ such that $\|\widehat{\bc}^{(k)}\|_0=k,$ which gives rise to the PDE:
\begin{equation}
u_t =f_{{\rm SP}(k)} \text{ where } f_{{\rm SP}(k)}:= \widehat c^{(k)}_1+\widehat c^{(k)}_{2}\partial_{x_1}u+\cdots+ \widehat c^{(k)}_m u\partial_{x_1}u+\cdots.
\label{stpde}
\end{equation}
We then numerically evolve each PDE $u_t =f_{{\rm SP}(k)}$, for $k=1,\ldots,K$ and calculate the corresponding MTEE. Among these PDEs, the one with the smallest MTEE is selected, then  let
$$ K_1 = \,\argmin_{k=1,2,\cdots,K}\text{MTEE}(\widehat{\mathbf{c}}^{(k)};w)
\text{ and } \mathcal{A}_1 = {\rm supp}(\widehat{\mathbf{c}}^{(K_1)})\;.
$$
If $\mathcal{A}_1 = \mathcal{A}_0$, the algorithm is terminated; otherwise, we continue to the second iteration. The proposed method requires solving the sparsity-constrained least-squares problems at least $K$ times.  These computations and the evaluation of MTEE can be computed  in parallel. 

At the second iteration, we refine the selection from the index set $\mathcal{A}_1$ with cardinality $K_1$.
For $k=0, \dots, K_1$, we run $\text{SP}(k;[F]_{\mathcal{A}_1},D_tU)$ to obtain a coefficient vector $\widehat{\mathbf{c}}^{(k)}\in\mathbb{R}^{K}$ such that
\begin{center}
	$\widehat{\mathbf{c}}^{(k)}_{\mathcal{A}_1} = \text{SP}(k;[F]_{\mathcal{A}_1},D_tU)\;,\text{ and } \widehat{\mathbf{c}}^{(k)}_{\mathcal{A}_1^\complement} = \mathbf{0}\;,$
\end{center} and the associated PDE $u_t =f_{{\rm SP}(k)}$ as in  \eqref{stpde}. Among these PDEs, the one with the smallest MTEE is selected, and we denote
$$ K_2 = \,\argmin_{k=1,2,\cdots,K_1}\text{MTEE}(\widehat{\mathbf{c}}^{(k)};w)\;,
\text{ and } \mathcal{A}_2 = {\rm supp}(\widehat{\mathbf{c}}^{(K_1)})\;.
$$
If $\mathcal{A}_2 = \mathcal{A}_1$, the algorithm is terminated; otherwise, we continue to the next iteration similarly.

The ST iteration will be terminated when the index set remains the same, i.e., $\mathcal{A}_j=\mathcal{A}_{j+1}$.  The ST outputs a recovered coefficient vector and the corresponding PDE denoted by ST($w$). A complete description of ST is given in Algorithm~\ref{STalgo}.

\begin{algorithm2e}
	\KwIn{$F\in\mathbb{R}^{NM^d\times K}$, $D_tU\in\mathbb{R}^{NM^d}$ and a positive integer $w$.}
	
	{\bf Initialization:} $j=0$, $K_0=K$ and $\mathcal{A}_0 =\{1,2,\cdots, K\}$.
	
	\While{$\mathcal{A}_{j+1} \neq \mathcal{A}_j$}{
		\textbf{Step 1.} For $k=1,2,\cdots,K_{j}$, run $\text{SP}(k;[F]_{\mathcal{A}_j},D_tU)$ to obtain a coefficient vector $\widehat{\mathbf{c}}^{(k)}\in\mathbb{R}^{K}$ such that
		\begin{center}
			$\widehat{\mathbf{c}}^{(k)}_{\mathcal{A}_j} = \text{SP}(k;[F]_{\mathcal{A}_j},D_tU)$ and $\widehat{\mathbf{c}}^{(k)}_{\mathcal{A}_j^\complement} = \mathbf{0}\;,$
		\end{center} and the associated PDE $u_t =f_{{\rm SP}(k)}$ given in \eqref{stpde}.
		
		\textbf{Step 2.} Among all the PDEs  $u_t =f_{{\rm SP}(k)}$ for $k=1,\ldots, K_j$, select the one with the minimum $\text{MTEE}(\widehat{\mathbf{c}}^{(k)};w)$ and update
		$$ K_{j+1} = \,\argmin_{k=1,2,\cdots,K_j}\text{MTEE}(\widehat{\mathbf{c}}^{(k)};w)
		\text{ and } \mathcal{A}_{j+1} = {\rm supp}(\widehat{\mathbf{c}}^{(k_{j+1})})\;.
		$$
		If $\mathcal{A}_{j+1}=\mathcal{A}_j$, terminate the algorithm; otherwise, update $j= j+1$.}
	
	\KwOut{Recovered coefficient $\widehat\bc^{K_{j+1}}$ and the corresponding PDE, denoted by ST($w$).}
	\caption{{\bf Subspace pursuit Time evolution (ST)} \label{STalgo}}
\end{algorithm2e}

\subsection{Subspace Pursuit Cross Validation (SC)}\label{subsec::SC}

Our second method utilizes the idea of cross-validation for the linear system in  \eqref{lineareq}. Cross-validation is commonly used in statistics for the choice of parameters in order to avoid overfitting~\cite{hastie2009elements}. We consider the two-fold cross-validation where data are partitioned into two subsets. One subset is used to estimate the coefficient vector, and the other one is used to validate the candidates. If a suitable coefficient vector is found within one subset, it should yield a small validation error for the other subset because of consistency.

For some fixed ratio parameter $\alpha\in(0,1)$, we split the rows of $D_t U \in \mathbb{R}^{NM^d}$ (and $F \in \mathbb{R}^{NM^d \times K}$) into two groups indexed by $\calT_1$ and $\calT_2$, such that  $\calT_1$ consists of the indices of the first $\lfloor \alpha NM^d\rfloor $ rows 
and $\calT_2$ consists of the indices of the rest of the rows.  Since we focus on PDEs with constant coefficients, the idea of cross validation is applicable: if a correct support is identified, the coefficient vector obtained from the data in $\calT_1$ should be compatible with the data in $\calT_2$.

We introduce our \textit{Subspace pursuit Cross-validation (SC)} algorithm where cross-validation is incorporated into the SP algorithm. SC consists of the following three steps:

\noindent
\textbf{Step 1}: For each sparsity level $k=1,2,...,K$, use SP to select a set of active features:
$$
\calA_k=\mathrm{supp}(\text{SP}(k;F,D_tU))\;.
$$
\textbf{Step 2}: Use the data in $\calT_1$ to compute the estimator for the coefficient vector, $\widehat{\bc}^{(k)}\in \mathbb{R}^K$, by the following least squares problem
\begin{align*}
\widehat{\bc}^{(k)}=\argmin \limits_{{\mathbf{c}}\in \mathbb{R}^K \text{such that}\  \bc_{\calA_k^{\complement}}=0} \|[F]_{\calA_k}^{\calT_1}{\mathbf{c}_{\calA_k}}-[D_tU]^{\calT_1}\|_2^2\;, 
\end{align*}
and then use the data in $\calT_2$ to compute a Cross-validation Estimation Error (CEE)
\begin{align}
\mathrm{CEE}(\calA_k;\alpha, \calT_1,\calT_2)=\|[D_tU]^{\calT_2}-[F]^{\calT_2}\widehat{\bc}^{(k)}\|_2\;.\label{eq.CEE}
\end{align}
\textbf{Step 3}: Set
$
k_{\min}=\argmin_k \mathrm{CEE}(\calA_k;\alpha, \calT_1,\calT_2)
$
and the estimated coefficient vector is given as
\begin{align*}
\widehat{\mathbf{c}}=
\argmin \limits_{{\mathbf{c}}\in \mathbb{R}^K \text{such that}\  \bc_{\calA_k^{\complement}}=0} \|[F]_{\calA_{k_{\rm min}}}^{\calT_1}{\mathbf{c}_{\calA_{k_{\rm min}}}}-[D_tU]^{\calT_1}\|_2^2\;. 
\end{align*} The identified PDE by SC is denoted as SC($\alpha$).

CEE in \eqref{eq.CEE} is an effective measure for consistency. If the estimated coefficient vector's support matches that of the true one, CEE is guaranteed to be small provided with sufficiently high resolution in time and space.

\begin{theorem}\label{CEEbound}
	Assume that $D_tU\rightarrow u_t$ and $F\rightarrow F_0$ pointwise as $\Delta t,\Delta x\to 0$. Let $\calA_0=\mathrm{supp}(\mathbf{c}_0)$ where $\mathbf{c}_0$ is the coefficient vector of the true PDE. For any set of support $\calA$, we have
	\begin{align*}
	\mathrm{CEE} (\calA;\alpha,\mathcal{T}_1,\mathcal{T}_2)\leq \left\|\left([F_0]_{\calA_0}^{\mathcal{T}_2} \big([F_0]_{\calA_0}^{\mathcal{T}_1}\big)^\dagger -[F_0]_{\calA}^{\mathcal{T}_2}\big([F_0]_{\calA}^{\mathcal{T}_1}\big)^\dagger\right) [u_t]^{\mathcal{T}_1}\right\|_2+g(\calA;\alpha,\mathcal{T}_1,\mathcal{T}_2)\;,
	\end{align*}
	where $g>0$ is a sum of residual terms of approximating the partial derivatives and feature matrix using data (see~\eqref{eq_g_remainder}), which is independent of $\calA_0$, such that $g\to 0$ as $\Delta t,\Delta x\to 0$.
\end{theorem}
\begin{proof}
	See Appendix \ref{ProofCEEbound}.
\end{proof}

In~\eqref{eq.CEE}, the data in $\calT_1$ serve as the training set, and the data in $\calT_2$ act as the validation set. One can also use the data in  $\calT_2$ for training and the data in $\calT_1$ for validation, which gives rise to the cross validation estimation error $\mathrm{CEE}(\calA_k;1-\alpha, \calT_2,\calT_1)$.
To improve the robustness of SC, we replace~\eqref{eq.CEE}  with the following averaged cross-validation error:
\begin{align*}
\mathrm{CEE}(\calA_k,\alpha)=\frac{1}{2}\left(\mathrm{CEE}(\calA_k;\alpha, \calT_1,\calT_2) + \mathrm{CEE}(\calA_k;1-\alpha, \calT_2,\calT_1)\right)\;.
\end{align*}

In general, one can randomly pick a part of the data as the training set and use the rest as the validation set. Our numerical experiments in Subsection \ref{subsecsccomparison} demonstrate that the splitting strategy does not affect the results.
For simplicity, we split the data according to the row index in this paper. 

The proposed SC algorithm is summarized in Algorithm \ref{SCalgo}.
In comparison with ST, SC does not involve any numerical evolution of the candidate PDE, so the computation of SC is faster.


\begin{algorithm2e}[t]
	\KwIn{$F\in\mathbb{R}^{NM^d\times K}$ and $D_tU\in\mathbb{R}^{NM^d}$;  $0<\alpha<1$ ratio of the training data.}
	
	\textbf{Step 1.} For $k=1,2,\cdots,K$, run $\text{SP}(k;F,D_tU)$ to obtain the support of the candidate coefficients
	$$
	\calA_k=\mathrm{supp}(\text{SP}(k;F,D_tU))\;.
	$$
	\textbf{Step 2.} For each $k$, compute the averaged cross validation error
	$$
	\mathrm{CEE}(\calA_k,\alpha)=\frac{1}{2}\left(\mathrm{CEE}(\calA_k;\alpha, \calT_1,\calT_2) + \mathrm{CEE}(\calA_k;1-\alpha, \calT_2,\calT_1)\right)\;.
	$$
	\textbf{Step 3.} Choose the $k$ which gives the smallest cross validation error and denote it by $k_{\min}$
	$$
	k_{\min}=\argmin_k \mathrm{CEE}(\calA_k,\alpha)\;.
	$$
	Estimate the coefficients by least squares as
	\begin{align*}
	\widehat{\mathbf{c}}=
	\argmin \limits_{{\mathbf{c}}\in \mathbb{R}^K \text{such that}\  \bc_{\calA_k^{\complement}}=0} \|[F]_{\calA_{k_{\rm min}}}^{\calT_1}{\mathbf{c}_{\calA_{k_{\rm min}}}}-[D_tU]^{\calT_1}\|_2^2\;.
	\end{align*}
	
	\KwOut{Recovered coefficient $\widehat\bc$ and the identified PDE denoted by SC($\alpha$). }
	\caption{{\bf Subspace pursuit Cross validation (SC) Algorithm} \label{SCalgo}}
\end{algorithm2e}

\section{Numerical Experiments}\label{sec::numexp}

In this section, we perform a systematic numerical study to demonstrate the effectiveness of ST and SC and compare them to IDENT \cite{kang2019ident}. To measure the identification error,
we use the following relative coefficient error $e_c$ and grid-dependent residual error $e_r$:
\begin{eqnarray}
	&e_c=\frac{\|\widehat{\mathbf{c}}-\mathbf{c}\|_1}{\|\mathbf{c}\|_1},\quad
	e_r=
	\begin{cases}
		\sqrt{\Delta x \Delta t}\|F(\widehat{\mathbf{c}}-\mathbf{c})\|_2 \mbox{ for 1D PDE}.\\
		\sqrt{\Delta x \Delta y\Delta t}\|F(\widehat{\mathbf{c}}-\mathbf{c})\|_2 \mbox{ for 2D PDE}.
	\end{cases}\;.
	\label{eq.cre}
\end{eqnarray}
The relative coefficient error $e_c$ measures the accuracy in the recovery of PDE coefficients, while the residual error $e_r$ measures the difference between the learned dynamics and the denoised one by SDD. Since each feature vector in $F$ may have different scales, $e_r$ can be different from $e_c$ in some cases.  When the given data contain noise, the features containing higher-order derivatives have greater magnitudes than the features containing lower order derivatives.  In this case, a small coefficient error in the high order terms may lead to a large $e_r$.  
We use both $e_c$ and $e_r$ to quantify the PDE identification error. 
To measure how well the solution of the identified PDE matches the dynamics of the correct PDE, we also use the following evolution error 
\begin{align}
  e_e=\Delta x \Delta t\left(\sum_n\sum_{\mathbf{i}} |u(\mathbf{x}_\mathbf{i},t^n)-\hat{u}(\mathbf{x}_\mathbf{i},t^n)|\right)
  \label{eqtimeerror}
\end{align}
where $u$ and $\hat{u}$ denote the solution of the exact and identified PDE from the same initial condition, respectively.

To generate data, we first solve the underlying PDE by forward Euler scheme using time and space step $\delta t$ and $\delta x$ (and $\delta y$) respectively, then downsample the data with time and space step $\Delta t$ and $\Delta x$ (and $\Delta y$). In the noisy case, we add Gaussian noise with standard deviation $\sigma$ to the clean data. We say that the noise is $p\%$ by setting $\sigma=\frac{p}{100}\sqrt{\frac{1}{NM^d}\sum_n \sum_{\mathbf{i}} (u(\mathbf{x}_\mathbf{i},t^n))^2}$. In the computation of $D_t U$ and the feature matrix $F$, we always use SDD with MLS with $h=0.04$ as the smoother. When MLS is used to denoise the data of two dimensional PDEs, one can either fit two-dimensional polynomials or fit one-dimensional polynomials in each dimension. In this work, we use the second approach. In ST, without specification, $\widetilde{\Delta t}=\Delta t/5$ is used.

We first consider PDEs containing partial derivatives up to the second order. Let the governing equation $f$ be a polynomial with degree up to 2. There are 10 features: $1,u,u^2,u_x,u_x^2,uu_x,u_{xx},u_{xx}^2,uu_{xx}, u_{x}u_{xx}$ in the dictionary for one dimensional PDEs. For two dimensional PDEs, there are 28 features, which contain $1,u,u_x,u_y,u_{xx},y_{xy},u_{yy}$ and their pairwise products. In the following examples, without specification, the spatial domain $[0,1]$ is used for one-dimensional PDEs and $[0,1]^2$ is used for two-dimensional PDEs. For both cases, zero Dirichlet boundary condition is used for all examples.


\subsection{Transport Equation}
Our first experiment is a transport equation with zero Dirichlet boundary condition:
\begin{equation}
u_t=-u_x\;,
\label{eq.transport1d}
\end{equation}
with an initial condition of
\[
u (x,0)=
\begin{cases}
\sin^2(2\pi x/(1-T))\cos(2\pi x/(1-T)), \mbox{ for } 0\leq x\leq 1-T,\\
0, \mbox{ otherwise}
\end{cases}\;,
\]
for $0<t\leq T$ and $x \in [0,1]$. The clean data $\mathbf{D}$ is generated by explicitly solving (\ref{eq.transport1d}) with $\delta x=\Delta x=1/256,\delta t=\Delta t=10^{-3}$ and $T=0.05$.  In theory, for the transport equation, the zero boundary condition should only be applied to the inflow boundary. We design our initial condition and choose the evolution time $T$ that the solution value is 0 at the outflow boundary during the evolution. The same setup is considered in the rest of this section.

\begin{table}[h]
	\centering
		\caption{Identification of the transport equation \eqref{eq.transport1d} with different noise levels.  In the noise-free case, applying SDD does not introduce a strong bias. The identification results (second column) by ST and SC are stable even with 30\% noise.  Here $w=20$ for ST, and $\alpha=1/200$ for SC.}\label{tab.transport1d}
	\begin{tabular}{r|c|c|c}
		\hline
		\textbf{Method}& \textbf{$0\%$ noise without SDD } & $e_c$ & $e_r$\\\hline
		ST&$u_t=-0.9994u_x$ & $6.20\times 10^{-4}$ & $4.89\times 10^{-4}$\\\hline
		SC &$u_t=-0.9993u_x-0.0010u_{xx}$ & $1.65\times 10^{-3}$ & $1.11\times 10^{-2}$\\\hline
		& \textbf{$0\%$ noise with SDD}&$e_c$&$e_r$\\\hline
		ST&$u_t=-0.9997u_x$ & $3.36\times 10^{-4}$ & $2.64\times 10^{-4}$\\\hline
		SC &$u_t=-0.9997u_x-0.0010u_{xx}$ & $1.34\times 10^{-3}$ & $1.11\times 10^{-2}$\\\hline \hline
			& \textbf{ $10\%$ noise without SDD}& $e_c$ & $e_r$\\\hline
		ST & $u_t=-3.028\times 10^{-4} u_{xx}$ & 1.00 &5.55 \\\hline
		SC & $u_t=9.4224u-2.9992u_{xx}$ & $1.04\times10$ & 5.62\\\hline
		& \textbf{ $10\%$ noise with SDD}& $e_c$ & $e_r$\\\hline
		ST, SC&$u_t=-1.0357u_x$& $3.57\times 10^{-2}$ & $2.67\times 10^{-2}$\\\hline \hline
		& \textbf{ $30\%$ noise without SDD }& $e_c$ & $e_r$\\\hline
		ST & $u_t=8.0587\times 10u-2.6316\times10^{-4}u_{xx}$ & $8.16\times10$ & $1.88\times 10$\\\hline
		SC & $u_t=8.2488\times 10u$ & $8.25\times10$ & $1.86\times 10$\\\hline
		& \textbf{ $30\%$ noise with SDD }& $e_c$ & $e_r$\\\hline
		ST, SC &$u_t=-0.9421u_x$& $5.79\times 10^{-2}$ & $4.31\times 10^{-2}$\\ \hline
	\end{tabular}
\end{table}

Table \ref{tab.transport1d} shows the results of ST(20) and SC(1/200) with various noise levels.   In practice, we have no a priori knowledge of whether the given data contain noise, so we conduct two experiments with and without SDD to check the effect of SDD on clean data. We observe that SDD makes a small difference in the noise-free case.  With clean data, SC identifies an additional  $u_{xx}$ term with a small coefficient, while ST can rule out all wrong terms. The corresponding $e_c$ and $e_r$ are both small. For 10\% or 30\% noise, the results by ST and SC with and without SDD are also shown. With SDD, both ST and SC identify the correct PDE with small $e_c$ and $e_r$ values.  SDD significantly improves the results.

\begin{figure}[h!]
	\centering
	\begin{tabular}{c c c c }
		(a) & (b) &(c) &(d)\\ \includegraphics[width=0.21\textwidth]{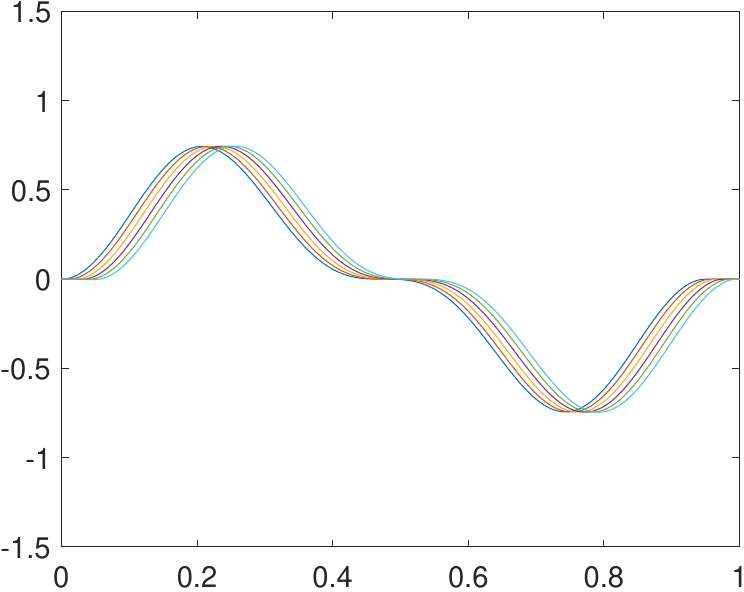}
		&
		\includegraphics[width=0.21\textwidth]{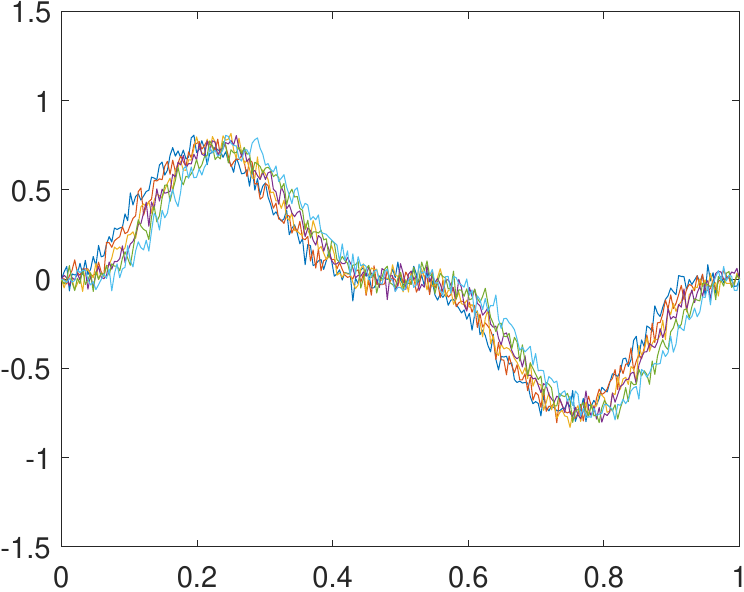} &
		\includegraphics[width=0.21\textwidth]{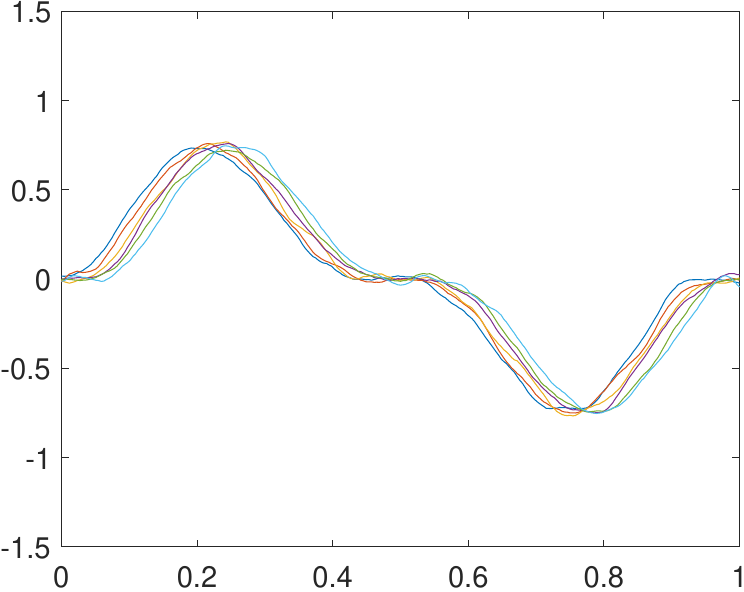} &
		\includegraphics[width=0.21\textwidth]{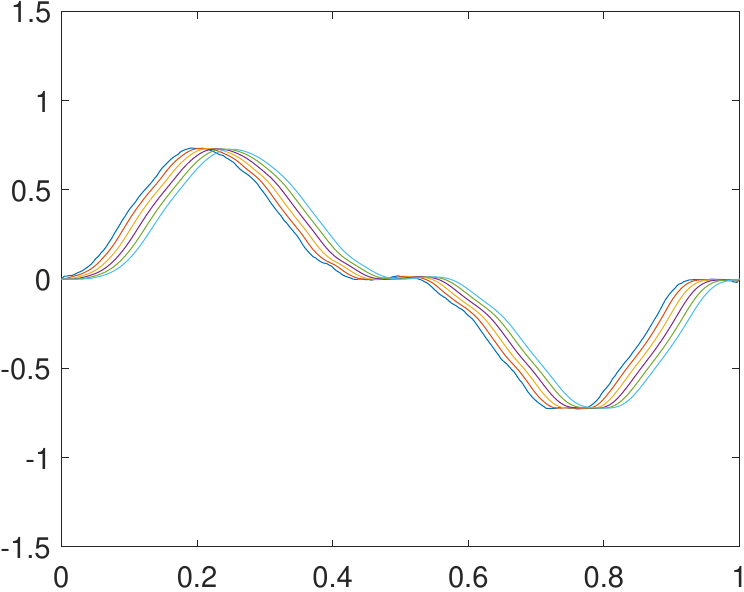}
		\\
		&  (e) & (f) & (g)  \\
		&  \includegraphics[width=0.21\textwidth]{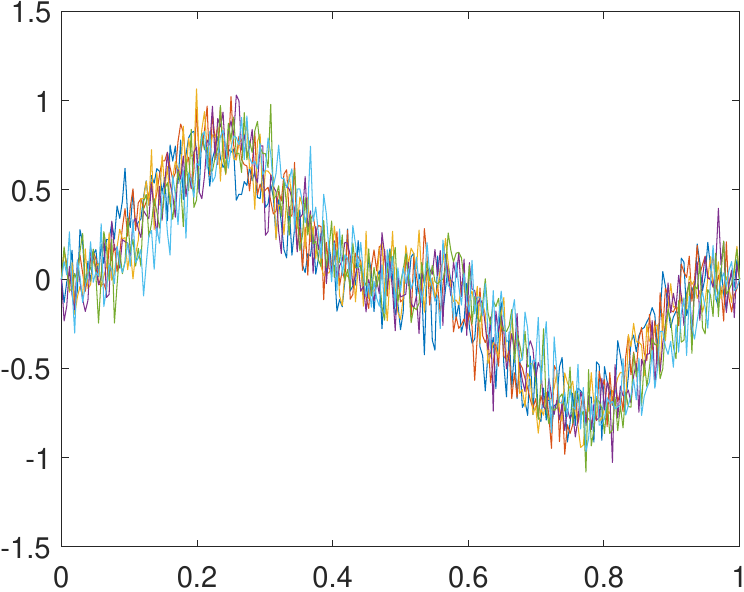} &
		\includegraphics[width=0.21\textwidth]{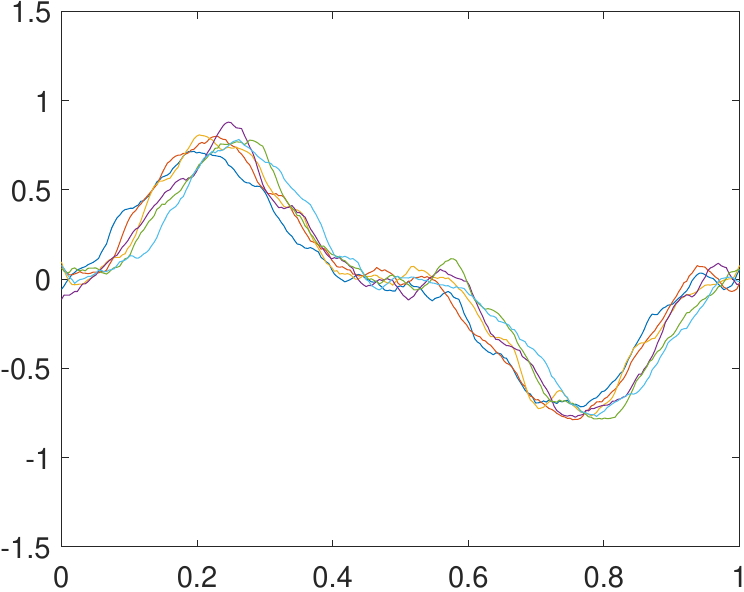} &
		\includegraphics[width=0.21\textwidth]{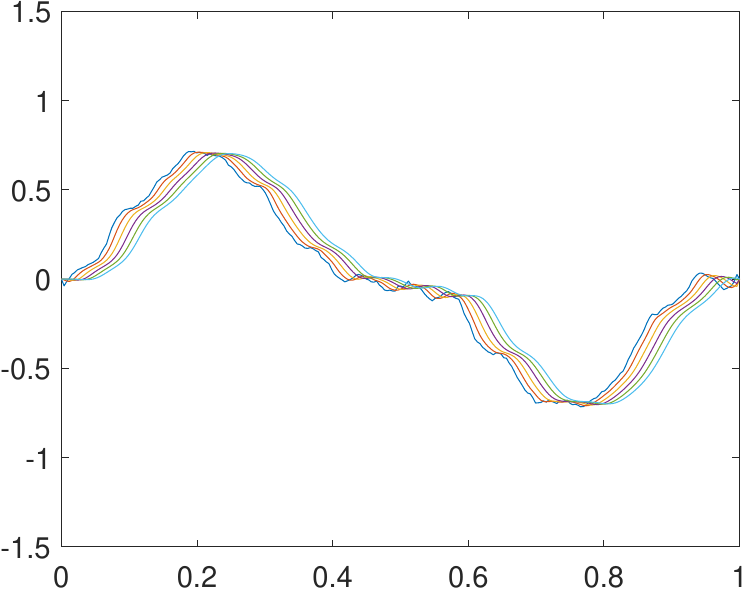}
	\end{tabular}
	\caption{Noisy and denoised data of the transport equation \eqref{eq.transport1d}, as well as simulations of the recovered PDE. (a) The clean data, (b) data with $10\%$ noise, (c) the denoised data $S_{\mathbf{x}}[U]$, (d) simulation of the PDE identified by ST and SC (identical). (e) Data with $30\%$ noise, (f) the denoised data $S_{(\mathbf{x})}[U]$, and (g) simulation of the PDE identified by ST and SC (identical).}\label{fig.transport.sol}
\end{figure}

To further demonstrate the significance of SDD and the effectiveness of ST and SC, we display the noisy data with $10\%$ and $30\%$ noise, the denoised data, and the recovered dynamics in Figure \ref{fig.transport.sol}. 
Even though the given data contain a large amount of noise, the recovered dynamics are close to the clean data. In the rest of the examples, SDD is always used for ST, SC and IDENT on noisy data.

\begin{figure}[h!]
	\centering
	\begin{tabular}{ccc}
		(a) & (b)  & (c)
		\\
		\hspace{-0.5cm}\includegraphics[width=0.33\textwidth]{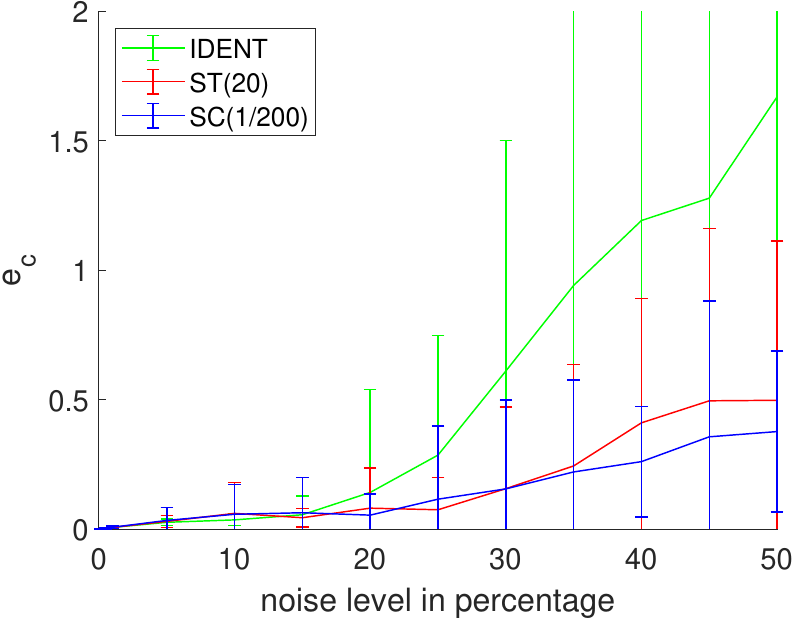} \hspace{-0.3cm}&
		\hspace{-0.3cm}\includegraphics[width=0.33\textwidth]{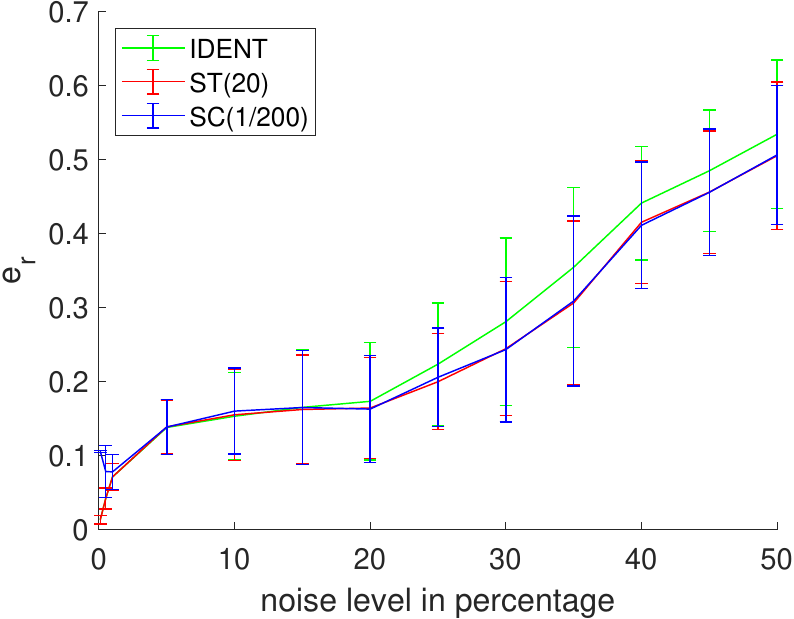}\hspace{-0.3cm}&
       \hspace{-0.3cm} \includegraphics[width=0.33\textwidth]{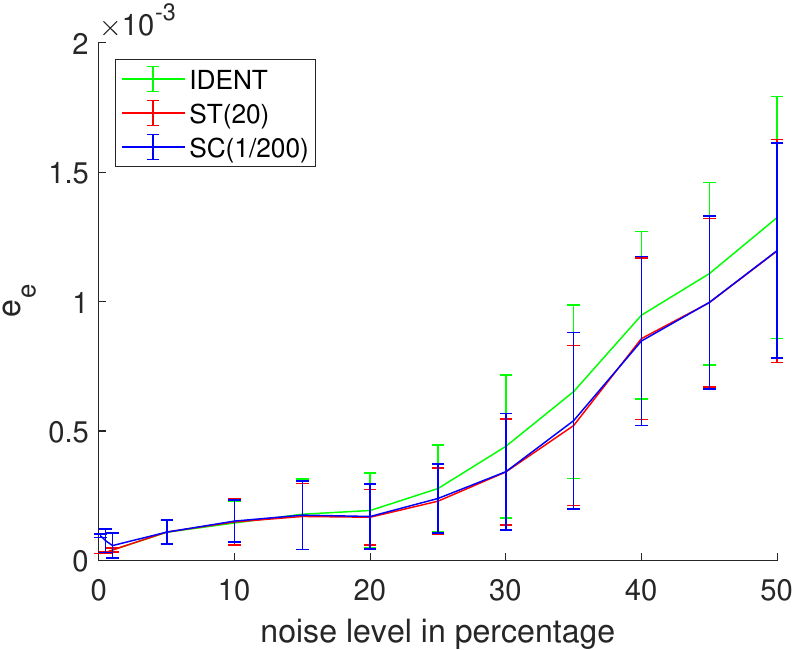}\hspace{-0.5cm}
	\end{tabular}
	\caption{The average error $e_c,e_r$ and $e_e$ over $50$ experiments for the transport equation \eqref{eq.transport1d} with respect to various noise levels.  (a) The curve represents the average $e_c$ for IDENT~\cite{kang2019ident} (Green), ST (Red) and SC (Blue), and the standard deviation is represented by vertical bars. (b) The average and variation of $e_r$ for IDENT (Green), ST (Red) and SC (Blue). (c) The average and variation of $e_e$ for IDENT (Green), ST (Red) and SC (Blue).
		The coefficient error $e_c$ by ST and SC is significantly smaller than that of IDENT.
	}\label{fig.transport1d.cre}
\end{figure}

Figure \ref{fig.transport1d.cre} shows how $e_c,e_r$ and $e_e$ change when the noise level varies. Each experiment is repeated 50 times and the error is averaged.  We test IDENT, ST(20) and SC(1/200). Figure \ref{fig.transport1d.cre} (a) shows that $e_c$ of ST or SC is much smaller than that of IDENT when the noise level is larger than 20\%.  Figure \ref{fig.transport1d.cre} (b) and (c) shows $e_r$ and $e_e$ versus noise, respectively. 
The coefficient error $e_c$ by ST and SC is significantly smaller than that of IDENT.

\begin{figure}[h!]
	\centering
	\begin{tabular}{cc}
		(a) & (b) \\
		\includegraphics[width=2.3in]{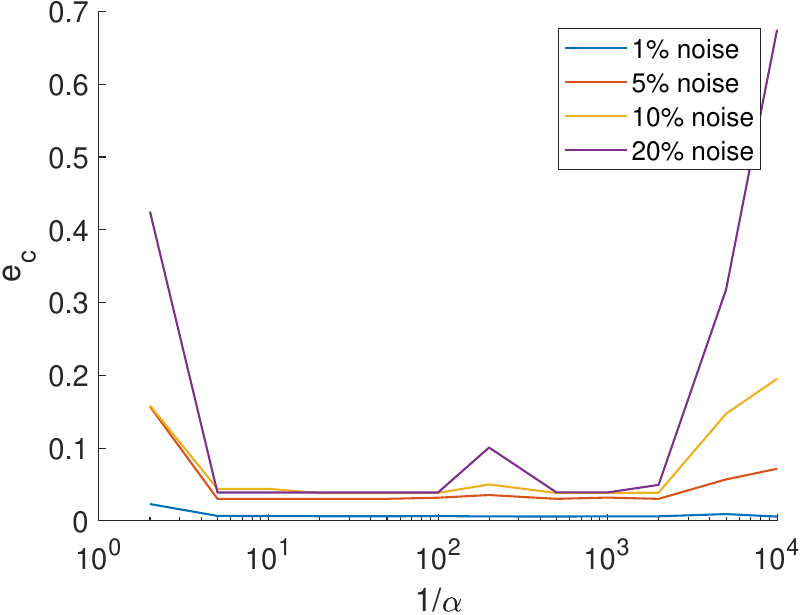}&
		\includegraphics[width=2.3in]{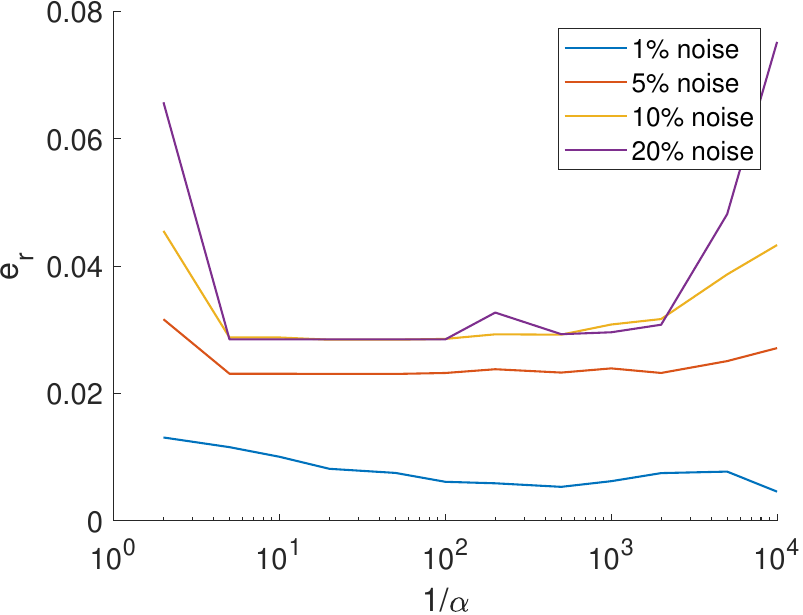}
	\end{tabular}
	\caption{Robustness of SC to the choice of  $\alpha$ for the recovery of the transport equation \eqref{eq.transport1d}.  (a) and (b) display $e_c$ and $e_r$ versus $1/\alpha$ respectively, with $1\%$ (Blue), $5\%$ (Red), $10\%$ (Orange), $20\%$ (Purple) noise. Each experiment is repeated 50 times, and the errors are averaged. We observe that SC is not sensitive to $\alpha$, and there is a wide range of values for $\alpha$ that give rise to a small error. }\label{fig.transport1d.alpha}
\end{figure}

In Figure \ref{fig.transport1d.alpha}, we explore the robustness of SC with respect to the choice of $\alpha$.  We present $e_c$ and $e_r$ versus $1/\alpha$ in (a) and (b) respectively, with $1\%, 5\%,10\%,20\%$ noise. Each experiment is repeated 50 times, and the error is averaged. The result shows that SC, in this case, is not sensitive to $\alpha$, and there are wide range choices of  $\alpha$ that give rise to a small error.

\begin{figure}[h!]
	\centering
	\begin{tabular}{ccc}
		(a) & (b)  & (c)
		\\
		\includegraphics[width=1.5in]{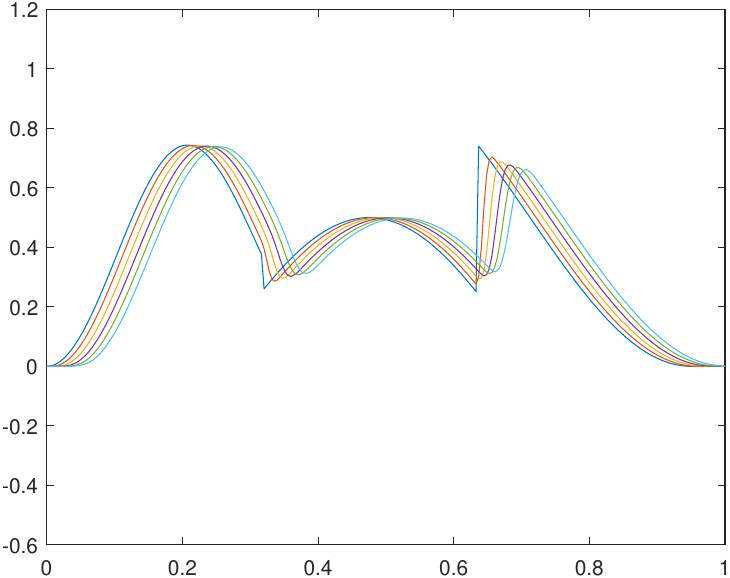}&
		\includegraphics[width=1.5in]{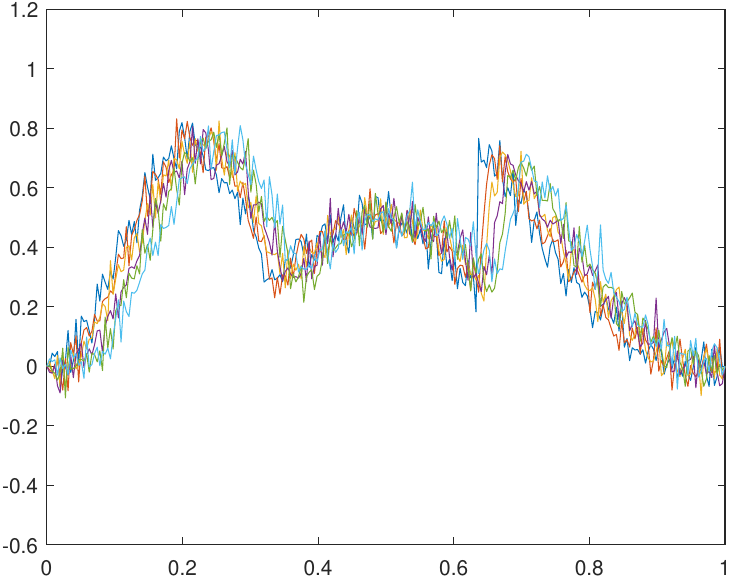}&
        \includegraphics[width=1.5in]{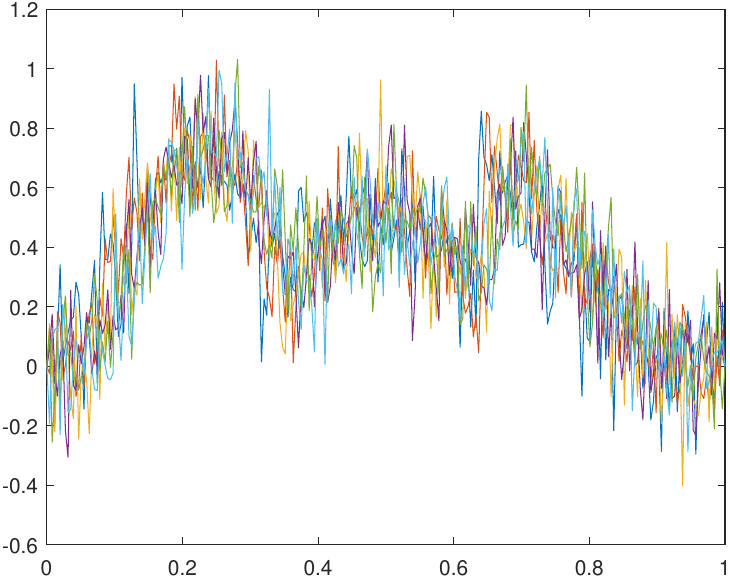}
	\end{tabular}
	\caption{Clean and noisy data of the transport equation \eqref{eq.transport1d} with the discontinuous initial condition (\ref{eq.transport.dis}).  (a) Clean data. (b) Noisy data with $10\%$ noise. (c) Noisy data with $30\%$ noise.
	}\label{fig.transport1d.dis}
\end{figure}

We next test ST and SC on data generated from the transport equation with a discontinuous initial condition. We set the initial condition as
\begin{equation}
u (x,0)=
\begin{cases}
\sin^2(2\pi x/(1-T))\cos(2\pi x/(1-T)), &\mbox{ for } 0\leq x< (1-T)/3,\\
-\cos^2(2\pi x/(1-T))+0.5, &\mbox{ for } (1-T)/3\leq x< 2(1-T)/3,\\
\sin^2(2\pi x/(1-T)), &\mbox{ for } 2(1-T)/3\leq x\leq (1-T),\\
0, &\mbox{ otherwise.}
\end{cases}
\label{eq.transport.dis}
\end{equation}
The clean data is generated by explicitly solving (\ref{eq.transport1d}) with $\delta x =\Delta x=1/256,\delta t =\Delta t=10^{-3}$ and $T=0.05$. After adding i.i.d. Gaussian noise, we have the noisy data. We show the clean data and the noisy data in Figure \ref{fig.transport1d.dis}. The identification results are shown in Table \ref{tab.transport1d.dis}. Even with the existence of discontinuities, ST and SC are stable and can identify the correct PDE with up to 30\% noise.

\begin{table}[h]
	\centering
		\caption{Identification of the transport equation \eqref{eq.transport1d} with the discontinuous initial condition (\ref{eq.transport.dis}) and different noise levels.  In the noise-free case, applying SDD does not introduce strong bias. The identification results (second column) by ST and SC are stable even with 30\% noise.  Here $w=20$ for ST, and $\alpha=1/200$ for SC.}\label{tab.transport1d.dis}
	\begin{tabular}{r|c|c|c}
		\hline
		\textbf{Method}& \textbf{$0\%$ noise without SDD } & $e_c$ & $e_r$\\\hline
		ST&$u_t=-1.0091u_x+9.65\times10^{-4}u_{xx}$ & $1.01\times 10^{-2}$ & $1.64\times 10^{-1}$\\\hline
		SC &$u_t=-1.0511u_x$ & $5.11\times 10^{-2}$ & $4.43\times 10^{-2}$\\\hline
		& \textbf{$0\%$ noise with SDD}&$e_c$&$e_r$\\\hline
		ST, SC&$u_t=-1.0274u_x$ & $2.74\times 10^{-2}$ & $1.95\times 10^{-2}$\\\hline\hline
		& \textbf{ $10\%$ noise }& $e_c$ & $e_r$\\\hline
		ST, SC&$u_t=-0.9913u_x$& $8.72\times 10^{-3}$ & $5.90\times 10^{-3}$\\\hline\hline
		& \textbf{ $30\%$ noise }& $e_c$ & $e_r$\\\hline
		ST, SC &$u_t=-0.9239u_x$& $7.61\times 10^{-2}$ & $5.36\times 10^{-2}$\\ \hline
	\end{tabular}
\end{table}
\subsection{Burgers' Equation}
In the second example, we test our methods on the Burgers' equation, which is a first-order nonlinear PDE:
\begin{equation}
u_t=-uu_x\;
\label{eq.burger1d}
\end{equation}
for $0<t\leq T$.
We use the initial condition
\begin{align}
  u(x,0)=\sin(4\pi x)\cos(\pi x)
  \label{eq.burgers.initial}
\end{align}
and zero Dirichlet boundary condition. Our data is generated by solving (\ref{eq.burger1d}) with $\delta x =\Delta x=1/256,\delta t =\Delta t=10^{-3}$ and $T=0.05$.
\begin{table}[h]
		\caption{Identification of the Burgers' equation \eqref{eq.burger1d} with initial condition (\ref{eq.burgers.initial}) and different noise levels. The identification results (second column) by ST and SC are  good with small $e_c$ and $e_r$ for a noise level up to $40\%$. Here $w =20$ for ST, and $\alpha=1/500$ for SC.
	}\label{tab.burger1d}
	\centering
	\begin{tabular}{r|c|c|c}
		\hline
		\textbf{Method}& \textbf{$0\%$ noise without SDD} & $e_c$ & $e_r$\\\hline
		ST &$u_t=-1.0023uu_{x}-2.38\times10^{-5}u_xu_{xx}$& $2.35\times 10^{-3}$ & $5.07\times 10^{-3}$\\\hline
		SC  &$u_t=-0.9960uu_{x}$ & $4.01\times 10^{-3}$ & $2.58\times 10^{-3}$\\\hline
		& \textbf{$0\%$ noise with SDD}& $e_c$ & $e_r$\\\hline
		ST &$u_t=-1.0079uu_{x}-0.0001u_xu_{xx}$& $7.97\times 10^{-3}$ & $1.43\times 10^{-2}$\\\hline
		SC &$u_t=-0.9888uu_{x}$ & $1.12\times 10^{-2}$ & $7.20\times 10^{-3}$\\\hline\hline
		& \textbf{$10\%$ noise}& $e_c$ & $e_r$\\\hline
		ST, SC&$u_t=-1.0246uu_{x}$& $2.46\times 10^{-2}$ & $1.52\times 10^{-2}$\\\hline\hline
		& \textbf{$40\%$ noise}& $e_c$ & $e_r$\\\hline
		ST, SC&$u_t=-0.7366uu_{x}$& $2.63\times 10^{-1}$ & $1.64\times 10^{-1}$\\\hline
	\end{tabular}
\end{table}

Table \ref{tab.burger1d} shows the results of ST(20) and SC(1/500) with various noise levels.
With clean data, ST identifies an additional term, but its coefficient is very small, and the corresponding $e_c$ and $e_r$ are small. SC works very well on clean data. With $10\%$ and $40\%$ noise, both methods identify the same PDE with small $e_c$ and $e_r$.

\begin{figure}[th!]
	\centering
	\begin{tabular}{ccc}
		(a)& (b) &(c) \\
		\hspace{-0.5cm}
		\includegraphics[width=0.33\textwidth]{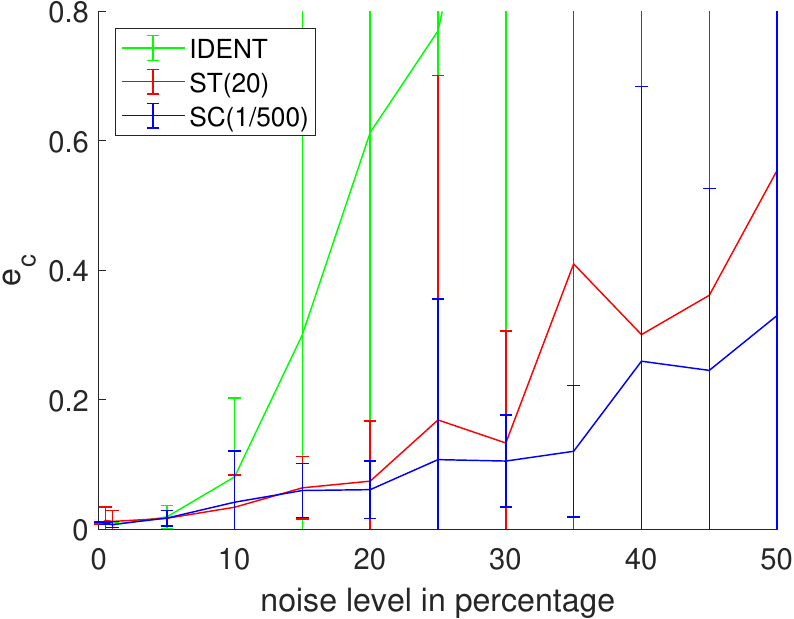}\hspace{-0.2cm} &
		\hspace{-0.2cm} \includegraphics[width=0.33\textwidth]{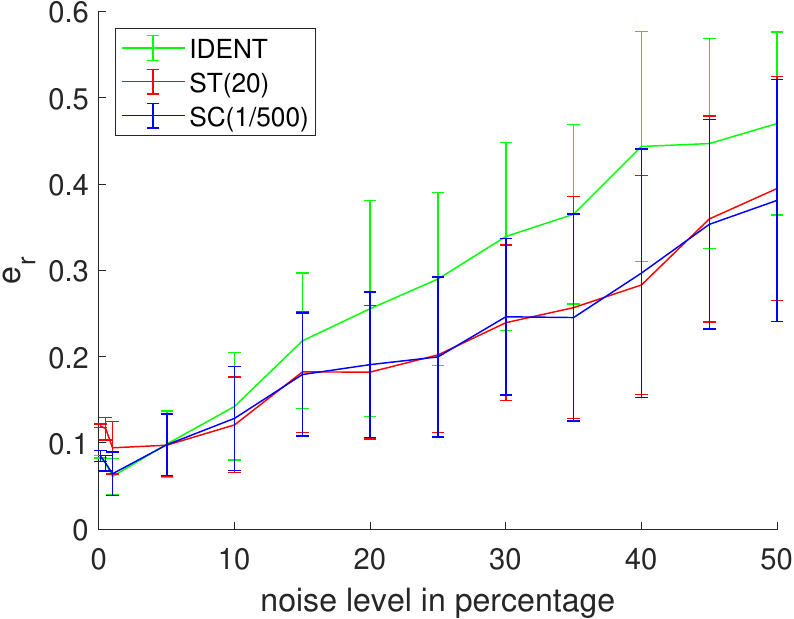} \hspace{-0.2cm} &
       \hspace{-0.2cm} \includegraphics[width=0.33\textwidth]{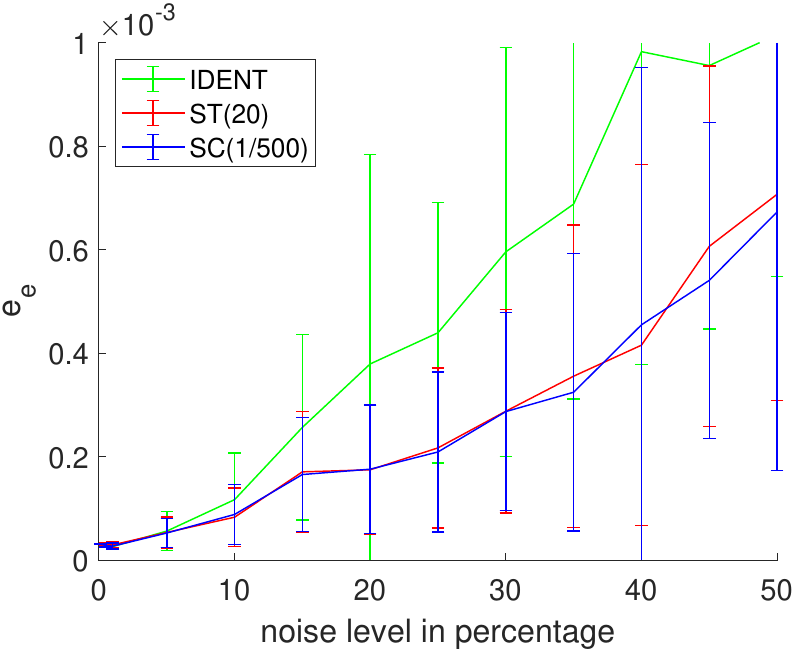} \hspace{-0.5cm}
	\end{tabular}
	\caption{The average error $e_c,e_r$ and $e_e$ over $50$ experiments for the Burgers' equation \eqref{eq.burger1d} with respect to various noise levels, where the initial condition is (\ref{eq.burgers.initial}). (a) The curve represents the average $e_c$ for IDENT~\cite{kang2019ident} (Geeen), ST (Red) and SC (Blue), and the  standard deviations are represented by vertical bars. (b) The average and variation of $e_r$ for IDENT (Geeen), ST (Red) and SC (Blue). (c) The average and variation of $e_e$ for IDENT (Geeen), ST (Red) and SC (Blue).
		The $e_c,e_r$ and $e_e$ of ST and SC are much smaller than those of IDENT. }\label{fig.burger1d.cre}
\end{figure}

Figure \ref{fig.burger1d.cre} shows how $e_c,e_r$ and $e_e$ change when the noise level varies. Each experiment is repeated 50 times and the errors are averaged.  We test IDENT, ST(20) and SC(1/500). The results in Figure \ref{fig.burger1d.cre} show that ST and SC perform better than IDENT.

In Table \ref{tab.burger1d.compare}, we compare SC, ST from this paper with IDENT in  \cite{kang2019ident}, the methods proposed in \cite{schaeffer2017learning} and \cite{rudy2017data}. The method from \cite{schaeffer2017learning} uses the spectral method to compute the spatial derivatives, which requires periodic boundary conditions. For a fair comparison, we use the initial condition
\begin{align}
u(x,0)=\sin(4\pi x)\cos(2\pi x)
\label{eq.burgers.initial.period}
\end{align}
and the periodic boundary condition (in which the boundary values are always 0).
Our data is generated by solving (\ref{eq.burger1d}) with $\delta x=\Delta x=1/256, \delta t=\Delta t=10^{-3}$ and $T=0.05$. We set $w =20$ for ST, and $\alpha=1/500$ for SC. For IDENT, we use SDD to denoise the data and to compute the partial derivatives, which improves the original IDENT in \cite{kang2019ident}. For the method in \cite{schaeffer2017learning}, we use the denoising method specified in  \cite[Example 3.9]{schaeffer2017learning}. The identification results are shown in Table \ref{tab.burger1d.compare}.  
Table \ref{tab.burger1d.compare} shows that ST, SC and IDENT are more robust than the method in \cite{schaeffer2017learning} at various noise levels. The errors given by ST, SC, and IDENT are also smaller. The results by the method in \cite{rudy2017data} are similar to those of ST and SC when the noise level is low. For large level of noise, for example $40\%$, ST and SC are more robust than the method in \cite{rudy2017data}. ST and SC can still identify the correct PDE with $40\%$ noise.

\begin{table}[h]

		\caption{Comparison of ST, SC with IDENT in  \cite{kang2019ident} and the methods in \cite{schaeffer2017learning} and \cite{rudy2017data} for the identification of the Burgers' equation \eqref{eq.burger1d} with the initial condition (\ref{eq.burgers.initial.period}), and various noise levels. 
		In this table, we only include the reconstructed terms with the coefficient magnitudes above $10^{-2}$. ST and SC are very stable compared to IDENT and the methods in \cite{schaeffer2017learning} and \cite{rudy2017data}. The coefficient error $e_c$ \eqref{eq.cre} and the time evolution error $e_e$ \eqref{eqtimeerror} are presented. With large noise, the errors given by ST, SC are smaller than the errors by other methods.
	}\label{tab.burger1d.compare}
	\centering
	\begin{tabular}{r|c|c|c}
		\hline
		\textbf{Method}& \textbf{$0\%$ noise} & $e_c$ &  $e_e$\\\hline
        \cite{schaeffer2017learning} & $u_t=-0.01u -0.95uu_x $& $6.49\times10^{-2}$  &$1.56\times10^{-4}$\\
        \hline
        \cite{rudy2017data} & 
        $\begin{aligned}
        u_t=-0.99uu_x
        \end{aligned}$ &
        $1.0\times 10^{-2}$ &  $3.46\times 10^{-5}$\\
        \hline
		ST, SC, IDENT &$u_t=-0.97uu_{x}$& $2.75\times 10^{-2}$ &  $8.01\times10^{-5}$\\\hline\hline
		& \textbf{$1\%$ noise}& $e_c$ & $e_2$\\\hline
        \cite{schaeffer2017learning} &$\begin{aligned}
        u_t=&-0.14u+0.01u^2\\
        &-0.89uu_x
        \end{aligned}
        $& $2.82\times 10^{-1}$ &$3.42\times10^{-4}$\\
        \hline
        \cite{rudy2017data} & 
        $\begin{aligned}
        u_t=-0.99uu_x
        \end{aligned}$ &
        $1.0\times 10^{-2}$ &  $3.46\times 10^{-5}$\\
        \hline
		ST, SC, IDENT &$u_t=-0.98uu_{x}$& $1.87\times 10^{-2}$ & $5.50\times 10^{-5}$\\\hline\hline
		& \textbf{$10\%$ noise}& $e_c$ & $e_e$\\\hline
        \cite{schaeffer2017learning} & $\begin{aligned}u_t=&-0.07+0.4u\\
        &\hspace{0.5cm}+0.44u^2-0.15uu_x\end{aligned}$ & $1.76$ & $2.59\times10^{-3}$\\
        \hline
        \cite{rudy2017data} & 
        $\begin{aligned}
        u_t=-0.94uu_x
        \end{aligned}$ &
        $6.0\times 10^{-2}$ &  $1.77\times 10^{-4}$\\
        \hline
        IDENT & $u_t=0.03u_x-1.00uu_x$&
        $3.0\times 10^{-2}$ &  $2.25\times 10^{-4}$\\
        \hline
		ST, SC &$u_t=-1.00uu_{x}$& $1.74\times 10^{-3}$ & $2.88\times 10^{-5}$\\\hline\hline
		& \textbf{$40\%$ noise}& $e_c$ & $e_e$\\\hline
        \cite{schaeffer2017learning} & $\begin{aligned}u_t=&-1+10.67u\\
        &\hspace{0.5cm}+1.84u^2-0.02uu_x\end{aligned}$ & $13.59$ & $7.16\times10^{-3}$\\
        \hline
        \cite{rudy2017data} & 
        $\begin{aligned}
        u_t=-0.93u-0.38uu_x
        \end{aligned}$ &
        1.56 &  $1.84\times 10^{-3}$\\
        \hline
		ST, SC, IDENT &$u_t=-1.02uu_{x}$& $2.39\times 10^{-2}$ & $8.27\times 10^{-5}$\\\hline
	\end{tabular}

\end{table}

\subsection{Burgers' Equation with Diffusion}

Our third example is the Burgers' equation with diffusion, which is a second order nonlinear PDE:
\begin{equation}
u_t=-uu_x+0.1u_{xx}\;.
\label{eq.burger1ddiff}
\end{equation}
We use the initial condition $u(x,0)=\sin(3\pi x)\cos(\pi x)$ and zero Dirichlet boundary condition. We first solve (\ref{eq.burger1ddiff}) with $\delta x=1/256, \delta t=10^{-5}$ and $T=0.05$. The given data is downsampled from the numerical solution such that $\Delta x=1/64$ and $\Delta t=10^{-4}$.

\begin{table}[h!]
	\centering
		\caption{Identification of the Burgers' equation with diffusion \eqref{eq.burger1ddiff} with different noise levels. The identification results (second column) by ST and SC are good with small $e_c$ and $e_r$ for a noise level up to $5\%$. Here $w =20$ for ST, and $\alpha=1/10$ for SC. }\label{tab.burger1ddiff}
	\begin{tabular}{r|c|c|c}
		\hline
		\textbf{Method}& \textbf{$0\%$ noise without SDD} & $e_c$ & $e_r$\\\hline
		ST, SC&$u_t=-1.0018uu_x+0.1001u_{xx}$& $1.67\times 10^{-3}$ & $8.14\times 10^{-4}$ \\\hline
		& \textbf{$0\%$ noise with SDD}& $e_c$ & $e_r$\\\hline
		ST, SC &$u_t=-0.9994uu_x+0.1009u_{xx}$& $1.36\times 10^{-3}$ & $7.68\times 10^{-3}$\\\hline\hline
		& \textbf{$1\%$ noise}& $e_c$ & $e_r$\\\hline
		ST, SC&$u_t=-0.9901uu_x+0.1013u_{xx}$& $1.02\times 10^{-2}$ & $1.19\times 10^{-2}$  \\\hline\hline
		& \textbf{$5\%$ noise}& $e_c$ & $e_r$\\\hline
		ST, SC&$u_t=-1.0170uu_x+0.0976u_{xx}$& $1.77\times 10^{-2}$ & $2.21\times 10^{-2}$ \\\hline
	\end{tabular}
\end{table}

Table \ref{tab.burger1ddiff} shows the results of ST(20) and SC(1/10) with various noise levels. With clean data, $1\%$ and $5\%$ noise, both methods identify the PDE with small $e_c$ and $e_r$.

\begin{figure}[h]
	\centering
	\begin{tabular}{ccc}
		(a)   &(b) & (c)\\
		\hspace{-0.7cm}
		\includegraphics[width=0.33\textwidth]{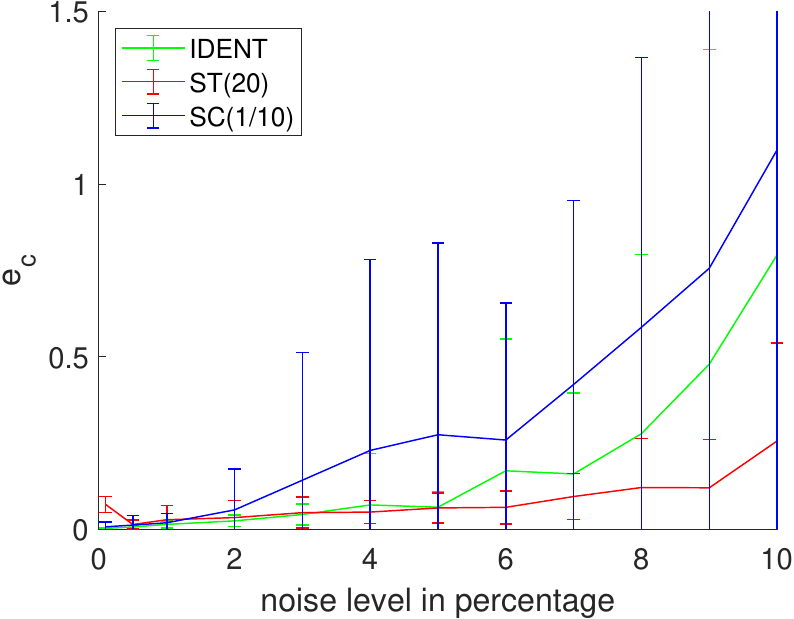}\hspace{-0.1cm}
		&\hspace{-0.1cm}\includegraphics[width=0.33\textwidth]{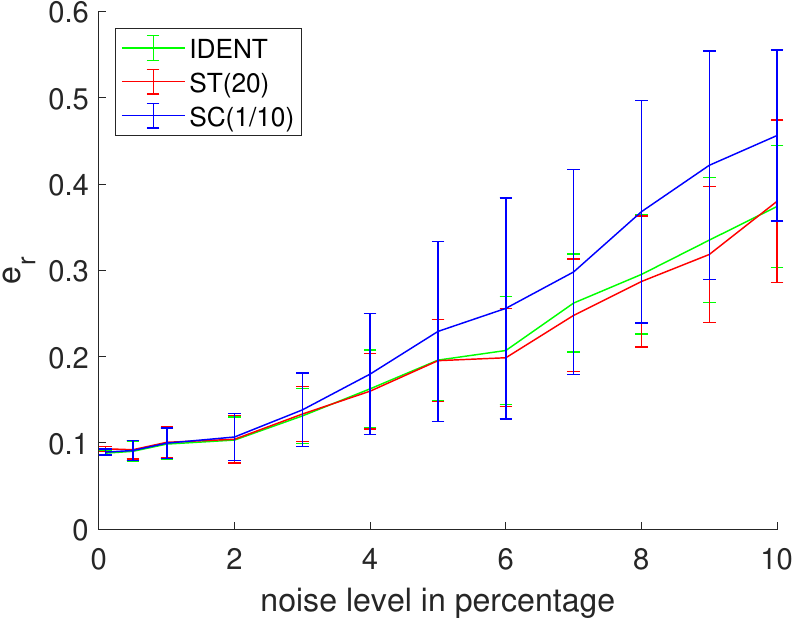}\hspace{-0.1cm}
        &\hspace{-0.1cm}\includegraphics[width=0.33\textwidth]{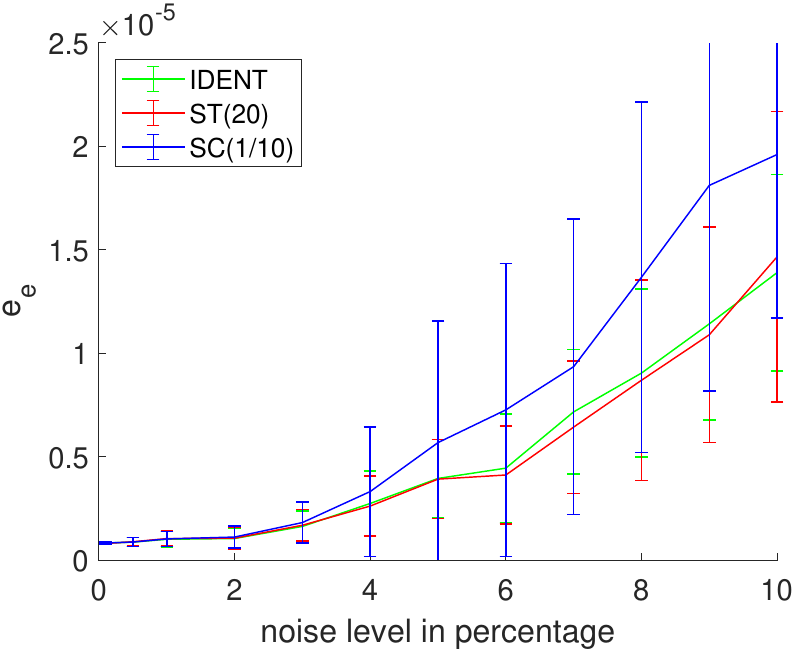}\hspace{-0.4cm}
	\end{tabular}
	\caption{The average error $e_c,e_r$ and $e_e$ over $50$ experiments of the Burgers' equation with diffusion \eqref{eq.burger1ddiff} with respect to various noise levels.  (a) The curve represents the average $e_c$ for IDENT~\cite{kang2019ident} (Geeen), ST (Red) and SC (Blue), and the standard deviations are represented by vertical bars. (b) The average and variation of $e_r$ for IDENT (Green), ST (Red) and SC (Blue). (b) The average and variation of $e_e$ for IDENT (Green), ST (Red) and SC (Blue). Among the three methods, ST gives the best result. }\label{fig.burger1ddiff.cre}
\end{figure}

Figure \ref{fig.burger1ddiff.cre} shows how $e_c,e_r$ and $e_e$ change when the noise level varies from $0.1\%$ to $10\%$. Each experiment is repeated 50 times, and the error is averaged.  We test IDENT, ST(20), and SC(1/10). Among the three methods, ST is the best. SC does not perform as well as ST and IDENT when the noise level is large. For high order PDEs, the high order derivatives are heavily contaminated by noise, even with SDD, which affects the accuracy of cross-validation. While ST and IDENT use time evolution, it is easier to pick correct features. In general, ST performs better than SC for high order PDEs when the given data contain heavy noise.

\begin{figure}[h!]
	\centering
	\begin{tabular}{cc}
		(a) & (b) \\
		\includegraphics[width=2.3in]{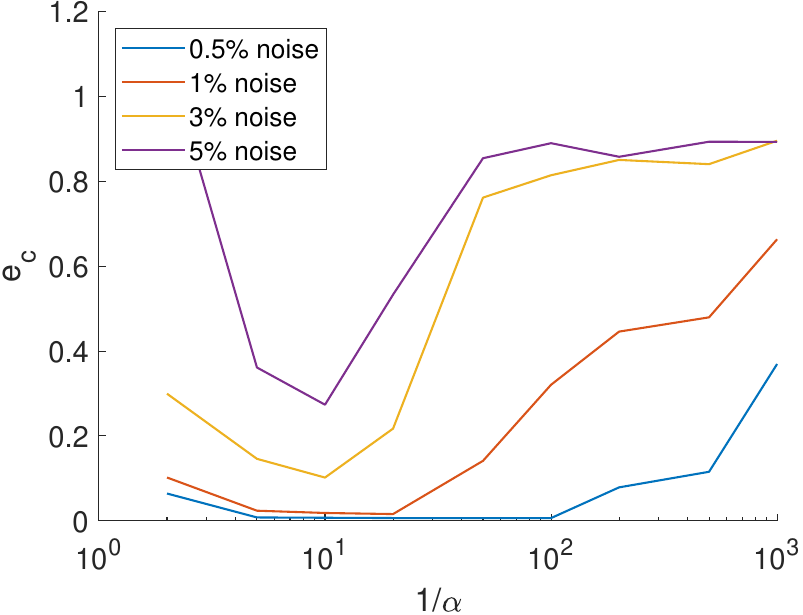}&
		\includegraphics[width=2.3in]{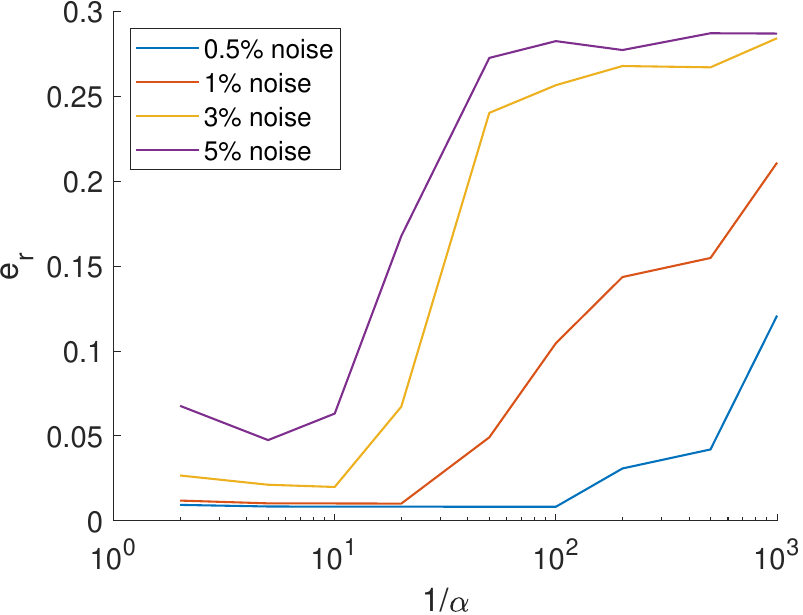}
	\end{tabular}
	\caption{Robustness of SC to the choice of  $\alpha$ for the recovery of the Burgers' equation with diffusion \eqref{eq.burger1ddiff}.  (a) and (b) display $e_c$ and $e_r$ versus $1/\alpha$ respectively, with $0.5\%$ (Blue), $1\%$ (Red), $3\%$ (Orange), $5\%$  (Purple) noise. Each experiment is repeated 50 times, and the errors are averaged. When the noise level is low, such as $0.5\%$ and $1\%$, there is a wide range of values for $\alpha$, which give a small error.  As the noise level increases, the range of the optimal $\alpha$ becomes narrow. }\label{fig.burger1DDiff.alpha}
\end{figure}

In Figure~\ref{fig.burger1DDiff.alpha}, we explore the effect of $\alpha$ in SC on the Burgers' equation with diffusion. Figure~\ref{fig.burger1DDiff.alpha} (a) and (b) show $e_c$ and $e_r$ versus $1/\alpha$ respectively, with $0.5\%$, $1\%$, $3\%$, and $5\%$ noise. When the noise level is low, such as $0.5\%$ and $1\%$, we have a wide range of good choices of $\alpha$ which gives rise to a smaller error. As the noise level increases, the range of the optimal $\alpha$ becomes narrow.

\subsection{The KdV Equation}

We test our algorithms on the KdV equation 
\begin{align}
	u_t+6uu_x+u_{xxx}=0,
	\label{eq.kdv}
\end{align}
on the spatial domain $[-10,10]$ and the time domain $0\le t \le T$ with $T=0.4$.
We use the initial condition $u(x,0)=5\sech^2(1.2x) $ and zero Dirichlet boundary condition. The data is generated with $\delta x=\Delta x=0.1, \delta t=10^{-5}.$ Data are downsampled in the time domain with $ \Delta t=10^{-3}$. Our dictionary contains $1,u,u_x,u_{xx}$ and $u_{xxx}$ and their pairwise products. There are 15 terms in the dictionary. The identified PDE by ST and SC from clean data is shown in Table \ref{tab.kdv}. In this example $w=20, \widetilde{\Delta t}=\Delta t/100$ is used in ST and $\alpha=1/1000$ is used in SC. Our results show that both ST and SC can identify the correct PDE.

\begin{table}[h!]
	\centering
	\caption{Identification of the KdV equation \eqref{eq.kdv}. Both ST and SC can identify the correct PDE. 
	}
	\label{tab.kdv}
	\begin{tabular}{r|c|c|c}
		\hline
		\textbf{Method}& Identified PDE & $e_c$ & $e_r$\\\hline
		ST, SC&$u_t=-6.135uu_x-1.0580 u_{xxx}$& $2.77\times 10^{-2}$ & $1.21$ \\\hline
	\end{tabular}
\end{table}

\subsection{A Larger Dictionary}

The examples above involve a dictionary which consists of the leading terms in the Taylor expansion of the governing equation $f(u,\partial_{\mathbf{x}}u,\partial_{\mathbf{x}}^2u)$. Our method is general and can be applied to other dictionaries.

We next test ST and SC on a larger dictionary, which includes $1,u,u_x, u_{xx}$ and $\sin(2\pi u), \cos(2\pi u)$ and their pairwise products. Since $\sin^2(2\pi u)+\cos^2(2\pi u)=1$, we exclude the term $\cos^2(2\pi u)$ to guarantee a set of linearly independent features. This dictionary contains 20 features. We consider the following PDE
\begin{align}
  u_t=u-0.1 u_x\sin(2\pi u)
  \label{eq.sin}
\end{align}
with the initial condition
$
  u(x,0)=0.8\sin(3\pi x)\cos(\pi x)
$
and zero Dirichlet boundary condition. The data are generated by solving (\ref{eq.sin}) with $\delta x=\Delta x=1/256, \delta t=\Delta t=4\times10^{-3}$ and $T=0.2$. The identified PDEs by ST and SC with various noise levels are shown in Table \ref{tab.sin}. On the clean data without SDD, ST identifies an additional term whose coefficient is very small. The corresponding $e_c$ and $e_r$ are very small. With up to $10\%$ noise, both ST and SC identify the correct PDE with a small $e_c$ and $e_r$.
\begin{table}[h!]
	\centering
		\caption{Identification of the \eqref{eq.sin} with different noise levels. The results (second column) by ST and SC are good with small $e_c$ and $e_r$ for up to $5\%$ noise. Here $w =20$ for ST, and $\alpha=1/500$ for SC. }\label{tab.sin}
	\begin{tabular}{r|c|c|c}
		\hline
		\textbf{Method}& \textbf{$0\%$ noise without SDD} & $e_c$ & $e_r$\\\hline
        ST&$\begin{aligned}u_t=&0.9994u-0.0995\sin(2\pi u)u_x\\ &\hspace{0.5cm}-2.90\times10^{-5}\cos(2\pi u)u_{xx}\end{aligned}$& $1.01\times 10^{-3}$ & $1.73\times 10^{-3}$ \\\hline
		SC&$u_t=0.9987u-0.0992\sin(2\pi u)u_x$& $1.88\times 10^{-3}$ & $1.13\times 10^{-3}$ \\\hline
		& \textbf{$0\%$ noise with SDD}& $e_c$ & $e_r$\\\hline
		ST, SC &$u_t=0.9903u-0.0895\sin(2\pi u)u_x$& $1.83\times 10^{-2}$ & $1.49\times 10^{-2}$\\\hline\hline
		& \textbf{$5\%$ noise}& $e_c$ & $e_r$\\\hline
		ST, SC&$u_t=0.9909u-0.0887\sin(2\pi u)u_x$& $1.85\times 10^{-2}$ & $1.56\times 10^{-2}$  \\\hline\hline
		& \textbf{$10\%$ noise}& $e_c$ & $e_r$\\\hline
		ST, SC&$u_t=1.0646u-0.1026\sin(2\pi u)u_x$& $6.11\times 10^{-2}$ & $1.33\times 10^{-2}$ \\\hline
	\end{tabular}
\end{table}

\subsection{Two Dimensional PDEs}
We next apply our methods to identify PDEs in a two-dimensional space.
The PDEs are solved with $\delta x=\delta y=0.02$ and $\delta t=8\times 10^{-4}$. Data are downsampled from the numerical solution with $\Delta x=0.04$ and $\Delta t=8\times 10^{-3}$. We fix $w=10$ for ST and $\alpha=3/200$ for SC.

The identification of two-dimensional PDEs is more challenging and more sensitive to noise. There are more features in two dimensions, and the directional variation of the data adds complexity to the problem.
We will show that both ST and SC are robust against noise.

We first consider the following PDE:
\begin{align}
\begin{cases}
u_t=0.02u_{xx}-uu_y\;\text{for}~(x,y,t)\in[0,1]^2\times [0,0.1],
\\
u(x,y,0) = \sin^2(\frac{3\pi x}{0.9})\sin^2(\frac{2\pi x}{0.9})\;\text{when}~(x,y)\in[0,0.9]^2~\text{and}~0~\text{otherwise}.
\end{cases}	\;,
\label{eq.2d.1}
\end{align}
which has different dynamics along the $x$ and $y$ directions. Table~\ref{tab.2dPDE1} shows the identification results of ST(10) and SC(3/200) with noise level $0\%, 5\%$ and $10\%$. Both methods identify the same features with small $e_c$ and $e_r$.
\begin{table}[h]
	\centering
	\caption{Identification of the PDE~\eqref{eq.2d.1} with different noise levels. The results (second column) by ST and SC have small $e_c$ and $e_r$ for up to $10\%$ noise. Here $w =10$ for ST, and $\alpha=3/200$ for SC. }\label{tab.2dPDE1}
	\begin{tabular}{r|c|c|c}
		\hline
		\textbf{Method}& \textbf{$0\%$ noise} & $e_c$ & $e_r$\\\hline
		ST, SC &$u_t=0.0189u_{xx}-0.9525uu_y$ & $4.75\times 10^{-2}$ & $2.48\times10^{-2}$ \\\hline\hline
		& \textbf{$5\%$ noise}& $e_c$ & $e_r$\\\hline
		ST, SC &$u_t=0.0178u_{xx}-0.9362uu_y$& $8.43\times 10^{-2}$ & $7.45\times 10^{-2}$\\ \hline\hline
		& \textbf{$10\%$ noise}& $e_c$ & $e_r$\\\hline
		ST, SC&$u_t=0.0134u_{xx}-0.8674uu_y$&$1.33\times 10^{-1}$ & $1.79\times 10^{-1}$ \\
		\hline
	\end{tabular}
\end{table}

\subsection{Identifiability Based on the Given Data}
For the PDE identification, especially in high dimensions, the given data $ U $ plays an important role.  When the initial condition has sufficient variations in each dimension,  the correct PDE can be identified. Otherwise, there may be multiple PDEs which generate the same dynamics.
For example, we consider the following transport equation:
\begin{align}
\begin{cases}
u_t=-0.5u_x+0.5u_y,~(x,y)\in[0,1]\times[0,1],~t\in[0,0.1]\\
u(x,y,0) = f(x,y),~(x,y)\in[0,1]\times[0,1]
\end{cases}\;,  \label{2DPDEexp1}
\end{align}
where $f$ denotes the initial condition.

We first choose the initial condition $f(x,y)=\sin(2\pi x/0.9))^2\sin(2\pi y/0.9)^2$ for $(x,y)\in[0,0.9]\times [0,0.9]$ and $0$ otherwise. The noise-free data are generated with  $\delta x=\delta y = 0.02$ and $\delta t = 7\times10^{-4}$, and downsampled in space by a factor of $2$ and in time by a factor of $10$.   The identified PDE by SC(1/200) is
$$u_t = -0.5001u_x+0.4800u_y\;,$$
where the recovered coefficients are very close to the true coefficients. The same result is identified by using ST(20). 

We next choose $f(x,y)=\sin(2\pi x/0.9))^2$ for $(x,y)\in[0,0.9]\times \mathbb{R}$ and $0$ otherwise. Our methods  SC(1/200) and ST(20) both identify
\begin{equation}
u_t = -0.4992u_x\;.\label{2DPDEexp1_eq}
\end{equation}
With this initial condition, the PDE in \eqref{2DPDEexp1} has  the exact solution:
\begin{align*}
u(x,y,t) = \begin{cases}\sin(\frac{2\pi(x-0.5t)}{0.9})^2,\quad x\in[0.5t,0.9+0.5t],~(x,y)\in\mathbb{R}\times[0,1],~t\in[0,0.1]\\
0,\quad\text{Otherwise}
\end{cases}
\end{align*}
which also satisfies $u_t=-0.5u_x$.
The identified PDE in (\ref{2DPDEexp1_eq}) approximates this simpler equation.
Since the given data only vary along the $x$ direction, the columns in the feature matrix related to $y$, e.g., $u_y$, $u_xu_y$, and $u_{yy}$, are mostly $0$. This explains why our method identifies the PDE in (\ref{2DPDEexp1_eq}), instead of \eqref{2DPDEexp1}. 

In this problem, the original PDE can be identified if the initial condition has sufficient variations. The identifiability of a PDE for a given dictionary under sparsity constraints can be defined as follows: Suppose the original PDE is associated with the coefficient vector $\mathbf{c}_0$ with sparsity $S$. This PDE is identifiable if there is a unique coefficient vector with sparsity no more than $S$, such that the evolution of the PDE associated with this coefficient vector, starting from the given initial condition, matches the given data. We believe it is an open question to investigate the theoretical conditions under which the PDE  is identifiable. Roughly speaking, the PDE problem is identifiable if the PDE solution with a given initial condition gives rise to the feature matrix $F$, which has a small pairwise coherence, in the sense that any two columns of $F$ have a small correlation. We refer to \cite[Theorem 1]{kang2019ident} for an identifiability condition in Lasso.

 \subsection{SC comparison}
\label{subsecsccomparison}

\begin{figure}
\begin{center}
	\includegraphics[height=1.8in]{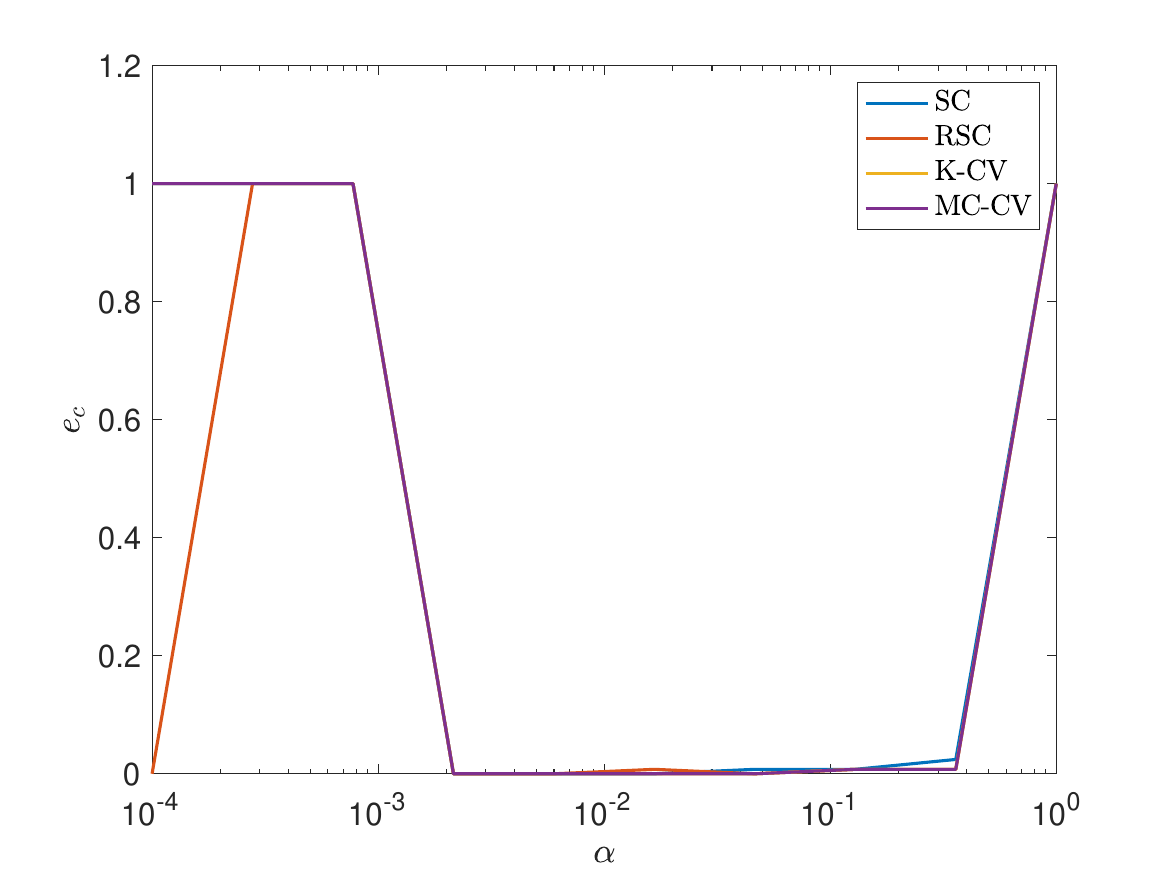}
	\caption{The coefficient errors $e_c$ for different sampling strategies for the identification of the PDE in (\ref{eq.random}) as the parameter $\alpha$ varies. This shows that in general, different sampling strategies in SC lead to similar identification results.}\label{fig_response}
\end{center}
\end{figure}

Our SC strategy is a two-fold strategy: In the first fold, we choose the first $\alpha$ fraction of the rows for training and the rest for testing. In the second fold, we choose the last alpha fraction of rows  for training and the rest for testing. Then we take the average of the two testing errors.  
We next compare the identification results using our current strategy, the random selection, $K$-fold cross validation and Monte-Carlo cross validation:

\begin{itemize}
	\item SC (our current strategy):  Without changing the time order of the data, for a fixed $0<\alpha<1$, we select the PDE by minimizing the average testing error of two types. (1) Head:  use the first $\alpha$ of the data for training, and the rest for testing. (2) Tail: use the last $\alpha$ of the data for training, and the rest for testing. 
	
	\item Random SC (RSC): Randomly permute the data in time, then the remaining procedure is the same as SC. 
	
	\item $K$-fold cross validation (K-CV): Randomly permute the data in time, then uniformly split the data into $K$ groups. The error for a candidate PDE is evaluated by taking the average of $K$ testing errors: for $k=1,2,\dots, K$, while the $k$-th group data is used for training, and the rest is for testing. 
	
	\item Monte-Carlo cross validation (MC-CV): Fix the number of simulation $N$ and a coefficient $0<\alpha<1$. For $n=1,\dots, N$, randomly permute the data and use the first $\alpha$ of the data for training, and the rest for testing.
\end{itemize}
We consider the following underlying PDE:
\begin{align} 
u_t = -0.5uu_x + 0.5uu_y
\label{eq.random}
\end{align}
with the initial condition $f(x,y)=\sin(2\pi(x+y))$ filtered by the Tukey window to comply with our zero-boundary requirement. 
We add $0.5\%$ noise to the data set. The methods above identify the same correct model, as shown in Table \ref{tab.random}.
\begin{table}
\begin{center}
	\begin{tabular}{c|c}
		\hline
		\multicolumn{2}{c}{With $0.5\%$ Noise}\\\hline
		\textbf{Sampling Strategy}& \textbf{Identified PDE}\\\hline
		SC: $\alpha=1/400$& $u_t=-0.5000uu_x+0.5002uu_y$\\\hline
		RSC: $\alpha=1/400$&$u_t=-0.5000uu_x+0.5002uu_y$\\\hline
		K-CV: $K=400$&$u_t=-0.5000uu_x+0.5002uu_y$\\\hline
		MC-CV: $N=100, \alpha=1/400$& $u_t=-0.5000uu_x+0.5002uu_y$\\\hline
	\end{tabular}
\caption{The PDE identification of (\ref{eq.random}) by SC with different sampling strategies.}\label{tab.random}
\end{center}
\end{table}
Moreover, the effective range for $\alpha$ (or equivalently $1/K$) for these sampling schemes are similar. This is demonstrated in Figure~\ref{fig_response}, where we vary $\alpha$ and record the coefficient errors $e_c$ of the identified PDEs, respectively.

\subsection{Choice of Smoother in SDD}\label{subsec::NoiseAnalysis}

\begin{center}
	\begin{figure}[h]
		\begin{tabular}{cllllll}
			& $u$ & $u_x$ & $u_{yy}$& $uu_x$\\
		 $0\%$& \includegraphics[height=0.8in]{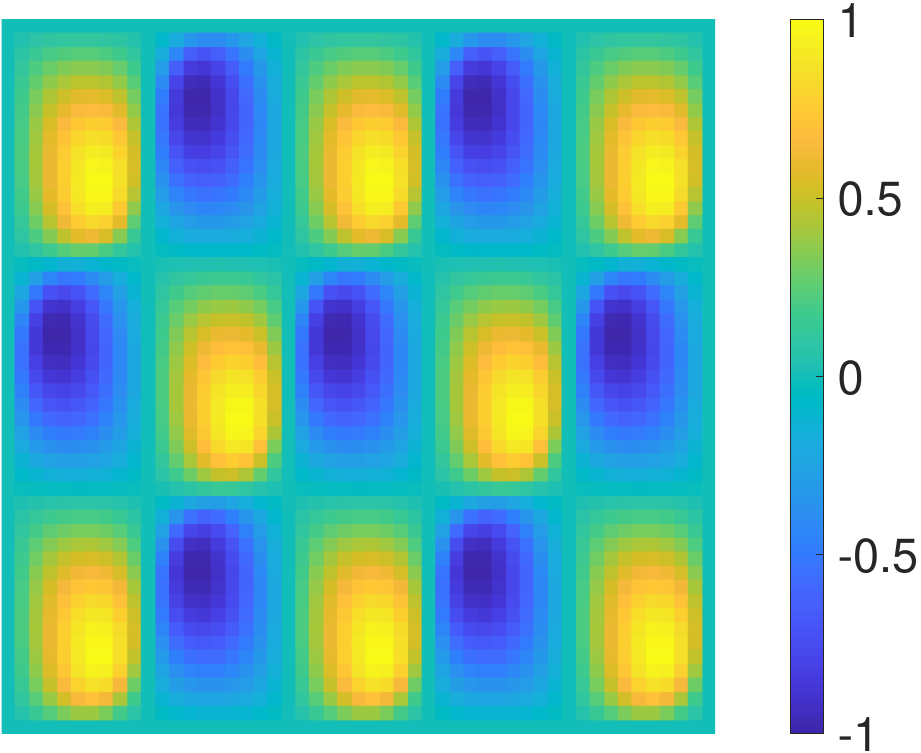}&
		\includegraphics[height=0.8in]{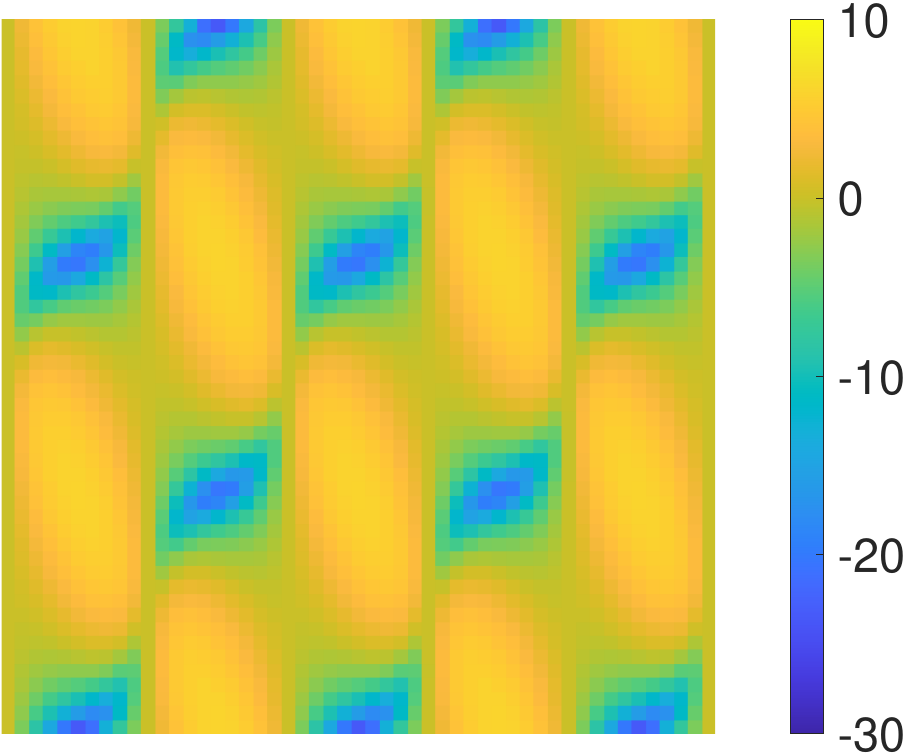}&
		  \includegraphics[height=0.8in]{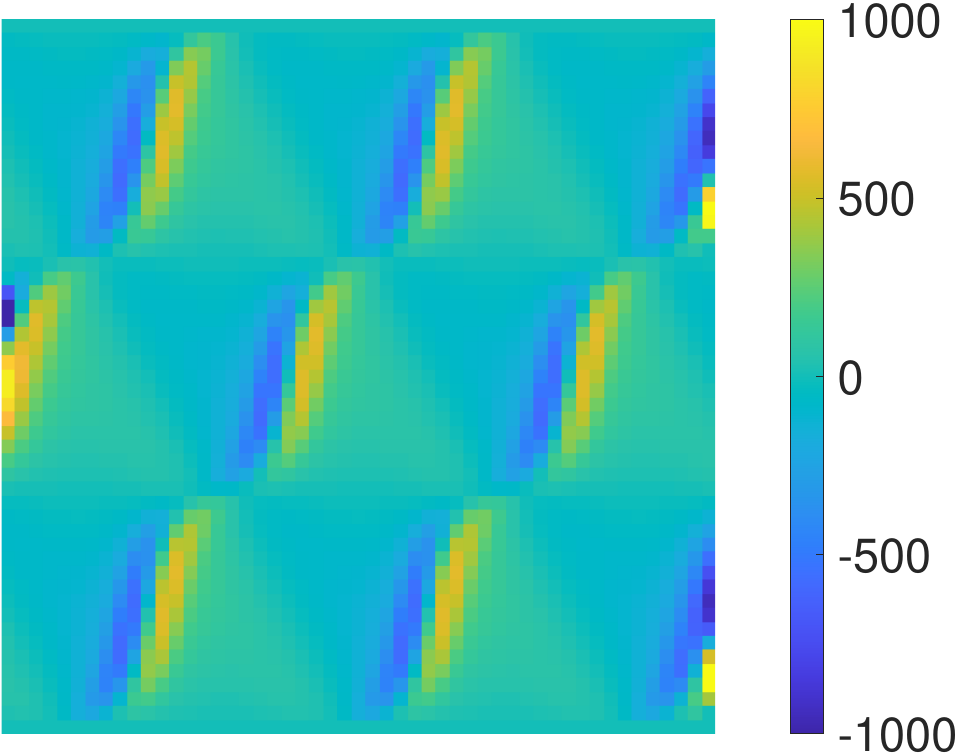}&
		  \includegraphics[height=0.8in]{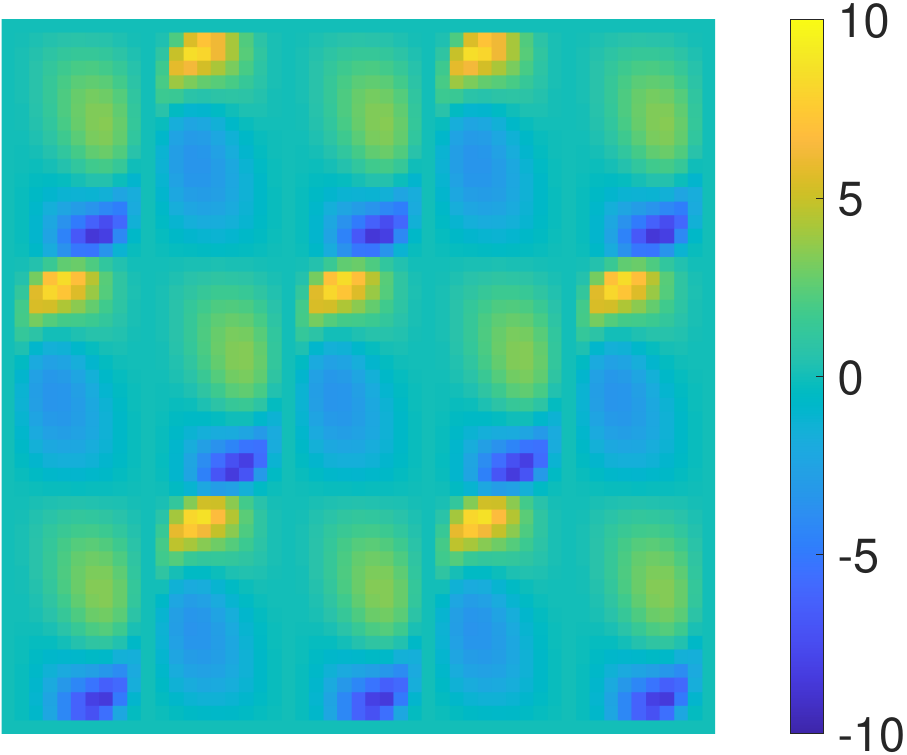}\\
		  $5\%$&\includegraphics[height=0.8in]{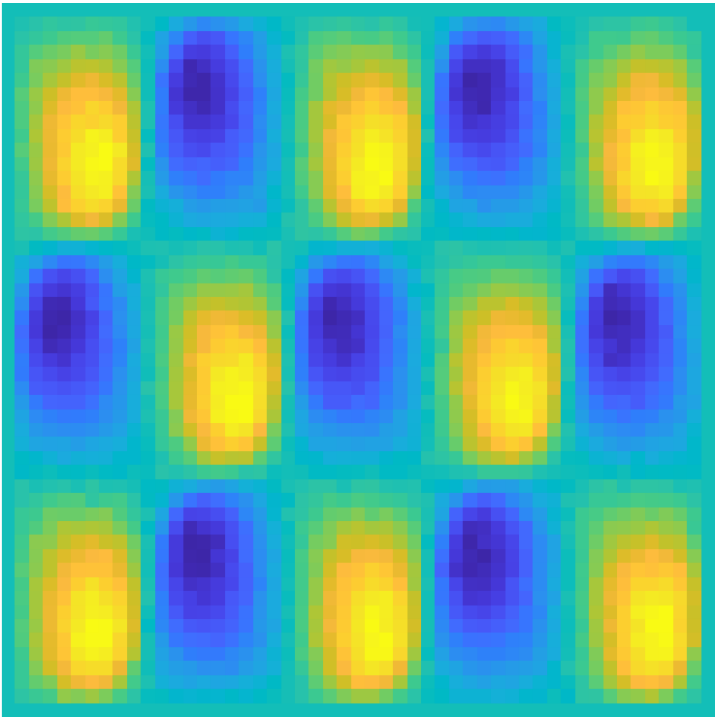}&
		  \includegraphics[height=0.8in]{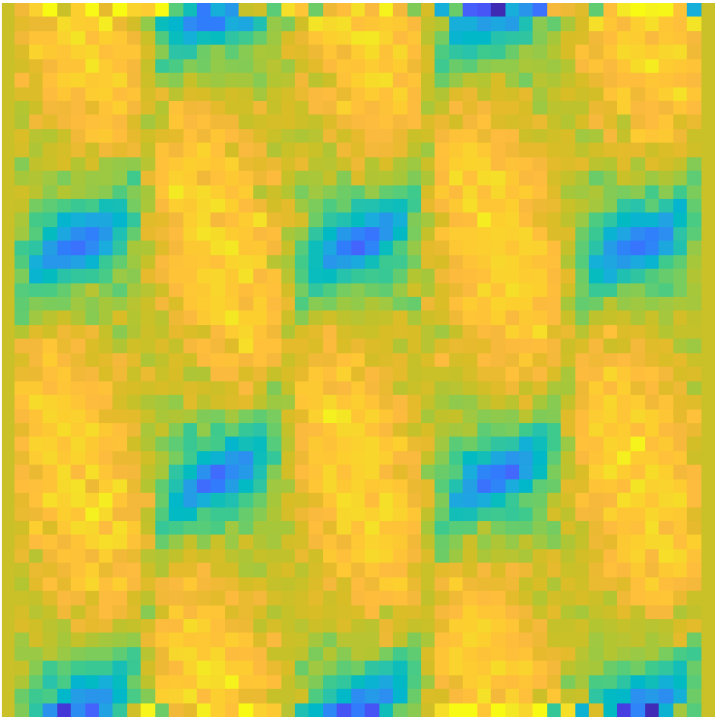}&
		  	\includegraphics[height=0.8in]{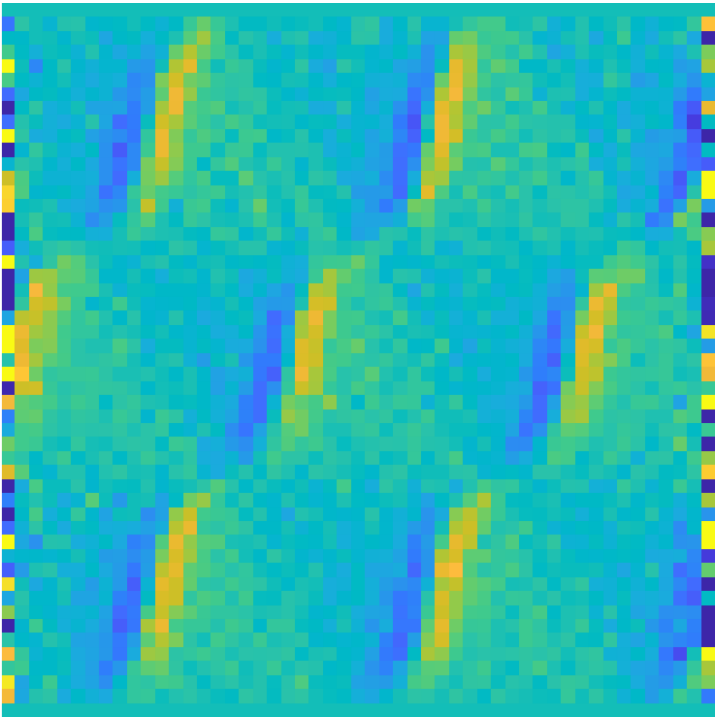}&
		  		\includegraphics[height=0.8in]{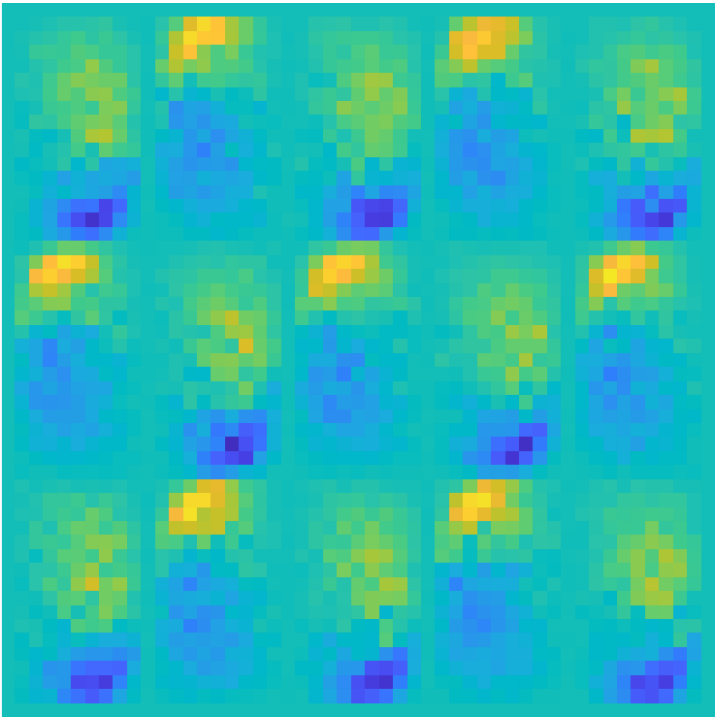}\\
		  	MA&\includegraphics[height=0.8in]{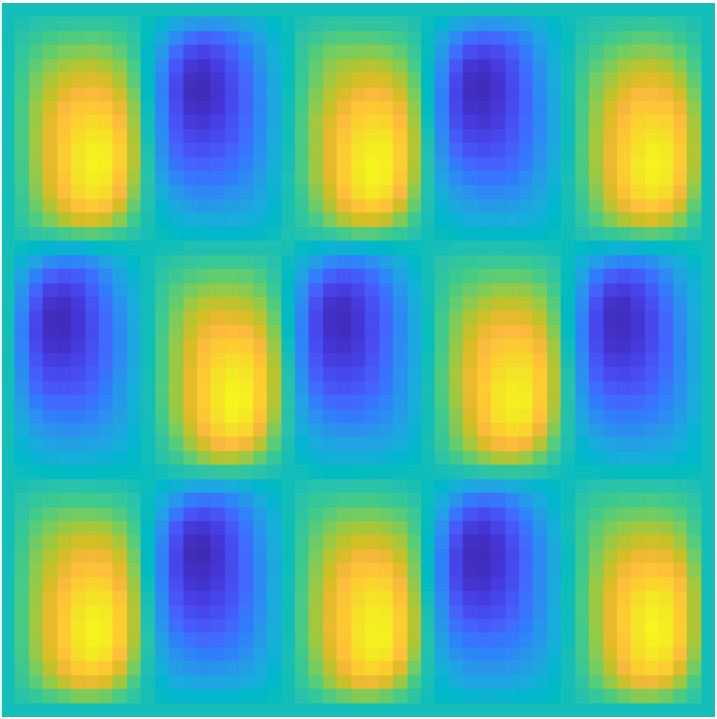} &
		  		\includegraphics[height=0.8in]{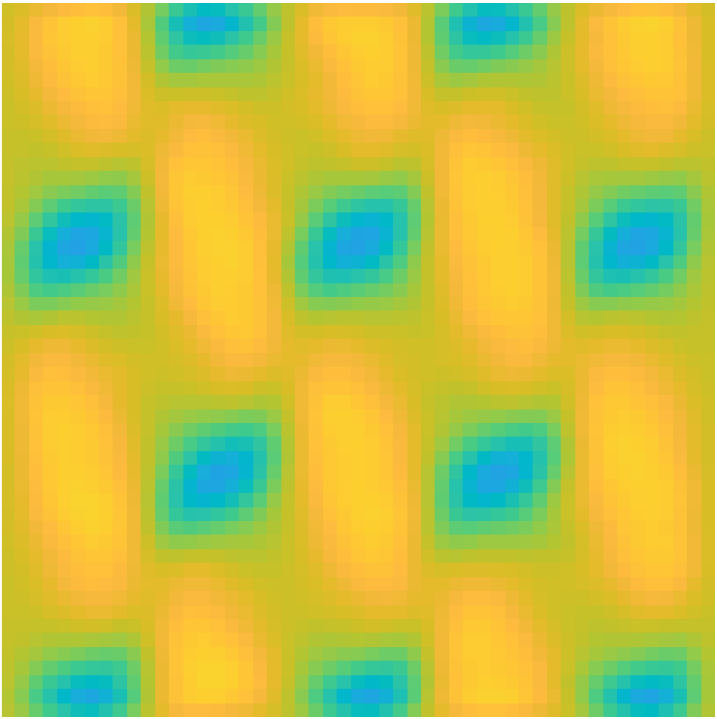} &
		  		\includegraphics[height=0.8in]{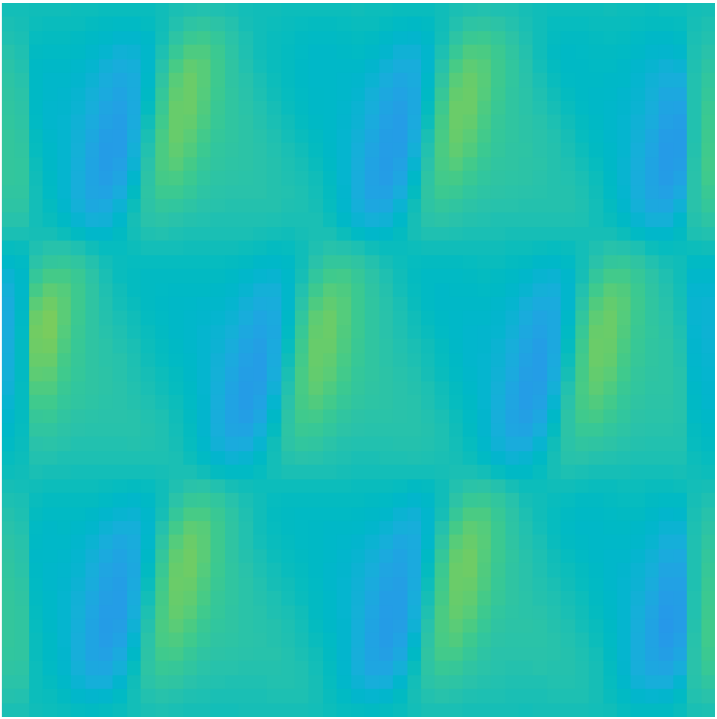} &
		  		\includegraphics[height=0.8in]{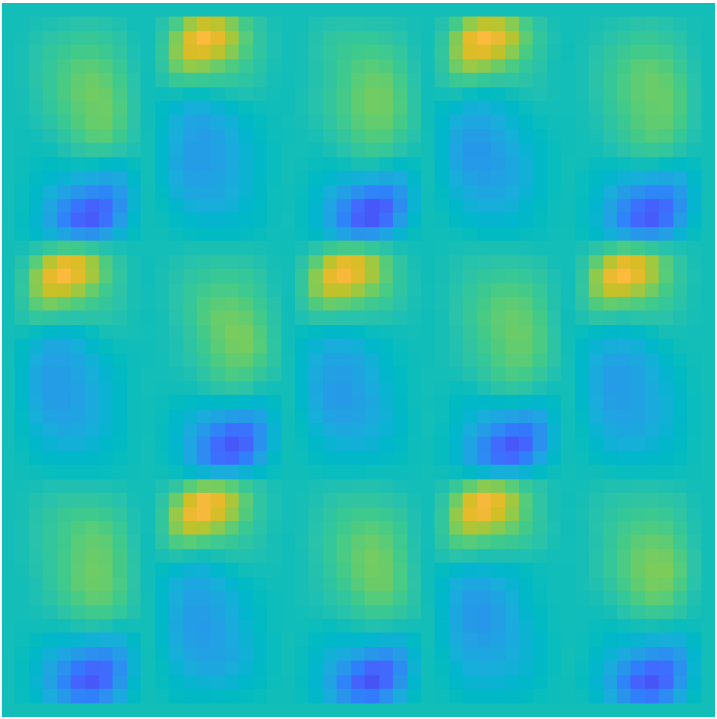}\\
		  	CS&\includegraphics[height=0.8in]{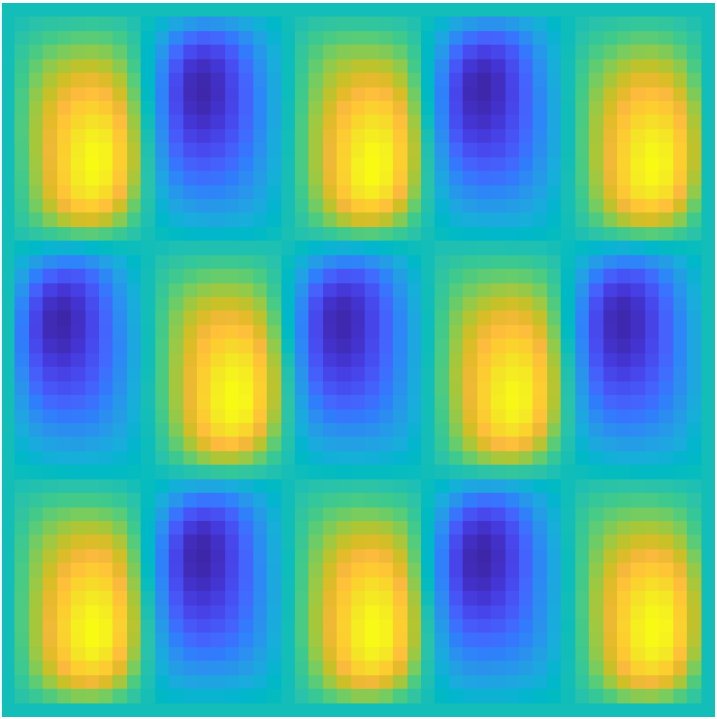} &
		  	\includegraphics[height=0.8in]{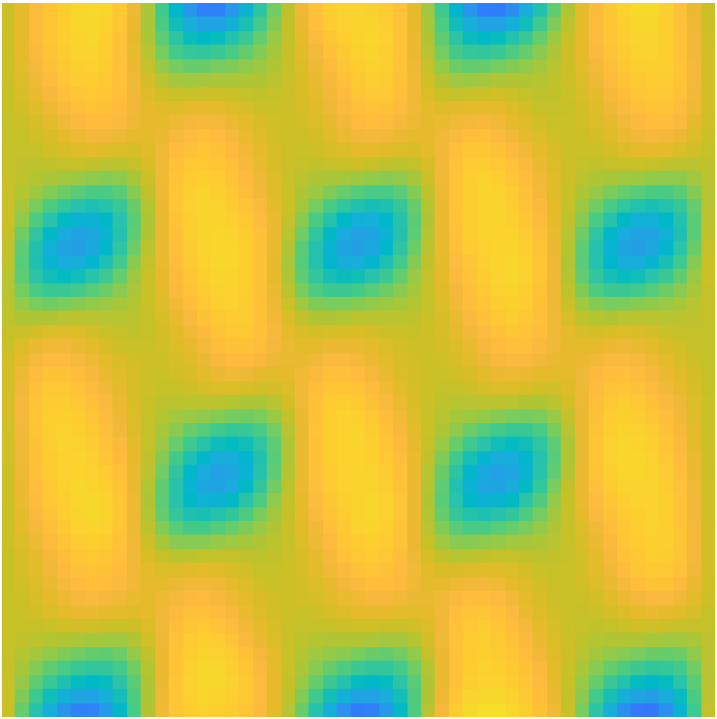} &
		  	\includegraphics[height=0.8in]{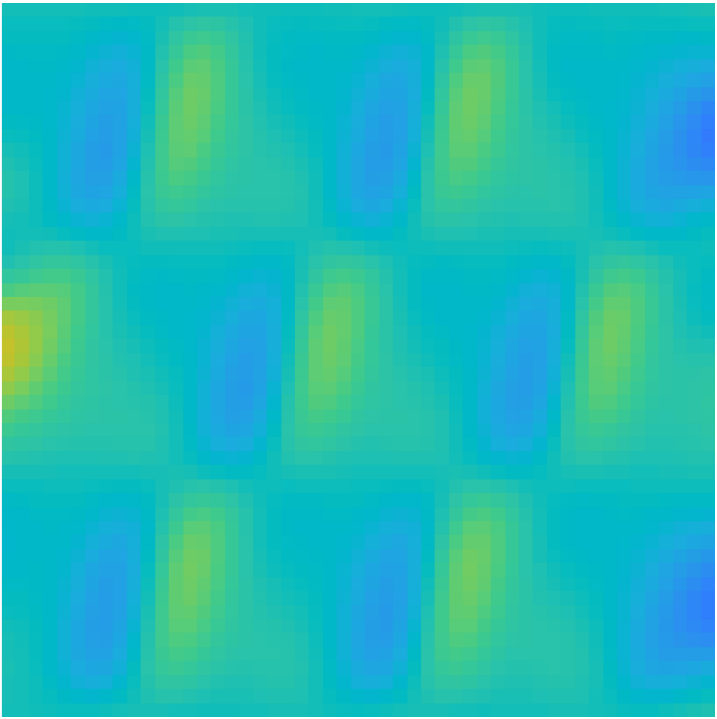} &
		  	\includegraphics[height=0.8in]{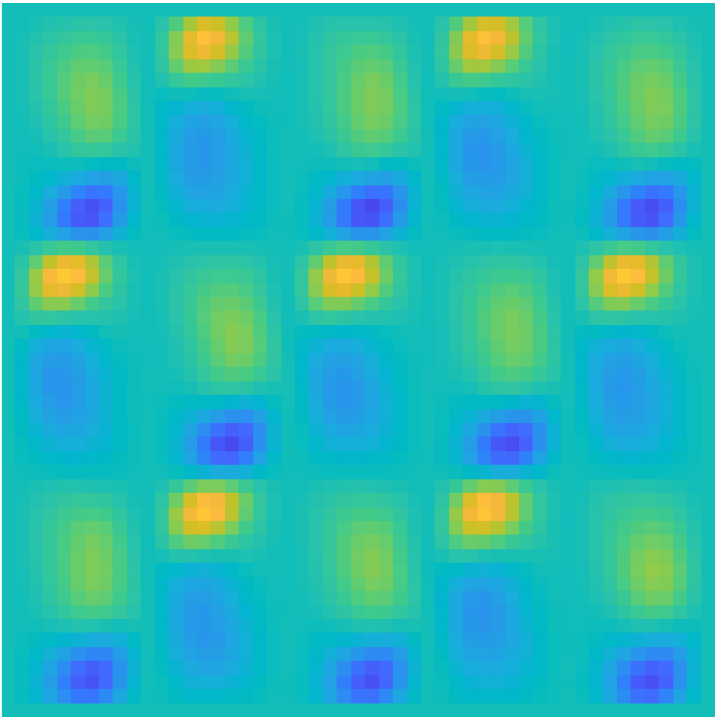} \\
		  	DF&\includegraphics[height=0.8in]{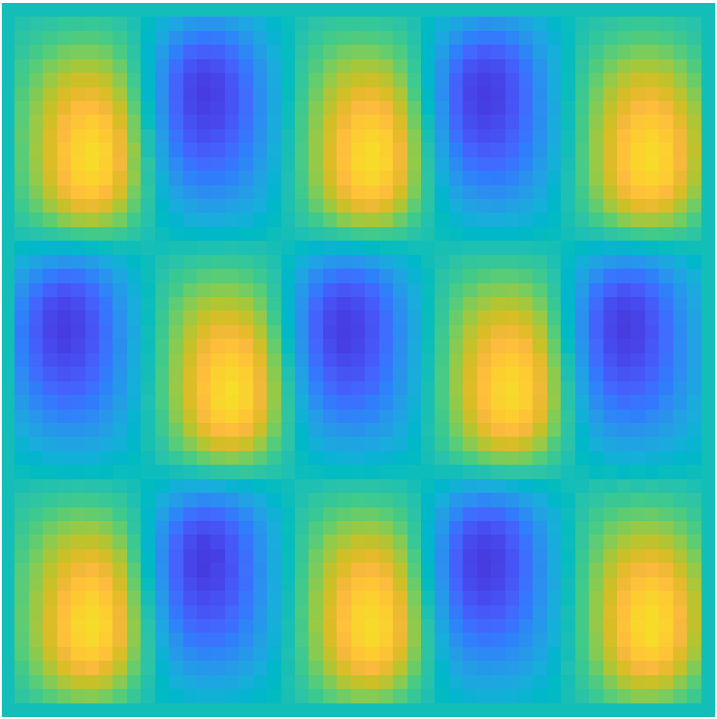} &
		  	\includegraphics[height=0.8in]{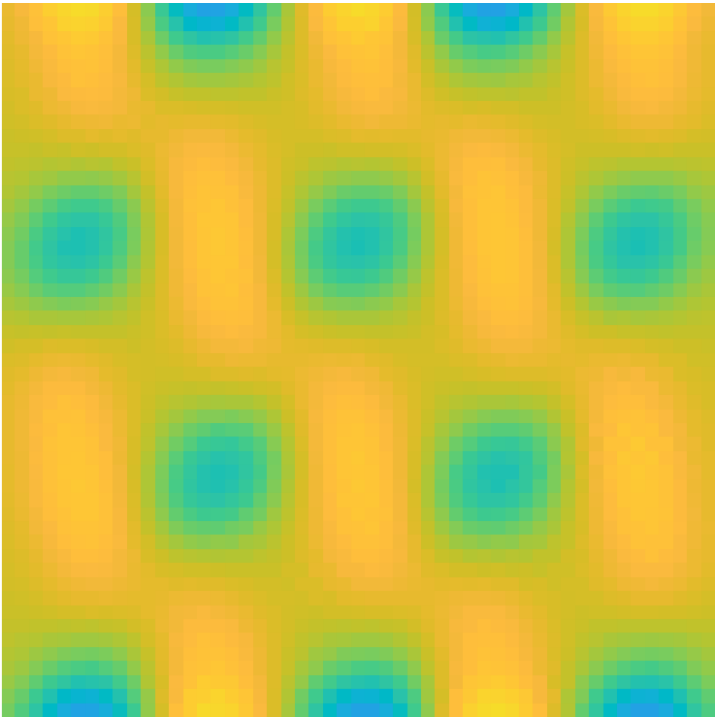} &
		  	\includegraphics[height=0.8in]{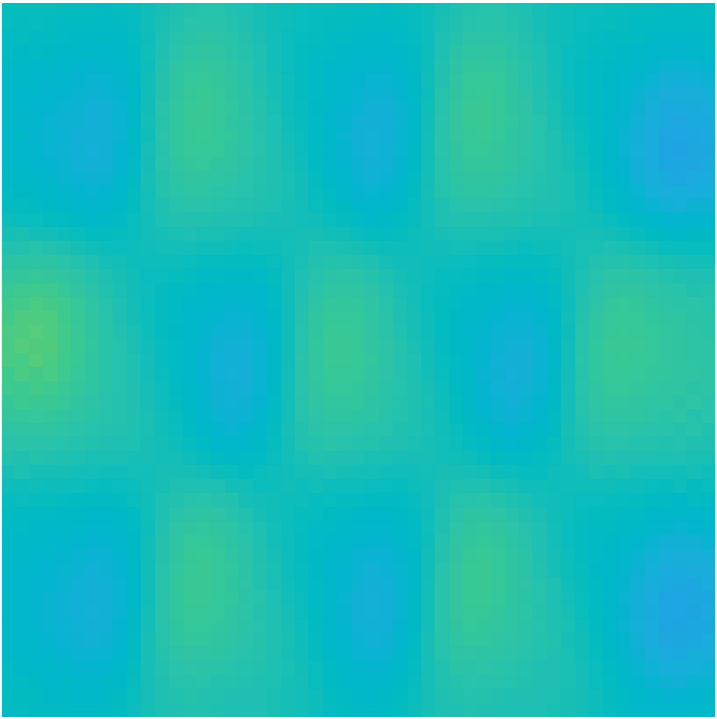} &
		  	\includegraphics[height=0.8in]{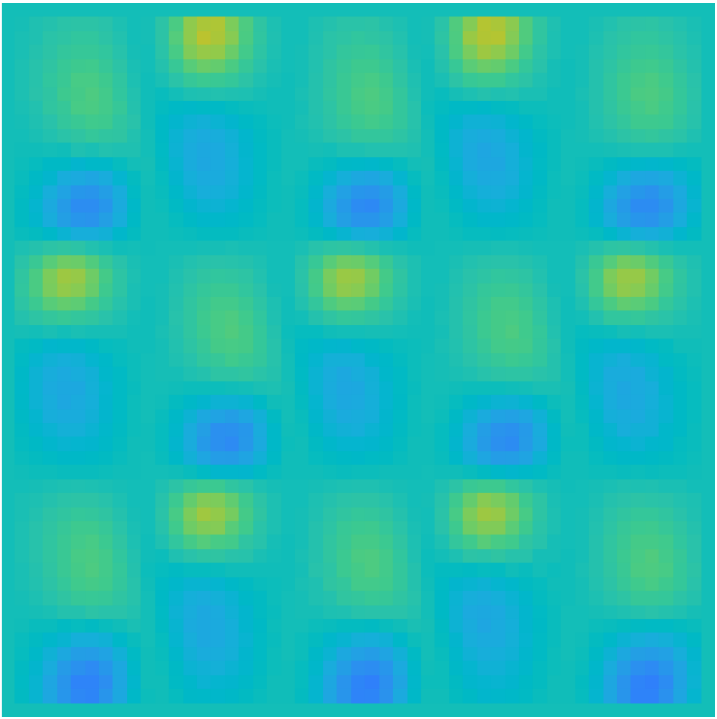} \\
		  	MLS&\includegraphics[height=0.8in]{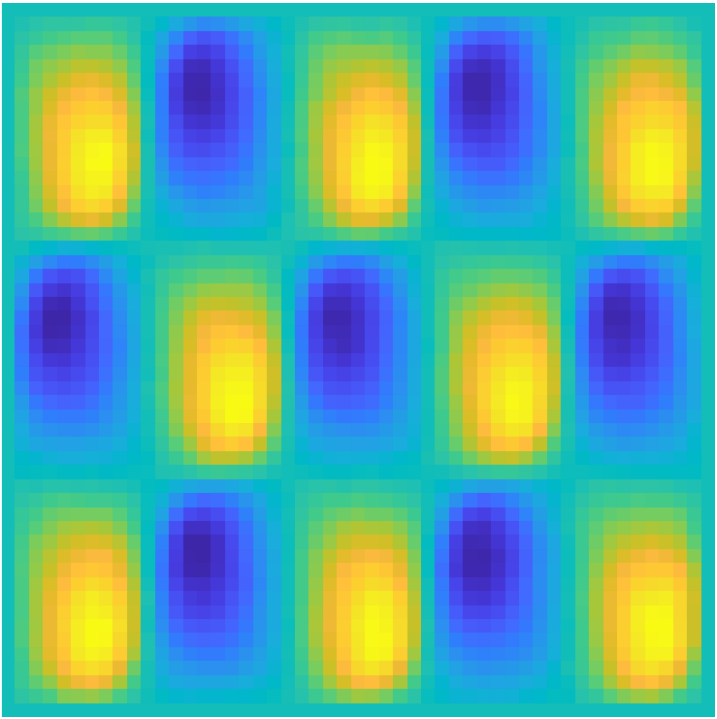} &
		  	\includegraphics[height=0.8in]{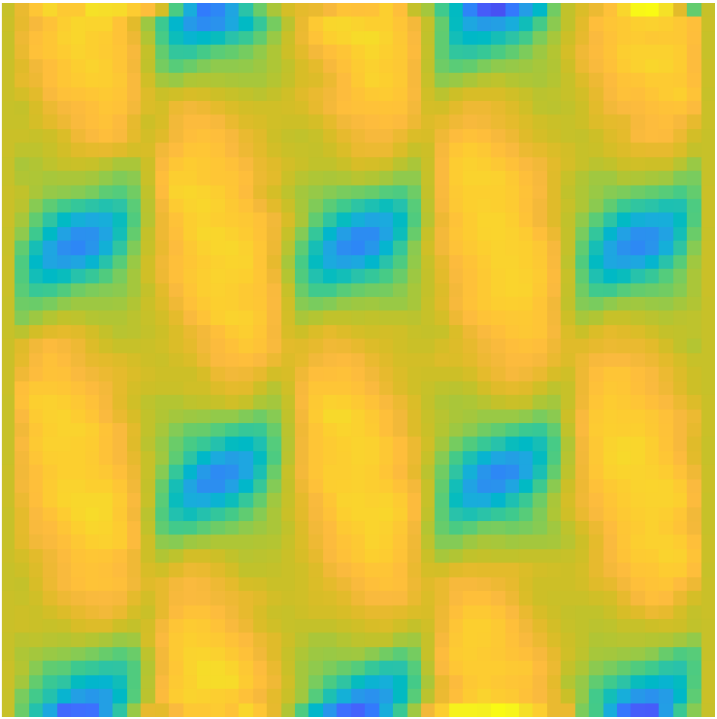} &
		  	\includegraphics[height=0.8in]{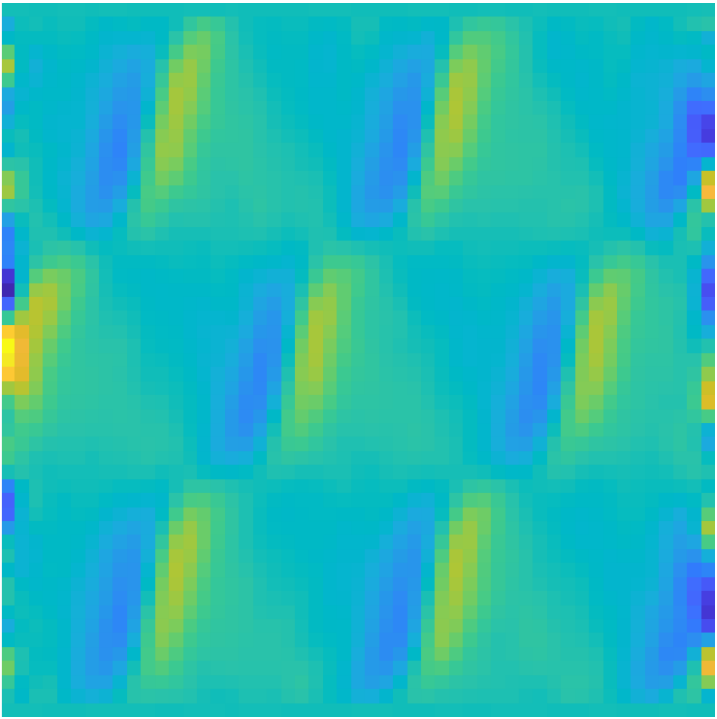} &
		  	\includegraphics[height=0.8in]{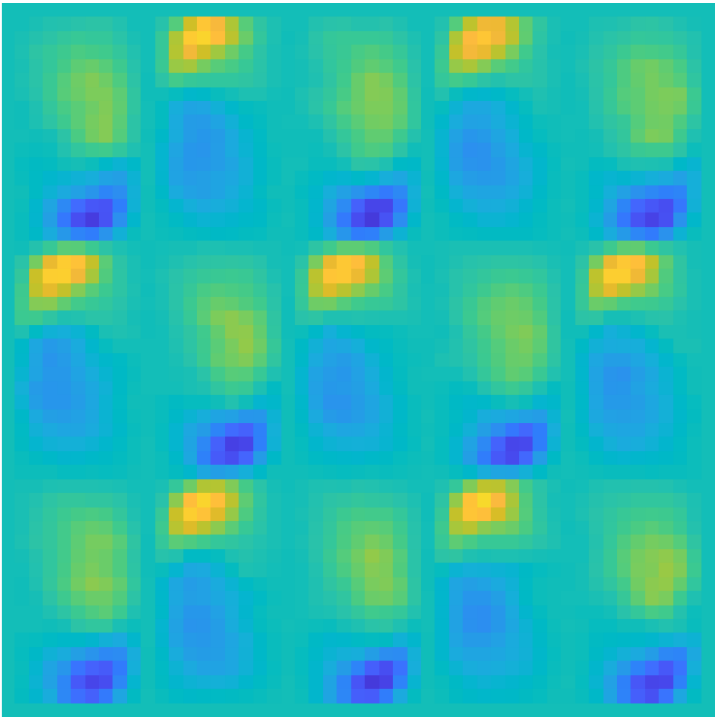}
		\end{tabular}
		\caption{SDD results with different smoothers. The first row  is the numerical solution of \eqref{PDEsmoother} at $t=0.15$ ($0\%$ noise) with the initial condition $u_0(x,y) =\sin(3\pi x)\sin(5\pi y)$ and its various partial derivatives. The second row shows the noisy data and its numerical derivatives when $5\%$ Gaussian noise is added to the clean data.  The bottom four rows are the SDD results at $t=0.15$ using MA, cubic spline (CS), diffusion (DF), and MLS in order.  While all methods recover $U$ (the first row), the dynamics of the derivatives, especially in the third and fourth rows, are best preserved by MLS.  }\label{fig_smoother}
	\end{figure}
\end{center}

In this paper, we use Moving Least Square (MLS) as the denoising in SDD.  To numerically justify this choice among Moving Average (MA)~\cite{tham1998dealing}, cubic spline interpolation~\cite{craven1978smoothing}, and diffusion smoothing~\cite{witkin1987scale},  we present the SDD results with these smoothers in Figure \ref{fig_smoother}. We first solve the PDE
\begin{equation}
u_t=-0.4uu_x-0.2 u u_y,~(x,y)\in[0,1]\times[0,1],~t\in[0,0.15]\;,
\label{PDEsmoother}
\end{equation}
with $\Delta t = 0.005$ and $\Delta x = \Delta y = 0.01$, where the initial condition is $u(x,y,0)=\sin(3\pi x)\sin(5\pi y)$. Then $5\%$ Gaussian noise is added to the numerical solution. 
Given the noisy data, we perform SDD denoising with different smoothers to obtain various partial derivatives. In MLS, we take the bandwidth $h =0.04$.
For MA, the window size for averaging is fixed to be $3$. For Cubic Spline, we use the MATLAB function $csaps$ with $p=0.5$. For the Diffusion denoising, we evolve the noisy surface following the heat equation $u_t=u_{xx}+u_{yy}$ with a time step size $(\Delta x)^2/4$ for $5$ iterations.
Figure \ref{fig_smoother} shows the SDD results of $u,u_x,u_{yy},uu_x$ at $t=0.15$ when different smoothers are used in SDD.
All of them recover $ U $ (the first row), while MLS preserves the underlying dynamics the best, i.e., the first and second-order derivatives.

\section{Conclusion}  \label{sec::conclusion}

This paper developed two robust methods for PDE identification from a single set of noisy data.  First, we proposed a Successively Denoised Differentiation (SDD) procedure to stabilize numerical differentiation, which significantly improves the accuracy in the computation of the feature matrix from noisy data.  We then proposed two new robust PDE identification algorithms called ST and SC. These algorithms utilize the Subspace Pursuit (SP) greedy algorithm to select a candidate set and then refine the results by time evolution or cross-validation.
We presented various numerical experiments to demonstrate the effectiveness of both methods.  SC is more computationally efficient, while ST performs better for PDEs with high order derivatives.  

\appendix

\section{Objectives of minimization}\label{Details_app}
 We discuss the error representation to compare different objective of PDE identification approaches.
We consider two ways to measure errors in PDE identification. The first one is the error between the identified numerical solution $\widehat{U}$ and the exact solution $u$, which is given by
$e(u):=\widehat{U}-u$.
The second error is $e(u_t):=D_t\widehat{U}-u_t$, which measures the difference between the numerical time derivative of $\widehat{U}$ and the ground truth $u_t$.  These two errors $e(u)$ and $e(u_t)$ are closely related, which relations are shown below (after Table \ref{tab:general}).

Many existing methods for the identification of PDEs or dynamical systems involve a minimization of $e(u)$ or $e(u_t)$. Consider the following decomposition of $e(u)$:
\begin{align}
e(u)=\underbrace{\widehat{U}-U}_{\text{Data fidelity}}+\underbrace{U-u}_{\text{Measurement error}},\label{eu}
\end{align}
where $U$ is the given data. In  (\ref{eu}),  the {\it Data fidelity} $\widehat{U}-U$ represents the accuracy of the identified PDE in comparison with the given data $U$. In literature, a class of {dynamic-fitting approaches} such as~\cite{bock1983recent,muller2004parameter,schmidt2009distilling,baake1992fitting} focus on controlling the data fidelity error in order to ensure if the numerical prediction is consistent with the evolution of the given data.  
The {\it Measurement error} $U-u$ comes from data acquisition where the given data are contaminated by noise. Denoising is an important step to reduce the measurement error.

\begin{table}[h]
	\centering
	\caption{Comparison of the objectives of PDE identification.  For parameter estimation problems (Type I), the  feature variables of the underlying PDEs are known.  For model identification problems (Type II),  such active set is unknown; hence sparsity is often imposed, or neural network is designed. 
	}~\label{tab:general}
	\begin{tabular}{ccc}
		\textbf{Problems}& \textbf{Objectives in Minimization}& \textbf{Methods}\\
		\hline
		\multirow{3}{*}{Type I} &
		Data fidelity&\cite{bock1983recent,baake1992fitting,timmer1998numerical,schmidt2009distilling,muller2004parameter,muller2002fitting}\\
		& Regression error& \cite{bock1981numerical,parlitz2000prediction,bar1999fitting,voss1999amplitude,liang2008parameter}\\
		& Regression error, Data fidelity& \cite{xun2013parameter}\\
		\hline
		\multirow{4}{*}{Type II} &
		Data fidelity&\cite{long2017pde}\\
		& Regression error& \cite{rudy2017data,schaeffer2017learning}\\
		& Regression error, Data fidelity& ST (Section~\ref{subsec::ST}),~\cite{kang2019ident}\\
		& Regression error, Coefficient error& SC (Section~\ref{subsec::SC})\\
		\hline
	\end{tabular}
\end{table}

The second error $e(u_t)$ can be  expressed as
\begin{align}
e(u_t)= \underbrace{D_t\widehat{U}-D_tU}_{\text{Response error}}+\underbrace{D_tU-F\widehat{\mathbf{c}}}_{\text{Regression error}}+ \underbrace{F(\widehat{\mathbf{c}}-\mathbf{c}_0)}_{\text{Coefficient error}}+\underbrace{(F-F_0)\mathbf{c}_0}_{\text{System error}}\;.\label{eut}
\end{align}
where $\widehat{\mathbf{c}}$ is the estimated coefficient. The first term  $D_t\widehat{U}-D_tU$ is called the {\it Response error}, which is the difference  between the numerical derivatives of the identified PDE and the given data. The  $L_2$ norm of the {\it Regression error} $D_tU-F\widehat{\mathbf{c}}$ is the most frequently used objective function in PDE identification for the regression-based methods \cite{bock1981numerical,parlitz2000prediction,bar1999fitting,voss1999amplitude,liang2008parameter}. In addition, one can introduce various types of regularization, such as the $L_1$ regularization~\cite{rudy2017data,schaeffer2017learning,kang2019ident} to induce sparsity.  The {\it  coefficient error}  $F(\widehat{\mathbf{c}}-\mathbf{c}_0)$ compares $\widehat{\bc}$ and  $\bc_0$. 
This term vanishes when $\widehat{\mathbf{c}}- \mathbf{c}_0$ lies in the null space of $F$, which can occur even when $\widehat{\mathbf{c}}\neq \mathbf{c}_0$. If the initial condition of the PDE is too simple, the null space of $F$ is very large, which makes the PDE identification problem ill-posed. See Equation (\ref{2DPDEexp1}) and \eqref{2DPDEexp1_eq} for an example. In order to guarantee a successful identification, the initial condition should have sufficient variations so that $F$ satisfies an incoherence or null space property \cite{donoho2001uncertainty}.
The final term $(F-F_0)\mathbf{c}_0 $ represents the {\it System error}, which is due to the numerical differentiation in the computation of $F$. Our SDD denoising technique can effectively reduce the system error.

We summarize the objectives considered by many existing methods in the literature in Table~\ref{tab:general} . These methods are categorized  according to which error term(s) that they aim at minimizing.
As for our proposed methods, ST minimizes the data fidelity, and SC focuses on the coefficient error and the regression error.

If the numerical scheme for the computation of $D_t\widehat{U}$ is consistent, then $\|e(u)\|_\infty\to 0$ and $\|e(u_t)\|_\infty\to 0$ are equivalent as $\Delta t,\Delta x\to 0$.
For $n=0,1,\dots,N$, we denote $e(u)^n$ and $e(u_t)^n$ as the values of $e(u)$ and $e(u_t)$ occurred at time $n\Delta t$, respectively. For $j=1,2,\dots,N$, we have
\begin{align*}
\frac{e(u)^{j}-e(u)^{j-1}}{\Delta t}&=\frac{\widehat{U}^{j}-\widehat{U}^{j-1}}{\Delta t}	-u_t^{j-1}+\mathbf{r}'=e(u_t)^{j-1}+\left(\frac{\widehat{U}^{j}-\widehat{U}^{j-1}}{\Delta t}-[D_t\widehat{U}]^{j-1}\right)+\mathbf{r}'\;,
\end{align*}
where $\|\mathbf{r}'\|_{\infty}=O(\Delta t)$. By induction, we obtain the following connection between $e(u)$ and $e(u_t)$:
\begin{align}
e(u)^{n}=e(u)^0+\sum_{j=0}^{n-1}e(u_t)^j\Delta t+\sum_{j=0}^{n-1}\left(\frac{\widehat{U}^{j+1}-\widehat{U}^j}{\Delta\ t}-[D_t\widehat{U}]^{j}\right)\Delta t+n\mathbf{r}\;,\label{th1eq1}
\end{align}
where the remainder $\|\mathbf{r}\|_{\infty}=O(\Delta t^2)$.  Equation~(\ref{th1eq1}) suggests that if the approximation $D_t\widehat{U}$ is consistent and $\|e(u)^0\|_\infty$ converges to $0$ as $\Delta x\to 0$ , $\|e(u)\|_\infty\to 0$  is equivalent to $\|e(u_t)\|_\infty\to 0$. Therefore, the PDE identification methods with the goal of having $\|e(u)\|_\infty$ or $\|e(u_t)\|_\infty$ approach to $0$ are equivalent.

It is often  practical to consider a grid-dependent $L_2$-norm of the errors, i.e., $\|\cdot\|_{2,\Delta}=\|\cdot\|_2\sqrt{\Delta x\Delta t}$ where $\|\cdot\|_2$ denotes the ordinary $L_2$ vector norm.  We provide an upper bound for  $\|e(u)\|_{2,\Delta}$.
\begin{theorem}\label{eu_2upper}
	Suppose $D_t\widehat{U}$ is computed using the forward difference. Then
	\begin{align}
	\|e(u)\|_{2,\Delta}^2 \leq X^dT^3\|e(u_t)\|^2_\infty +O(\|e(u_t)\|_\infty +\Delta t)+ O(\Delta t)\;.\label{eu_2}
	\end{align}
\end{theorem}
\begin{proof}
	Recall that $U\in\mathbb{R}^{M^dN}$ is the vectorization of the data. By the definition of the grid-dependent norm,  $\|U\|_{2,\Delta}^2=\Delta x^d\Delta t\|U\|_2^2=\frac{X^dT}{M^dN}\|U\|_2^2$. Using (\ref{th1eq1}), we have
	\begin{align*}
	\|e(u)\|_2^2&=\|e(u)^0\|_2^2+\sum_{n=1}^{N}\|e(u)^n\|_{2}^2\\
	&\leq \|e(u)^0\|^2_{2}+\sum_{n=1}^{N}\left(\sum_{j=0}^{n-1}\|e(u_t)^j\|_{2}\right)^2\Delta t^2+M^d\sum_{n=1}^{N}n^2O(\Delta t^4)+\\
	&\sum_{n=1}^{N}\|e(u)^0\|_{2}\sum_{j=0}^{n-1}\|e(u_t)^j\|_{2}\Delta t+M^{d/2}\sum_{n=1}^{N}\|e(u)^0\|_{2}nO(\Delta t^2)\\
	&+M^{d/2}\sum_{n=1}^{N}\sum_{j=0}^{n-1}\|e(u_t)^j\|_{2}nO(\Delta t^3)\\
	&\leq \|e(u)^0\|^2_{2}+\sum_{n=1}^{N}\left(\sum_{j=0}^{n-1}\|e(u_t)^j\|_{2}\right)^2\Delta t^2+M^dO(T^3\Delta t)\\
	&+\|e(u)^0\|_{2}\sum_{n=1}^{N}\sum_{j=0}^{n-1}\|e(u_t)^j\|_{2}\Delta t+M^{d/2}\|e(u)^0\|_{2}O(T^2)\\
	&+M^{d/2}\sum_{n=1}^{N}\sum_{j=0}^{n-1}\|e(u_t)^j\|_{2}nO(\Delta t^3)\;.
	\end{align*}
	Since $\|e(u_t)^j\|_{2}\leq M^{d/2}\|e(u_t)\|_\infty$, we can simplify the expression above as:
	\begin{align*}
	\|e(u)\|_{2}^2&\leq \|e(u)^0\|^2_{2}+M^{d}T^2N\|e(u_t)\|^2_\infty+M^dO(T^3\Delta t)\\
	&+TM^{d/2}N\|e(u)^0\|_{2}\|e(u_t)\|_\infty +M^{d/2}\|e(u)^0\|_{2}O(T^2)+M^{d}\|e(u_t)\|_\infty O(T^3)\;.
	\end{align*}
	Thus
	\begin{align*}
	\|e(u)\|_{2,\Delta}^2&=\Delta x^d\Delta t\|e(u)\|_{2}^2\\
	&\leq \Delta t\|e(u)^0\|^2_2+X^dT^3\|e(u_t)\|^2_\infty + O(X^dT^3\Delta t^2)\\
	&(\|e(u_t)\|_\infty +\Delta t)\|e(u)^0\|_2 O(T^2X^{d/2})+ X^d \|e(u_t)\|_\infty O(T^3\Delta t)\;.
	\end{align*}
\end{proof}
The upper bound expressed in~\eqref{eu_2} depends on several properties of the computational domain $\Omega$ and the sampling grid: the resolution $\Delta t$ and the domain size $X,T$. 
To derive useful information from Theorem~\ref{eu_2upper}, we assume that  $\|e(u_t)\|_\infty=O(\Delta t)$. This condition holds, for example, when we use first order forward difference and the underlying data is noiseless. 
\begin{corollary}\label{eu_0}
	When the time-space domain is fixed, i.e., $T>0$ and $X>0$, if $||e(u_t)||_\infty=O(\Delta t)$, we have
	\begin{align}
	\|e(u)\|_{2,\Delta}\to 0\;,\label{conv_eu_0}\quad\Delta t,\Delta x\to 0\;,
	\end{align}
\end{corollary}
This result suggests that, with the assumptions satisfied, increasing both the time and space resolutions is a sufficient condition for controlling $\|e(u)\|_{2,\Delta} \to 0$. The convergence of $\|e(u)\|_{2,\Delta}$ as $\Delta t,\Delta x\to0$ guarantees the success of the methods which minimize the data fidelity term, e.g., ST and IDENT in \cite{kang2019ident}.

\section{Proof of Proposition \ref{CEEbound}}\label{ProofCEEbound}
\begin{proof}
	\begin{align*}
	&[D_tU]^{\mathcal{T}_2}-[F]^{\mathcal{T}_2}_\calA\big([F]^{\mathcal{T}_1}_\calA\big)^\dagger [D_tU]^{\mathcal{T}_1}\\
	&=[D_tU]^{\mathcal{T}_2}-[u_t]^{\mathcal{T}_2}+[u_t]^{\mathcal{T}_2}-[F]^{\mathcal{T}_2}_\calA\big([F]^{\mathcal{T}_1}_\calA\big)^\dagger  [D_tU]^{\mathcal{T}_1}\\
	&=\underbrace{[D_tU]^{\mathcal{T}_2}-[u_t]^{\mathcal{T}_2}}_{E_1}+[u_t]^{\mathcal{T}_2}-[F]^{\mathcal{T}_2}_\calA\big([F]^{\mathcal{T}_1}_\calA\big)^\dagger [u_t]^{\mathcal{T}_1}\underbrace{-[F]^{\mathcal{T}_2}_\calA\big([F]^{\mathcal{T}_1}_\calA\big)^\dagger ([D_tU]^{\mathcal{T}_1}-[u_t]^{\mathcal{T}_1})}_{E_2}\\
	&=[u_t]^{\mathcal{T}_2}-([F_0]^{\mathcal{T}_2}_{\calA}+[F]^{\mathcal{T}_2}_{\calA}-[F_0]^{\mathcal{T}_2}_{\calA})\big([F]^{\mathcal{T}_1}_\calA\big)^\dagger [u_t]^{\mathcal{T}_1}+E_1+E_2\\
	&=[u_t]^{\mathcal{T}_2}-[F_0]^{\mathcal{T}_2}_{\calA}\big([F]^{\mathcal{T}_1}_\calA\big)^\dagger [u_t]^{\mathcal{T}_1}\underbrace{-([F]^{\mathcal{T}_2}_{\calA}-[F_0]^{\mathcal{T}_2}_{\calA})\big([F]^{\mathcal{T}_1}_\calA\big)^\dagger [u_t]^{\mathcal{T}_1}}_{E_3}+E_1+E_2\\
	&=\underbrace{[u_t]^{\mathcal{T}_2}-[F_0]^{\mathcal{T}_2}_{\calA_0}\big([F_0]^{\mathcal{T}_1}_{\calA_0}\big)^\dagger[u_t]^{\mathcal{T}_1}}_{=0}+\big([F_0]^{\mathcal{T}_2}_{\calA_0}\big([F_0]^{\mathcal{T}_1}_{\calA_0}\big)^\dagger-[F_0]^{\mathcal{T}_2}_{\calA}\big([F]^{\mathcal{T}_1}_\calA\big)^\dagger \big)[u_t]^{\mathcal{T}_1}\\
	&+E_1+E_2+E_3\\
	&=\big([F_0]^{\mathcal{T}_2}_{\calA_0}\big([F_0]^{\mathcal{T}_1}_{\calA_0}\big)^\dagger-[F_0]^{\mathcal{T}_2}_{\calA}\big([F]^{\mathcal{T}_1}_{\calA}\big)^\dagger\big)[u_t]^{\mathcal{T}_1}+E_1+E_2+E_3\\
	&=\big([F_0]^{\mathcal{T}_2}_{\calA_0}\big([F_0]^{\mathcal{T}_1}_{\calA_0}\big)^\dagger-[F_0]^{\mathcal{T}_2}_{\calA}\big([F_0]^{\mathcal{T}_1}_{\calA}\big)^\dagger\big)[u_t]^{\mathcal{T}_1}\\
	&\underbrace{-[F_0]^{\mathcal{T}_2}_{\calA}\big(\big([F]^{\mathcal{T}_1}_\calA\big)^\dagger -\big([F_0]^{\mathcal{T}_1}_{\calA}\big)^\dagger\big)[u_t]^{\mathcal{T}_1}}_{E_4}+E_1+E_2+E_3\\
	&=\big([F_0]^{\mathcal{T}_2}_{\calA_0}\big([F_0]^{\mathcal{T}_1}_{\calA_0}\big)^\dagger-[F_0]^{\mathcal{T}_2}_{\calA}\big([F_0]^{\mathcal{T}_1}_{\calA}\big)^\dagger\big)[u_t]^{\mathcal{T}_1}+E_1+E_2+E_3+E_4\;.
	\end{align*}
	Then we have:
	\begin{align*}
	&\mathrm{CEE} (\calA_k;\alpha,\mathcal{T}_1,\mathcal{T}_2)
	\leq \|\big([F_0]^{\mathcal{T}_2}_{\calA_0}\big([F_0]^{\mathcal{T}_1}_{\calA_0}\big)^\dagger-[F_0]^{\mathcal{T}_2}_{\calA}\big([F_0]^{\mathcal{T}_1}_{\calA}\big)^\dagger\big)[u_t]^{\mathcal{T}_1}\|_2
	\\&+\|[D_tU]^{\mathcal{T}_2}-[u_t]^{\mathcal{T}_2}\|_2
	+\|\big([F]^{\mathcal{T}_1}_\calA\big)^\dagger \|_2\,\big(\|[F]^{\mathcal{T}_2}_{\calA}\|_2\,\|[D_tU]^{\mathcal{T}_1}-[u_t]^{\mathcal{T}_1}\|_2\\
	&+\|[F]^{\mathcal{T}_2}_{\calA}-[F_0]^{\mathcal{T}_2}_{\calA}\|_2\,\|[u_t]^{\mathcal{T}_1}\|_2\big)
\\
&+\|[F_0]^{\mathcal{T}_2}_{\calA}\|_2\,\|\big([F]^{\mathcal{T}_1}_{\calA}\big)^\dagger\|_2\,\|\big([F_0]^{\mathcal{T}_1}_{\calA}\big)^\dagger\|_2\,\|[F]^{\mathcal{T}_1}_{\calA}-[F_0]^{\mathcal{T}_1}_{\calA}\|_2\,\|[u_t]^{\mathcal{T}_1}\|_2\;.
	\end{align*}
		In the last term on the right hand side of the inequality, we applied the norm bound in Theorem 4.1 of \cite{wedin1973perturbation}.  Then by setting 
	\begin{align}
	g(\mathcal{A};\alpha,\mathcal{T}_1,\mathcal{T}_2)&=\|[D_tU]^{\mathcal{T}_2}-[u_t]^{\mathcal{T}_2}\|_2
	+\|\big([F]^{\mathcal{T}_1}_\calA\big)^\dagger \|_2\,\big(\|[F]^{\mathcal{T}_2}_{\calA}\|_2\,\|[D_tU]^{\mathcal{T}_1}-[u_t]^{\mathcal{T}_1}\|_2\nonumber\\
	&+\|[F]^{\mathcal{T}_2}_{\calA}-[F_0]^{\mathcal{T}_2}_{\calA}\|_2\,\|[u_t]^{\mathcal{T}_1}\|_2\big)
	\nonumber\\
	&+\|[F_0]^{\mathcal{T}_2}_{\calA}\|_2\,\|\big([F]^{\mathcal{T}_1}_{\calA}\big)^\dagger\|_2\,\|\big([F_0]^{\mathcal{T}_1}_{\calA}\big)^\dagger\|_2\,\|[F]^{\mathcal{T}_1}_{\calA}-[F_0]^{\mathcal{T}_1}_{\calA}\|_2\,\|[u_t]^{\mathcal{T}_1}\|_2\label{eq_g_remainder}
	\end{align}
	we have proved the theorem.

\end{proof}

\bibliographystyle{siamplain}
\bibliography{references}
\end{document}

%% file: ex_shared.tex
\usepackage{lipsum}
\usepackage[noadjust]{cite}
\usepackage{amssymb,amsmath}
\usepackage{amsfonts}
\usepackage{multirow}
\usepackage{graphicx}
\usepackage{epstopdf}
\usepackage{algorithmic}
\usepackage{arydshln}
\usepackage{enumitem}
\setlist[enumerate]{leftmargin=.5in}
\setlist[itemize]{leftmargin=.5in}
\usepackage[algo2e,ruled,vlined]{algorithm2e}
\ifpdf
  \DeclareGraphicsExtensions{.eps,.pdf,.png,.jpg}
\else
  \DeclareGraphicsExtensions{.eps}
\fi

\newsiamremark{remark}{Remark}
\newsiamremark{hypothesis}{Hypothesis}
\crefname{hypothesis}{Hypothesis}{Hypotheses}
\newsiamthm{claim}{Claim}
\newcommand{\bc}{\mathbf{c}}
\newcommand{\bb}{\mathbf{b}}

\newcommand{\calT}{\mathcal{T}}
\newcommand{\calA}{\mathcal{A}}

\headers{Robust IDENT from A Single set of Noisy observation}{Y. He, S.H. Kang, W. Liao, H. Liu, Y. Liu}

\title{Robust Identification of differential equations by numerical techniques 
  from A Single Set of Noisy Observation}

\author{Yuchen He \thanks{Institute of Natural Sciences, Shanghai Jiao Tong University. Email: yuchenroy@sjtu.edu.cn}
\and Sung-Ha Kang \thanks{School of Mathematics, Georgia Institute of Technology. Email: kang@math.gatech.edu. Research is supported in part by Simons Foundation grant 282311 and 584960.}
\and Wenjing Liao \thanks{School of Mathematics, Georgia Institute of Technology. Email: wliao60@gatech.edu. Research is supported in part by NSF grant NSF-DMS 1818751 and NSF-DMS 2012652.}\and
Hao Liu \thanks{Department of Mathematics, Hong Kong Baptist University. Email: haoliu@ hkbu.edu.hk}
\and
Yingjie Liu \thanks{School of Mathematics, Georgia Institute of Technology. Email: yingjie@math.gatech.edu. Research is supported in part by NSF grants DMS-1522585 and DMS-CDS\&E-MSS-1622453.}}

\usepackage{amsopn}

\DeclareMathOperator{\sech}{sech}
